\documentclass[11pt,reqno]{amsart}
\usepackage[left=26mm,right=26mm,top=26mm,bottom=28mm]{geometry}
\usepackage{mathtools,amssymb,amsthm,mathrsfs,color}
\usepackage{txfonts}
\usepackage{microtype}
\usepackage{booktabs,tabularx,array}
\allowdisplaybreaks[2]
\usepackage[colorlinks,linkcolor=red,anchorcolor=green,citecolor=blue]{hyperref}

\usepackage[T1]{fontenc}
\usepackage[utf8]{inputenc}

\definecolor{bleu1}{RGB}{0,57,128}
\def\bleu1{\color{bleu1}}

\usepackage{etoolbox}
\patchcmd{\section}{\normalfont}{\normalfont \bleu1}{}{}
\patchcmd{\subsection}{\normalfont}{\normalfont \bleu1}{}{}
\patchcmd{\subsubsection}{\normalfont}{\normalfont \bleu1}{}{}

\newtheorem{The}{\bleu1 Theorem}[section]
\newtheorem{THEalpha}[The]{\bleu1 Theorem}
\newtheorem{Pro}[The]{\bleu1 Proposition}
\newtheorem{Lem}[The]{\bleu1 Lemma}

\theoremstyle{definition}
\newtheorem{defn}[The]{\bleu1 Definition}

\numberwithin{equation}{section}

\newcommand{\T}{\mathbb{T}}
\newcommand{\R}{\mathbb{R}}
\newcommand{\Z}{\mathbb{Z}}

\makeatletter
\@namedef{subjclassname@2020}{\textup{2020} Mathematics Subject Classification}
\makeatother

\newcommand{\Num}{\mathsf N}
\newcommand{\cF}{\mathcal F}
\newcommand{\cH}{\mathcal H}

\newcommand{\spec}{\operatorname{spec}}
\newcommand{\Tr}{\operatorname{Tr}}
\newcommand{\gap}{\operatorname{gap}}
\newcommand{\diam}{\operatorname{diam}}
\newcommand{\Span}{\operatorname{span}}

\title[Hubbard--Wigner states and thin-torus quantum Hall interactions]
{The Hubbard--Wigner model: from Aubry--Mather theory to Fractional Quantum Hall Effects}
\author{Tianqi Shi, Jinxin Xue}
\date{}
\address[Tianqi Shi]{School of Mathematics and Statistics, Northeast Normal University, Changchun 130024, China}
\email{tqshi@nenu.edu.cn}

\address[Jinxin Xue]{New Cornerstone Sciences Laboratory, Department of Mathematical Sciences, Tsinghua Unviersity, Room 305, Wenbei Building, Haidian, Beijing, 100084}
\email{jxue@tsinghua.edu.cn}
\keywords{Fractional quantum Hall effect, $V_1$ interaction, thin torus, density plateau, spectral gap, Chern class, Hubbard--Wigner model, Coulomb interaction}

\begin{document}
\begin{abstract}
The Hubbard--Wigner model describes binary lattice configurations
with convex repulsive pair interactions. We give a unified
classification of its periodic minimizers, minimizing invariant
measures, and Aubry sets through a positive block-count identity,
recovering explicit plateau widths and a complete devil's staircase
for the minimizing density. This variational structure parallels
Aubry--Mather theory and the authors' work on Aubry's completeness
conjecture, where the staircase is related to the sliding and
pinning mechanisms underlying the conducting and insulating
interpretations of particle chains. In the binary model,
completeness follows from discrete convexity and an exact sum
of plateau widths.

We then study how this classical structure persists for interactions
arising in the fractional quantum Hall effect on a thin torus.
At rational density $p/q$, the diagonal reference has a
$q$-dimensional lowest eigenspace, a spectral gap uniform in the
number of repeated periods, and an associated bundle of Chern
number $p$. For the full $V_1$ operator, we prove persistence of
a density-$1/3$ plateau and a positive spectral gap, uniformly in
system size at sufficiently small fixed transverse circumference;
the lowest-eigenspace bundle has rank three and Chern number one.
For the full periodic Coulomb operator, we establish multiplicity,
spectral-gap, and Chern-number conclusions at bounded denominators,
and quantitative bounds on plateau-endpoint displacement and loss
of width under explicit joint conditions on size and geometry.
\end{abstract}
\maketitle

\section{Introduction}\label{sec:introduction}
The Hubbard--Wigner model describes particles occupying sites of
an integer lattice and repelling one another. A configuration is a
binary sequence $\omega\in\Sigma_2:=\{0,1\}^{\mathbb Z}$: $\omega_i=1$
means that site $i$ is occupied. On a finite interval $\Lambda$,
its energy is
\begin{equation}\label{eq:intro-Hamiltonian}
 H_{\Lambda,c}(\omega)=
 \sum_{i,j\in\Lambda, i<j}U(j-i)\omega_i\omega_j
 -c\sum_{i\in\Lambda}\omega_i.
\end{equation}
The pair potential $U$ penalizes nearby particles, while increasing
$c$ favors configurations with more occupied sites. Hubbard studied
this one-dimensional model in connection with charge ordering in
organic conductors, and classified its periodic minimum-energy
configurations for convex repulsive interactions
\cite{Hubbard1978}. We use the name Hubbard--Wigner for this
classical occupation model; its energy contains no hopping term.

Repulsion favors evenly spaced particles, but the lattice permits
only integer separations. At density $1/3$, for example, the most
uniform arrangements are the three translates of
$\ldots100100\ldots$. At density $p/q$, the analogous arrangements
have $p$ occupied sites in each period of length $q$, distributed
as uniformly as possible. These are the balanced periodic words.
As $c$ varies, a fixed rational density can continue to minimize
the energy throughout a nontrivial interval. The resulting response
is a continuous nondecreasing function with horizontal intervals
at rational values: a \emph{devil's staircase}.
It is \emph{complete} if the interiors of these intervals cover
the transition interval up to a set of Lebesgue measure zero.
Thus completeness asks whether almost every parameter selects a
periodic arrangement, even though irrational minimizing densities
also occur. The staircase and its completeness for convex lattice
interactions were developed in
\cite{MR710121,Aubry_1983exact}; explicit plateau widths make this
model particularly accessible.

The authors' work on Aubry's completeness conjecture
\cite{ShiXue2026} places this question in a broader variational
setting. The underlying example is the Frenkel--Kontorova chain:
particles with real positions $x_j$ are connected by elastic
interactions and subjected to a periodic substrate potential.
A typical local energy is
\[
 \frac12(x_{j+1}-x_j)^2+V(x_j),\qquad V(x+1)=V(x).
\]
Its equilibrium equation defines an area-preserving twist map.
Aubry--Mather theory studies configurations that minimize the
energy of every finite segment with its endpoints fixed.
Their mean spacing is the rotation number. The theory produces
ordered minimizing sets at every rotation number: periodic orbits
at rational values, and invariant circles or ordered Cantor sets
at irrational values. In this way, a global family of minimizers
remains available even where invariant circles have broken down.

The physical distinction between an invariant circle and a broken
circle is central to Aubry's interpretation. A continuous family
of minimizing configurations permits the chain to change its
phase relative to the substrate without crossing an energy
barrier. This is the sliding regime. When such a family is broken,
intermediate phases can require positive energy, measured by a
Peierls barrier; the chain is pinned. For a chain carrying charge,
these are the mechanisms behind the conducting and insulating
interpretations, respectively
\cite[\S1.1.3]{ShiXue2026}. 

The common mathematical language is convex duality. If
$\beta_{\rm AM}(\rho)$ denotes minimum mean action at rotation
number $\rho$, then
\[
 \alpha_{\rm AM}(c)=\sup_\rho\{c\rho-\beta_{\rm AM}(\rho)\},
 \qquad \rho(c)=\alpha_{\rm AM}'(c).
\]
A corner of $\beta_{\rm AM}$ gives an interval of $c$ on which
$\rho(c)$ is constant. The authors prove that uniform hyperbolicity
of the union of the Mather sets over a parameter interval forces
$\alpha_{\rm AM}'$ to be purely singularly continuous
\cite[Theorem~1.5]{ShiXue2026}: its variation has no absolutely
continuous part. 

In the Hubbard--Wigner model, density takes the place of rotation
number, and the parameter $c$ in \eqref{eq:intro-Hamiltonian}
plays the role of the convex-dual parameter. Here a positive
block-count decomposition replaces the dynamical hyperbolicity
argument: it identifies the minimizing configurations and gives
plateau widths whose sum fills the transition interval.
Our first theorem gathers the periodic, invariant-measure, and
Aubry-set descriptions into this common variational framework.
We then use these classical minimizing words as reference states
for the quantum Hall operators introduced after the theorem.

\subsection{The classical variational structure}\label{sec:classical structure}
Let $U:\mathbb N\to[0,\infty)$ satisfy
\begin{enumerate}\renewcommand{\theenumi}{U\arabic{enumi}}
\item\label{item:U1} $\Delta U(r):=U(r)-U(r+1)\ge0$;
\item\label{item:U2} $\sum_{r\ge1}U(r)<\infty$;
\item\label{item:U3} $\Delta^2U(r):=U(r-1)-2U(r)+U(r+1)>0$ for $r\ge2$.
\end{enumerate}
We consider the functional \eqref{eq:intro-Hamiltonian}. Equip $\Sigma_2$ with the product topology and let $\sigma$ be the left shift, $(\sigma\omega)_j=\omega_{j+1}$.
Let $\mathcal M_\sigma$ be set of all $\sigma$-invariant probability measures on $\Sigma_2$
and define the density and mean energy by
\[
 \rho(\mu)=\int_{\Sigma_2}\omega_0\,d\mu,\qquad
 e_U(\mu)=\sum_{r\ge1}U(r)\int_{\Sigma_2}\omega_0\omega_r\,d\mu.
\]
The minimum mean energy and its convex dual are
\begin{equation}\label{eq:beta-def}
 \beta_U(\rho)=\min_{\substack{\mu\in\mathcal M_\sigma\\\rho(\mu)=\rho}}e_U(\mu),
 \qquad 0\le\rho\le1,
\end{equation}
\begin{equation}\label{eq:alpha}
 \alpha_U(c)=\max_{0\le\rho\le1}\{c\rho-\beta_U(\rho)\}.
\end{equation}
For an $L$-periodic word, write
$N(\omega)=\sum_{j=0}^{L-1}\omega_j$ and $\rho(\omega)=N(\omega)/L$.
Its mean energy, with the word extended periodically, is
\begin{equation}\label{eq:periodic-mean-energy}
 e_U(\omega)=\frac1L\sum_{j=0}^{L-1}\sum_{r\ge1}U(r)\omega_j\omega_{j+r}.
\end{equation}

\begin{defn}[Mechanical words and their orbit closures]\label{defn:mechanical}
For $\rho\in[0,1]$ and $\vartheta\in[0,1)$, the mechanical word is
\begin{equation}\label{eq:lower-mechanical}
 \omega_j^{\rho,\vartheta}
 =\lfloor(j+1)\rho+\vartheta\rfloor-\lfloor j\rho+\vartheta\rfloor.
\end{equation}
At a reduced rational density $p/q\in(0,1)$, write
\begin{equation}\label{eq:per}
 X_{p/q}^{\rm per}=\{\omega^{p/q,\vartheta}:0\le\vartheta<1\}.
\end{equation}
For irrational $\rho$, the Sturmian subshift is
\begin{equation}\label{eq:strumian}
 X_\rho^{\rm St}=\overline{\{\omega^{\rho,\vartheta}:0\le\vartheta<1\}}.
\end{equation}
\end{defn}
A binary sequence is \emph{balanced} if the numbers of ones in any
two blocks of the same length differ by at most one. Mechanical
words have this property by \eqref{eq:block-telescope}.
Lemma~\ref{lem:rational-balanced} shows that
$X_{p/q}^{\rm per}$ is a single $q$-point shift orbit.

For the pointwise classification, put

\begin{equation}\label{eq:Njn}
 S_n(j;\omega):=\sum_{r=0}^{n-1}\omega_{j+r},\qquad j\in\mathbb Z,\ n\ge1.
\end{equation}
For a reduced rational $p/q\in(0,1)$, define
\begin{align}
 \mathcal C_{p/q}^+
 &:=\{\omega:S_n(j;\omega)\in
       \{\lfloor np/q\rfloor,\lfloor np/q\rfloor+1\}\text{ for all }j\in\mathbb Z,n\ge1\},\label{eq:Cplus}\\
 \mathcal C_{p/q}^-
 &:=\{\omega:S_n(j;\omega)\in
       \{\lceil np/q\rceil-1,\lceil np/q\rceil\}\text{ for all }j\in\mathbb Z,n\ge1\}.\label{eq:Cminus}
\end{align}
\smallskip
\noindent\emph{Example: an irrational Sturmian configuration.}
Take $\rho=(\sqrt5-1)/2$ and $\vartheta=0$ in
\eqref{eq:lower-mechanical}. Its entries starting at $j=0$ are
\[
 0\,1\,0\,1\,1\,0\,1\,0\,1\,1\,0\,1\,1\,0\,1\,0\,1\,1\,0\,1\,0\,\cdots.
\]
One can construct this word by rotating a point on $\mathbb R/\mathbb Z$
through the angle $\rho$ at each step: the entry is one when
$\{j\rho+\vartheta\}\in[1-\rho,1)$ and zero otherwise.
Here $\{\cdot\}$ denotes fractional part. Telescoping the floor
differences shows that every length-$n$ block contains either
$\lfloor n\rho\rfloor$ or $\lceil n\rho\rceil$ ones. Thus its
density is $\rho$, and its particles are as evenly distributed as
the integer lattice permits. The sequence cannot be periodic,
since a periodic binary sequence has rational density. Varying
$\vartheta$ and taking limits gives $X_\rho^{\rm St}$; the closure
includes the limiting codings at the endpoints of the coding interval.
This is the irrational counterpart of the balanced periodic orbit.

\smallskip
\noindent\emph{Example: the two rational endpoint classes.}
The difference between \eqref{eq:Cplus} and \eqref{eq:Cminus}
occurs only when $q$ divides the block length $n$. At such a length,
the periodic configurations have exactly $np/q$ ones in every block;
$\mathcal C_{p/q}^+$ permits one additional one, whereas
$\mathcal C_{p/q}^-$ permits one fewer. At all other lengths,
both classes require the usual two balanced counts. In particular,
both contain $X_{p/q}^{\rm per}$, and their intersection is exactly
that periodic orbit by Lemma~\ref{lem:rational-balanced}.

For a concrete nonperiodic example at density $1/3$, specify a
configuration by listing its occupied sites. Consider
\[
\begin{array}{c|l}
 \mathcal C_{1/3}^+ & \ldots,-9,-6,-3,0,2,5,8,11,\ldots\\
 \mathcal C_{1/3}^- & \ldots,-9,-6,-3,0,4,7,10,13,\ldots .
\end{array}
\]
In the first row all successive spacings are three except for
the spacing two between $0$ and $2$; in the second they are three
except for the spacing four between $0$ and $4$. More generally,
the distance between particles $m$ steps apart is $3m$ or $3m-1$
in the first row, and $3m$ or $3m+1$ in the second. These distance
bounds give precisely the respective block constraints
\eqref{eq:Cplus} and \eqref{eq:Cminus}. For example, the block
$\{0,1,2\}$ in the first row contains two particles, while
$\{1,2,3\}$ in the second contains none, instead of the single
particle in every length-three block of a periodic minimizer.

Each example has density $1/3$ and agrees with a translate of
$\ldots100100\ldots$ on each tail, but the two tails have different
phases. Thus the signs $+$ and $-$ describe a local excess or
deficit in block counts, not a change of asymptotic density.
The theorem below places these additional configurations at the
upper and lower plateau endpoints, respectively; inside the plateau
only the periodic orbit remains in the Aubry set.

The Aubry set $\mathscr A(c)$ is defined by zero return Peierls barrier  for the normalized local energy
$\sum_{r\ge1}U(r)\omega_0\omega_r-c\omega_0+\alpha_U(c)$;
see \eqref{eq:Peierls-barrier}--\eqref{eq:natural-Aubry}.

\begin{THEalpha}[Classical classification and the complete staircase]
\label{intro:variational}\label{thm:exact-beta}\label{thm:periodic-minimizers}
Under the assumptions above:
\begin{enumerate}
\item The minimizing measure at each density is unique. At irrational
$\rho$ it is the unique invariant measure on the minimal Sturmian
subshift $X_\rho^{\rm St}$; at reduced $p/q\in(0,1)$ it is the uniform
measure on $X_{p/q}^{\rm per}$. For every $q\mid L$, these $q$ words
are exactly the $L$-periodic minimizers at density $p/q$.
The endpoint measures are the point masses on the constant words.
\item The function $\beta_U$ is strictly convex and differentiable
exactly at irrational interior densities. Its dual belongs to
$C^1(\mathbb R)$. The density $\alpha_U'(c)$ equals $p/q$ precisely
on $[c_{U,-}(p/q),c_{U,+}(p/q)]$, where
$c_{U,\pm}(\rho)=\beta'_{U,\pm}(\rho)$, and the plateau width is
\begin{equation}\label{intro:width}
 W_{U,p/q}:=q\sum_{k\ge1}k\Delta^2U(kq)>0.
\end{equation}
These plateau interiors exhaust
$(0,2\sum_{r\ge1}U(r))$ up to a Lebesgue null set.
Outside this interval $\alpha'_U$ equals zero or one, respectively.
\item If also $\sum_{r\ge1}rU(r)<\infty$ and
$\rho=\alpha'_U(c)\in(0,1)$, then $\mathscr A(c)=X_\rho^{\rm St}$
at irrational density. At rational density $p/q$, it equals
$X_{p/q}^{\rm per}$ in the plateau interior and
$\mathcal C_{p/q}^-$ or $\mathcal C_{p/q}^+$ at the lower or upper
endpoint. Every additional endpoint configuration has $q$-periodic
tails in both directions; there are no additional invariant minimizing
measures.
\end{enumerate}
\end{THEalpha}
This result is remarkably parallel to the classification of minimal configurations in Aubry-Mather theory. 
The periodic classification and width formula go back to
\cite{Hubbard1978,MR710121,Aubry_1983exact}; the formula also appears,
with a different pair-counting normalization, in
\cite[Eq.~(9)]{RotondoEtAl2016}. Here, we give a complete picture using the modern language of Aubry-Mather theory. Table \ref{tab:AM-comparison} records the analogy with the
classification in \cite[Theorem~2.2]{ShiXue2026}. Completeness here
follows by summing the explicit widths, independently of the
hyperbolicity assumption in the twist-map theorem.

\begin{table}[htbp]
\centering\small
\renewcommand{\arraystretch}{1.12}
\begin{tabularx}{\textwidth}{@{}p{.17\textwidth}XX@{}}
\toprule
 & Twist-map Aubry--Mather theory & Hubbard--Wigner model\\
\midrule
Parameter & Rotation number and dual cohomology class.
 & Density and its convex-dual parameter.\\
Irrational case & Invariant circle or Cantorus.
 & Unique minimizing measure on a Sturmian subshift.\\
Rational case & Periodic states and two connecting classes.
 & One periodic orbit; two endpoint classes with periodic tails.\\
complete devil's staircase & Aubry's completeness conjecture solved in \cite{ShiXue2026}
 & holds in general\\
\bottomrule
\end{tabularx}
\caption{The variational analogy; block balance replaces the ordering
of real configurations. The binary model has no invariant-circle alternative.}
\label{tab:AM-comparison}
\end{table}

\subsection{From the classical staircase to quantum Hall operators}\label{sec:classical-quantum}
Rotondo, Molinari, Ratti, and Gherardi
\cite{RotondoEtAl2016} connected the Hubbard--Wigner staircase
to the fractional quantum Hall effect by studying a torus with
one short circumference. We first explain the physical setting,
then formulate the finite-dimensional spectral problem used here.

Consider electrons confined to a plane in a perpendicular magnetic
field. An electric field can produce a current transverse to it;
the proportionality coefficient is the Hall conductivity
$\sigma_{xy}$. Experiments show intervals of magnetic field or
carrier density over which this coefficient remains at an integer
or fractional multiple of $e^2/h$, where $e$ is the magnitude of
the electron charge and $h$ is Planck's constant. These constant
values are Hall plateaux. At the same time, the longitudinal
response becomes very small. The fractional values reflect
collective effects of the interaction between electrons
\cite[\S1.3]{Tong2016}.

The magnetic one-particle operator has discrete energy levels,
called Landau levels. On a torus containing $L$ magnetic flux
quanta, the lowest level has dimension $L$. Restricting $N$
electrons to it leaves $L$ available orbitals, with filling
fraction $\nu=N/L$. The kinetic energy is then constant at
fixed $N$, so the interaction selects the lowest-energy states.
At $\nu=1/3$ there is one electron per three flux quanta; the
corresponding fractional Hall state has Hall conductivity of
magnitude $e^2/(3h)$ \cite[\S\S1.4, 3.1]{Tong2016}.

When one circumference is short, these orbitals can be indexed
along a one-dimensional cycle. Their occupation vectors are
binary words, and the diagonal part of the interaction assigns
an ordinary pair energy to each word. Rotondo et al. use the
thin-torus reduction of the Coulomb interaction to relate this
minimization to Hubbard's model and its staircase. The familiar
$1/3$ patterns are again the translates of $100100\ldots$
\cite{BergholtzKarlhede2008,RotondoEtAl2016}.
This is a bridge between two-dimensional magnetic interactions
and one-dimensional variational problems.

Our question is what remains when the non-diagonal terms,
which couple different occupation words, are retained. We study
three precise quantities: the intervals selecting a fixed degree,
the gap above the lowest eigenspace, and the Chern number of that
space as boundary conditions vary. The first measures persistence
of the density staircase. The latter two supply spectral and
geometric information that a classical occupation energy alone
does not contain. A plateau in the parameter $c$ will always mean
a density plateau in this paper; identifying it with an experimental
Hall-conductivity plateau requires further transport information.

We begin with the diagonal reference and its dependence on boundary
conditions. Section \ref{sec:classical structure} gives a single classical minimization problem.
We put its energies on the diagonal of an operator, and then realize
its occupation basis using one-particle eigenfunctions with variable
boundary phases. The matrix coefficients remain fixed; the changes
of frame determine a vector bundle over the parameter torus.
The two subsequent subsections treat the full $V_1$ and Coulomb
interactions in this framework.

We use the conventions of \cite[\S\S1.4, 3.1]{Tong2016},
with magnetic length and cyclotron energy equal to one. The
spatial torus has sides $a=2\pi/\kappa$ and $b=\kappa L$.
Increasing $\kappa$ shortens the first circumference; increasing
$L$ at fixed $\kappa$ enlarges the second. The $V_1$ estimates
hold uniformly in $L$ for sufficiently large fixed $\kappa$,
whereas the Coulomb results use explicit joint conditions on
$L$ and $\kappa$.

At fixed boundary parameters we use the graded exterior algebra
\begin{equation}\label{intro:occupation-space}
 \cF_L=\bigoplus_{N=0}^L\cH_{L,N},\qquad
 \cH_{L,N}=\bigwedge^N\mathbb C^L.
\end{equation}
Equip $\mathbb C^L$ with the orthonormal basis $e_0,\ldots,e_{L-1}$.
The index $N$ is the exterior degree. Throughout the paper, $L$
denotes the number of sites and $N$ the number of occupied sites.
For $0<N<L$, we write the density in lowest terms as
\[
 \rho=\frac NL=\frac pq,\qquad \gcd(p,q)=1,
 \qquad L=Mq,\quad N=Mp,\quad M=\gcd(N,L).
\]
Thus $p,q$ specify the reduced density, whereas $M$ specifies the
size of its finite realization. The minimizing configurations have
least period $q$, with $p$ occupied sites in each period, and repeat
this period $M$ times. An arbitrary configuration of the same density
need not have period $q$.
A binary
word $\omega=(\omega_0,\ldots,\omega_{L-1})$ with $\sum_j\omega_j=N$
corresponds to the unit occupation vector
\[
 \Psi_\omega=e_{j_1}\wedge\cdots\wedge e_{j_N},
 \qquad j_1<\cdots<j_N,\quad\omega_{j_a}=1.
\]
These vectors form an orthonormal basis of $\cH_{L,N}$.
Write $c_j^*$ for exterior multiplication by $e_j$, $c_j$ for its
adjoint, $n_j=c_j^*c_j$, and $\Num=\sum_jn_j$.
Thus $n_j\Psi_\omega=\omega_j\Psi_\omega$ and
$\Num|_{\cH_{L,N}}=NI$. Indices are read
modulo $L$, written $[j]$ when needed. The periodized classical operator is
\begin{align}
 K_L(d)&=\sum_{k\ge0}[U(d+kL)+U(L-d+kL)],\quad1\le d<L,
       \label{pc:periodization}\\
 D_L&=\sum_{0\le i<j<L}K_L(j-i)n_i n_j.
       \label{pc:diagonal-energy}
\end{align}
To compare this operator precisely with Section \ref{sec:classical structure}, set
$s_L=\sum_{k\ge1}U(kL)$. For an $L$-periodic word with $N$ ones,
grouping the interaction terms by their residues modulo $L$ gives
\[
 L e_U(\omega)
 =\langle\Psi_\omega,D_L\Psi_\omega\rangle+Ns_L.
\]
The last term accounts for interactions with periodic copies of the
same occupied site, which are excluded from $D_L$. It is constant
at fixed $N$, so the classical periodic minimizers are exactly the
occupation vectors spanning the lowest eigenspace of $D_L$.
When exterior degrees are compared, the corresponding identity is
\[
 L e_U(\omega)-cN
 =\langle\Psi_\omega,[D_L-(c-s_L)\Num]\Psi_\omega\rangle.
\]
Thus omitting the self-image term only translates the finite-period
convex-dual parameter; moreover $s_L\to0$ as $L\to\infty$.

For a self-adjoint operator $A$ commuting with $\Num$, define its
minimum eigenvalue and lowest eigenspace in degree $N$ by
\[
 E_A(N)=\min\operatorname{spec}(A|_{\cH_{L,N}}),\qquad
 \mathcal G_A(N)=\ker(A|_{\cH_{L,N}}-E_A(N)I).
\]
The spectral gap above this eigenspace is
\[
 \Delta_A(N)=\inf\{\lambda-E_A(N):
 \lambda\in\operatorname{spec}(A|_{\cH_{L,N}}),\ \lambda>E_A(N)\},
\]
with $\inf\varnothing=+\infty$.
A uniform gap in
Theorem~\ref{intro:periodic} means a lower bound
$\Delta_{D_L}(N)\ge\delta(U,q)>0$ for $L=Mq$ and $N=Mp$,
independent of the repetition number $M$.

For $0<N<L$, define the two secant-slope bounds
\[
 c_{A,-}(N)=\max_{0\le N'<N}\frac{E_A(N)-E_A(N')}{N-N'},\qquad
 c_{A,+}(N)=\min_{N<N'\le L}\frac{E_A(N')-E_A(N)}{N'-N}.
\]
Their open interval $I_A(N)=(c_{A,-}(N),c_{A,+}(N))$ is precisely
the set of parameters $c$ for which $N$ is the unique minimizer of
$N'\mapsto E_A(N')-cN'$. We call a nonempty such interval a
\emph{density plateau}. Thus a plateau compares different degrees,
whereas $\Delta_A(N)$ compares eigenvalues within one degree.
The associated convex function is
$\alpha_{A,L}(c)=L^{-1}\max_{0\le N\le L}\{cN-E_A(N)\}$.

For a periodic word, the balance condition of
Section~\ref{sec:classical structure} applies to its periodic extension,
including blocks crossing the period boundary or longer than $L$.
At density $\rho$, it is equivalent to
\begin{equation}\label{eq:periodic-balance}
 S_\ell(j;\omega)\in\{\lfloor\ell\rho\rfloor,\lceil\ell\rho\rceil\}
 \qquad(j\in\mathbb Z,\ \ell\ge1).
\end{equation}
The associated vector $\Psi_\omega$ is called a
\emph{balanced periodic occupation vector}.

Indeed, the average block count over the $L$ starting positions is
$\ell\rho$; integer block counts differ by at most one precisely
when they take the two indicated values. At reduced density $p/q$,
Lemma~\ref{lem:rational-balanced} shows that these words have least
period $q$ and form a single orbit of $q$ cyclic translates.
On $L=Mq$ sites they are the $M$-fold repetitions of these primitive
words. For $p/q=1/3$, they are the repetitions of $100$, $010$,
and $001$.

The magnetic boundary parameters are angles
$\theta=(\theta_x,\theta_y)$, each defined modulo $2\pi$. They form
the parameter torus
\begin{equation}\label{mg:twisttorus}
 \mathbb T^2_{\rm tw}=(\mathbb R/2\pi\mathbb Z)^2,
\end{equation}
oriented by $d\theta_x\wedge d\theta_y>0$. This is distinct from
the spatial torus on which the differential operator acts.
For the magnetic differential expression
\begin{equation}\label{mg:operator}
 H_0=\tfrac12[(-i\partial_x+y)^2+(-i\partial_y)^2],
\end{equation}
we impose
\begin{align}
 \psi(x+a,y)&=e^{-i\theta_x}\psi(x,y),\label{mg:bcx}\\
 \psi(x,y+b)&=e^{i\theta_y-ibx}\psi(x,y).\label{mg:bcy}
\end{align}
The differential expression stays fixed; the angles change the
boundary conditions and hence its domain.

Equivalently, one can fix the domain and vary the magnetic
potential by adding a closed one-form. In our convention the
covariant momentum is $-i\,d-A$, and the displayed operator
has local magnetic potential $A=-y\,dx$, with
$dA=dx\wedge dy$. Since its total flux is $ab=2\pi L\ne0$,
$A$ is not a globally defined one-form on the spatial torus:
it is a local expression for a connection on a line bundle.
The difference between two connections on that bundle, however,
is a globally defined one-form.

Indeed, write
\[
 \psi(x,y)=e^{-i\theta_x x/a+i\theta_y y/b}u(x,y).
\]
Then $u$ satisfies the boundary conditions with $\theta=0$,
and conjugating the operator gives the potential
\[
 A_\theta=A+A_{\mathrm{flat}},\qquad
 A_{\mathrm{flat}}=\frac{\theta_x}{a}\,dx
                  -\frac{\theta_y}{b}\,dy,
 \qquad dA_{\mathrm{flat}}=0.
\]
The signs follow from the boundary-phase convention above.
Thus $dA_\theta=dA$: varying $\theta$ changes the connection
but not the magnetic field. The added form determines a de Rham
\emph{cohomology} class. Its periods on the two oriented
coordinate loops are $\theta_x$ and $-\theta_y$; in other words,
its cohomology class pairs with the homology classes of these loops.
Exact changes are removed by single-valued real gauge phases,
and general $U(1)$ gauge transformations identify periods
that differ by multiples of $2\pi$. The resulting parameter
space is therefore
\[
 H^1(\mathbb T^2_{\rm space};\mathbb R)
 \big/2\pi H^1(\mathbb T^2_{\rm space};\mathbb Z)
 \simeq\mathbb T^2_{\rm tw}.
\]
Physically, these parameters describe Aharonov--Bohm fluxes
through the two noncontractible cycles, detected by the phase
acquired around a loop, while the local magnetic field is fixed.

Let $\mathcal E_L(\theta)$
be its $L$-dimensional lowest eigenspace. The corresponding
exterior power is $\bigwedge^N\mathcal E_L(\theta)$, locally
identified with $\cH_{L,N}$ by a choice of orthonormal eigenfunctions.
In these local frames, $D_L$ has the occupation-basis expression
\eqref{pc:diagonal-energy}.
More explicitly, writing $E_L(\omega)=
\sum_{i<j}K_L(j-i)\omega_i\omega_j$, the same word determines
\[
 \Psi_\omega(\theta)=\bigwedge_{j:\omega_j=1}^{\uparrow}e_j(\theta),
 \qquad
 D_L(\theta)\Psi_\omega(\theta)
 =E_L(\omega)\Psi_\omega(\theta).
\]
The coefficients $E_L(\omega)$ do not depend on $\theta$; the
notation $D_L(\theta)$ records the space on which the operator acts.
Consequently the spectrum is the same at every boundary parameter.
At $N=pL/q$, if $\omega^{(0)},\ldots,\omega^{(q-1)}$ are the
classical minimizers from Section \ref{sec:classical structure}, then
\[
 \mathcal G(\theta)=\operatorname{span}
 \{\Psi_{\omega^{(a)}}(\theta):0\le a<q\}.
\]
The magnetic changes of frame preserve this span and determine
its global gluing, as shown in Section \ref{mg:section}.
The two parameters therefore have different roles: $c$ selects a
minimizing exterior degree, while $\theta$ varies the boundary
conditions at fixed exterior degree and fixed operator coefficients.

As $\theta$ varies, an isolated lowest eigenspace
$\mathcal G(\theta)$ of constant dimension $r$ forms a Hermitian
vector bundle $\mathcal G\to\mathbb T^2_{\rm tw}$: its fiber over
$\theta$ is precisely $\mathcal G(\theta)$, and ``rank $r$'' means
that each fiber has dimension $r$. Local orthonormal bases of these
fibers are related on overlaps by unitary matrices. Even when the
operator matrix is independent of $\theta$ in local frames, these
changes of frame can make the bundle nontrivial.

The flat spatial connection supplies a point of
$H^1(\mathbb T^2_{\rm space};\mathbb R)/
2\pi H^1(\mathbb T^2_{\rm space};\mathbb Z)$; the topology
of this family of eigenspaces is instead measured by a class
in $H^2(\mathbb T^2_{\rm tw};\mathbb Z)$.
Its \emph{Chern number} is the integer
\[
 C(\mathcal G)=\frac{i}{2\pi}
 \int_{\mathbb T^2_{\rm tw}}\operatorname{tr}F_\nabla.
\]
Here $\nabla$ is any unitary connection on $\mathcal G$; in a local
orthonormal frame, $\nabla=d+\mathcal A$ and
$F_\nabla=d\mathcal A+\mathcal A\wedge\mathcal A$ is its curvature.
The integral is independent of this choice and equals the first
Chern number of the determinant line bundle
$\det\mathcal G=\bigwedge^r\mathcal G$.
It measures the topology of the entire lowest eigenspace family,
not of a selected eigenvector within a degenerate fiber
\cite{Husemoller1994}. Thus the last assertion of the next theorem
specifies both the fiber dimension $q$ and the integer $p$ associated
with its variation over the boundary-parameter torus.

\begin{THEalpha}[Lowest eigenspaces of the periodic diagonal operator]
\label{intro:periodic}\label{pc:main}\label{mg:main}
Fix a reduced density $0<p/q<1$. For every integer $M\ge1$, set
$L=Mq$ and $N=Mp$. The lowest eigenspace
$\mathcal G_{D_L}(N)$ is spanned by the $q$ exterior-product
vectors associated with the balanced periodic words. There exists
$\delta(U,q)>0$, independent of $M$, such that
\[
 \Delta_{D_L}(N)\ge\delta(U,q).
\]
Over $\mathbb T^2_{\rm tw}$, these eigenspaces form a Hermitian
vector bundle of rank $q$ and first Chern number $p$, both independent
of $M$.
\end{THEalpha}

The gap follows from an integer block defect at some length at most
$q$. The Chern calculation depends only on the primitive binary
word and the transition maps; it supplies the reference
bundle for both operator families.

\subsection{The full \texorpdfstring{$V_1$}{V1} interaction}\label{sec:V1}
The $V_1$ model is a positive two-body operator that isolates the
short-distance exclusion responsible for the Laughlin state at filling
$1/3$. To explain its meaning, first consider two spin-polarized
fermions in the lowest Landau level on the plane. Their relative
angular momentum is odd. A rotationally invariant pair interaction,
after projection to this level, is specified by its energies $v_m$
in the relative-angular-momentum subspaces. The pure $V_1$ model
keeps only the coefficient $v_1>0$:
\[
 H_{V_1}=v_1\sum_{i<j}P_{ij}^{(1)},
\]
where $P_{ij}^{(1)}$ projects the relative motion of the pair $i,j$
onto angular momentum one. Thus the subscript refers to relative
angular momentum, not to a nearest-neighbor interaction on a lattice
\cite{Haldane1983,TrugmanKivelson1985,Tong2016}.

Antisymmetry already forces a zero when two coordinates coincide.
The $V_1$ zero-energy condition removes the linear term of the
holomorphic relative-coordinate expansion, so that the zero has
order at least three. The cubic pair factors of the Laughlin state
satisfy this condition. On the torus, the corresponding periodic
model has three independent Laughlin zero modes at filling $1/3$
\cite{Laughlin1983,HaldaneRezayi1985}. It is therefore a parent
Hamiltonian for this state, meaning that the state is an exact
zero-energy minimizer. The model retains one specific component of
a repulsive interaction; it is not the full Coulomb operator. Its
value here is that exact kernel vectors can be combined with a
spectral estimate. Their explicit existence alone gives neither a
uniform gap nor a global density plateau.

Recent OpenAI Math Release preprints report progress on the full
$V_1$ spectral-gap problem beyond thin geometries. The
preprint establishes an inequality $H_Q^2\ge\gamma H_Q$ on the
entire lowest-Landau-level exterior algebra, with a fixed
$\gamma>1/25$ for all sufficiently large $Q$. These results concern expanding round spheres, and the first also
treats the planar limit. Our geometric regime is a torus with one
short circumference, and our argument starts from a diagonal
Hubbard--Wigner reference. It addresses the three-dimensional
lowest eigenspace, its Chern number, and the comparison of different
exterior degrees needed for an interval of constant minimizing density.
The spherical estimates do not directly give these torus statements,
and their numerical gap constants use a different normalization.
We cite the OpenAI works as recent preprints; the proofs below do
not use their estimates.

Our full torus normalization, for $L\in6\mathbb N$, is
\begin{equation}\label{ml:eq:QH}
 Q_R=\sum_{\substack{x>0\\x-R\in\Z}}
 x e^{-\kappa^2x^2}c_{[R+x]}c_{[R-x]},\qquad
 H_{L,\kappa}=\sum_{R\in\frac12\Z/L\Z}Q_R^*Q_R=\sum_{R\in\Z/L\Z}Q_R^*Q_R+\sum_{R\in(\Z+\frac{1}{2})/L\Z}Q_R^*Q_R.
\end{equation}

We take this standard torus expression from
\cite[\S1.1, Eq.~(1.2)]{WY}; see also \cite[\S3]{LNWY}.
Here $c_j$ removes a particle from the $j$th one-particle orbital.
Thus $c_{[R+x]}c_{[R-x]}$ removes a pair whose orbital indices have
midpoint $R$ and separation $2x$ before reduction modulo $L$.
The index $R$ labels the pair's center in orbital coordinates,
rather than its position as a point of the two-dimensional torus;
the conserved index sum is $2R$ modulo $L$. The operator $Q_R$
combines these pair-removal amplitudes with weights
$x e^{-\kappa^2x^2}$, including all periodic images. These are
amplitudes, so different separations can cancel each other.

The physical role of $Q_R$ is to test the pair component penalized
by the $V_1$ interaction at center index $R$. Mathematically,
$Q_R$ maps the $N$-particle space to the $(N-2)$-particle space, and
\[
 \langle\psi,H_{L,\kappa}\psi\rangle
 =\sum_R\|Q_R\psi\|^2,\qquad
 H_{L,\kappa}\psi=0\ \Longleftrightarrow\ Q_R\psi=0
 \text{ for every }R.
\]
Although $Q_R$ removes two particles, $Q_R^*Q_R$ preserves their
total number. Its expansion removes a pair at one separation and
recreates it at another with the same center index. This explains
both the diagonal pair energies and the non-diagonal pair
transitions in the decomposition below.

Its connection with the Hubbard--Wigner model becomes transparent
in the basis of binary occupation vectors. Write
\begin{equation}\label{ml:eq:Hubbard-perturbation}
 H_{L,\kappa}=D_{L,\kappa}+T_{L,\kappa},\qquad
 D_{L,\kappa}=\operatorname{diag}H_{L,\kappa}.
\end{equation}
The diagonal part is the finite cyclic pair interaction
\[
 D_{L,\kappa}=\sum_{0\le i<j<L}V_{L,\kappa}(j-i)n_i n_j.
\]
Its exact coefficient is
\begin{equation}\label{ml:eq:finite-V1-potential}
 V_{L,\kappa}(r)
 =f_{L,\kappa}(r/2)^2+f_{L,\kappa}((L-r)/2)^2,
 \qquad 1\le r<L,
\end{equation}
where
\begin{equation}\label{intro:periodic-amplitude}
 f_{L,\kappa}(x)
 =\sum_{j\in\mathbb Z}(x+jL)e^{-\kappa^2(x+jL)^2}.
\end{equation}
These formulas follow by collecting all terms in each $Q_R$ that
remove the same unordered pair, and then squaring their combined
coefficient; see Lemma~\ref{ml:lem:exact-diagonal-pair}.
In particular, $V_{L,\kappa}(r)=V_{L,\kappa}(L-r)$.
Thus $D_{L,\kappa}$ is an operator of the same pair-interaction form
as the classical operator in Section \ref{sec:classical-quantum}, with its finite cyclic
kernel specified directly by the full $V_1$ operator.
For the limiting variational problem, we use
\begin{equation}\label{ml:eq:classical-U}
 V_\kappa(r):=\lim_{\substack{L\to\infty\\6\mid L}}V_{L,\kappa}(r)
 =\frac{r^2}{4}e^{-\kappa^2r^2/2},\qquad r\ge1.
\end{equation}
This limiting potential satisfies the decreasing, summable, strictly
discretely convex hypotheses of Section \ref{sec:classical structure} for sufficiently large
$\kappa$. Proposition~\ref{ml:prop:limits} identifies the limiting
diagonal energy with its classical variational energy.

The non-diagonal part $T_{L,\kappa}$ replaces two occupied indices
by two others with the same sum modulo $L$. These are the off-diagonal matrix entries.
For $\kappa$ large, their Gaussian coefficients decay rapidly. In this precise sense,
the full $V_1$ operator is a small non-diagonal perturbation of a
Hubbard--Wigner-type diagonal operator. Smallness is measured by
the local interaction norm \eqref{ml:eq:interaction-norm}, rather
than by an operator norm on the whole system uniform in $L$.

For the spectral argument, the two leading diagonal coefficients
have different scales:
\[
 D_{L,\kappa}=a_\kappa A_L+\varepsilon_\kappa B_L
 +\text{smaller diagonal terms},\qquad
 a_\kappa=\tfrac14e^{-\kappa^2/2},\quad
 \varepsilon_\kappa=e^{-2\kappa^2},
\]
where $A_L=\sum_jn_jn_{j+1}$ and $B_L=\sum_jn_jn_{j+2}$.
We therefore implement the perturbation argument through the
positive comparison operator $K_{L,\kappa}$ of
\eqref{ml:eq:comparison}, which satisfies $K_{L,\kappa}\le H_{L,\kappa}$
and has the same kernel. After rescaling,
$\varepsilon_\kappa^{-1}K_{L,\kappa}=2A_L+B_L+\mathcal R_{L,\kappa}$,
and Lemma~\ref{ml:lem:tail} gives an interaction norm of order
$e^{-2\kappa^2}$ for the remainder. This is the fixed-reference
perturbation estimate used to prove the theorem below.

Let $\alpha_\kappa(c)=\lim_{L\to\infty,\,6\mid L}
\alpha_{H_{L,\kappa},L}(c)$, whose existence is part of the
following theorem. Let
$W^{\rm full}_{1/3}(\kappa)$ be the maximal interval length starting
at zero on which $\alpha_\kappa(c)=c/3$.
Denote the classical width associated with $V_\kappa$ in
\eqref{intro:width} by $W^{\rm cl}_{1/3}(\kappa)$.

\begin{THEalpha}[Degree selection and a spectral gap for the full $V_1$ operator]
\label{intro:quantum}\label{ml:thm:main}\label{ml:cor:selective}
\label{vb:full-chern}
For every $0<t<3$, there exist $\kappa_t>0$ and $g_t>0$
such that, for all $\kappa\ge\kappa_t$ and $L\in6\mathbb N$,
\[
 (0,t e^{-2\kappa^2})\subset I_{H_{L,\kappa}}(L/3),\qquad
 \Delta_{H_{L,\kappa}}(L/3)\ge g_t e^{-2\kappa^2}.
\]
The lowest eigenspaces of $H_{L,\kappa}|_{\mathcal H_{L,L/3}}$
form a rank-three Hermitian vector bundle over
$\mathbb T^2_{\rm tw}$ with first Chern number one. The gap bound
is uniform over this torus. The functions $\alpha_\kappa$ exist, and
\begin{equation}\label{ml:eq:relative-width}
 \liminf_{\kappa\to\infty}
 \frac{W^{\rm full}_{1/3}(\kappa)}{W^{\rm cl}_{1/3}(\kappa)}\ge1.
\end{equation}
For every $\kappa>0$, the convex function $\alpha_\kappa$ has no
nonempty open affine interval with slope strictly between zero and $1/3$.
\end{THEalpha}

\subsection{Coulomb stability and comparison of exterior degrees}\label{eq:coulomb stab}
The Coulomb interaction provides a second test of the classical
variational picture. It retains the long-range electrostatic
repulsion $v(r)=g/r$, $g>0$, between electrons confined to a plane.
After projection to the lowest Landau level, this gives an operator
on the same exterior algebra as the $V_1$ model. On the spatial
torus of Section~\ref{sec:classical-quantum}, define
\[
 k_{mn}=\left(\kappa m,\frac{2\pi n}{\kappa L}\right),\qquad
 \mathcal R_{mn}=\sum_{j=0}^{L-1}
 e^{2\pi i nj/L}c_{[j+m]}^*c_j,\qquad (m,n)\in\mathbb Z^2.
\]
We use the complete projected Coulomb operator of
\cite[Eq.~(2) and Supplementary Material, \S I]{RotondoEtAl2016},
written in the Fourier basis with magnetic length one:
\begin{equation}\label{co:representation}
 H^{\mathrm C}_{L,\kappa}
 =\sum_{(m,n)\ne(0,0)}w_{mn}
       (\mathcal R_{mn}^*\mathcal R_{mn}-\Num),
 \qquad
 w_{mn}=\frac{g}{2L|k_{mn}|}e^{-|k_{mn}|^2/2}.
\end{equation}
The vectors $k_{mn}$ are the Fourier frequencies of the spatial
torus. Up to a scalar phase and the Gaussian factor from projection,
$\mathcal R_{mn}$ represents a Fourier component of particle density:
$c_{[j+m]}^*c_j$ moves one occupied orbital index by $m$, and the
exponential supplies its Fourier phase. The factor $|k_{mn}|^{-1}$
comes from the two-dimensional Fourier transform
$\widehat v(k)=2\pi g/|k|$; the Gaussian comes from the lowest
Landau level. Subtracting $\Num$ removes the contribution of a
particle paired with itself. The resulting operator preserves
particle number and couples pairs of particles
\cite{BergholtzKarlhede2008,RotondoEtAl2016}.
Only the zero Fourier vector is omitted; all other periodic modes
are retained, and the series converges in operator norm at finite
$L$ and $\kappa$.

Its relation to the Hubbard--Wigner model is the exact decomposition
\[
 H^{\mathrm C}=D^{\mathrm C}+T^{\mathrm C},\qquad
 D^{\mathrm C}=\operatorname{diag}H^{\mathrm C}
 =\sum_{0\le i<j<L}W(j-i)n_i n_j.
\]
Here $W$ is the finite cyclic pair kernel in \eqref{co:W}, including
all terms of the exact periodic sum, while $T^{\mathrm C}$
has zero diagonal in the exterior-product basis.
The diagonal operator is a classical binary pair model of
Hubbard--Wigner type. We prove the discrete convexity needed to
select its balanced periodic minimizers directly on the finite
cycle, without imposing the infinite-lattice summability assumption
of Section \ref{sec:classical structure}. At fixed density $p/q$, these minimizers again span
a $q$-dimensional reference subspace. For $\kappa\ge3$, the estimate
\[
 \|T^{\mathrm C}\|\le g\kappa L e^{-2\kappa^2}
\]
quantifies the smallness of the non-diagonal part. Its dependence
on $L$ explains why the Coulomb results use a joint parameter
regime rather than the fixed-$\kappa$, uniform-in-$L$ regime of
our $V_1$ theorem.

There are two distinct objectives. At fixed exterior degree, we
ask whether the classical periodic space continues to an isolated
lowest eigenspace with the same dimension and Chern number. This
is the first assertion below, simultaneously for rational densities
with bounded denominator. To obtain density plateaux, we must
also compare the minimum eigenvalues at different exterior degrees.
Write $E_D^{\mathrm C}(N)=E_{D^{\mathrm C}}(N)$ and
$E^{\mathrm C}(N)=E_{H^{\mathrm C}}(N)$. The second assertion controls the displacement of plateau endpoints
and the possible loss of plateau width. Its proof uses a
second-order estimate for the energy difference, uniformly in $N$.

\begin{THEalpha}[Coulomb spectral estimates and minimizing intervals]
\label{intro:coulomb}
Let $\kappa\ge4$ and $L\ge6$.
\begin{enumerate}
\item If $Q\ge2$, $L\ge3Q$, and
$LQ^2e^{-2\kappa^2}\le(4\pi)^{-1}$, then for each reduced
$0<p/q<1$ with $q\le Q$ and $q\mid L$, the path
$D^{\mathrm C}+sT^{\mathrm C}$, $0\le s\le1$, restricted to $\mathcal H_{L,pL/q}$ has a lowest eigenspace of
dimension exactly $q$. The gap above it is at least
$g\kappa/(2\pi q^2)$. These eigenspaces form a vector bundle
over $\mathbb T^2_{\rm tw}$ with first Chern number $p$.
\item If $L^3e^{-2\kappa^2}\le(4\pi)^{-1}$, put
$\zeta=4\pi g\kappa L^3e^{-4\kappa^2}$.
Then, for every $0<N<L$,
\[
 |c_{H^{\mathrm C},\pm}(N)-c_{D^{\mathrm C},\pm}(N)|\le\zeta.
\]
In particular, if $|I_{D^{\mathrm C}}(N)|>2\zeta$, then
$I_{H^{\mathrm C}}(N)$ is nonempty and has length at least
$|I_{D^{\mathrm C}}(N)|-2\zeta$.
\end{enumerate}
The fixed-degree spectral estimates hold uniformly over
$\mathbb T^2_{\rm tw}$.
\end{THEalpha}

\subsection{Organization}
Sections~\ref{sec:physical-setting}--\ref{sec:minimizing-sets} establish
the classical formulas. Section~\ref{pc:section} obtains gaps from
integer block defects and records the finite-cycle identity used for
Coulomb. Section~\ref{mg:section} supplies the common magnetic geometry
and spectral-bundle continuation. Sections~\ref{ml:section} and~\ref{co:section} prove the operator estimates.
Appendix \ref{app:aubry-mather} gives the Aubry-set classification and restate the completeness of devil's staircase. Appendix \ref{co:analysis} provides the analytic
Coulomb estimates. Section \ref{sec:classical-quantum} introduces the exterior-algebra notation used throughout
the operator arguments. Throughout, $L$ is the dimension of the underlying vector space,
$N$ is the exterior degree, and $c$ is the convex-dual parameter; $H_{L,\kappa}$ denotes $V_1$ and
$H^{\mathrm C}_{L,\kappa}$ denotes Coulomb.

\section*{Acknowledgment}
J. X. is supported by the National Natural Science Foundation of China
(Grant No. 12271285, 11790273), and by the New Cornerstone Investigator program and Xiaomi
Foundation. Moreover, J. X. is a member of IAS supported by the Shiing-Shen Chern Member Fund this semester. T. S. is partly supported by the National Natural Science Foundation of China (Grant No.~12601348). 

ChatGPT was used in the sections related to Hall effect when working on detailed proofs. The paper was written by the authors. 

\section{The classical lattice model}
\label{sec:physical-setting}

Throughout Sections~\ref{sec:physical-setting}--\ref{pc:section},
we retain assumptions \eqref{item:U1}--\eqref{item:U3} and the
notation of Section~\ref{sec:classical structure}. Introduce the local energy
\begin{equation}\label{eq:periodic-energy}
 A_c(\omega)=\omega_0\sum_{r\ge1}U(r)\omega_r-c\omega_0.
\end{equation}
Then \eqref{eq:periodic-mean-energy} is the cyclic average of $A_0$.
For a finite interval $\Lambda\subset\mathbb Z$, we use the Hamiltonian $H_{\Lambda,c}$ defined in \eqref{eq:intro-Hamiltonian}.
For $\Lambda=\{0,\ldots,L-1\}$ and the periodic extension of $\omega$,
\[
 0\le L e_U(\omega)-H_{\Lambda,0}(\omega)
 \le b_L:=\sum_{r\ge1}\min(r,L)U(r),\qquad \frac{b_L}{L}\longrightarrow0.
\]
Dominated convergence gives the limit, so periodic and interval energies have the same limiting minimum per site.

\begin{Lem}\label{lem:U}
Under \eqref{item:U1}--\eqref{item:U3},
\[
 \lim_{r\to\infty} rU(r)=0,\qquad
 \sum_{r\ge1}(r+1)\Delta U(r)<\infty,\qquad
 \lim_{r\to\infty}r^2\Delta U(r)=0,\qquad
 \sum_{r\ge2}r^2\Delta^2U(r)<\infty.
\]
\end{Lem}
\begin{proof}
Monotonicity and summability give $(r/2)U(r)\le\sum_{j=\lceil r/2\rceil}^{r}U(j)\to0$ as $r\to\infty$. Note that $\Delta U(r)=U(r)-U(r+1)$ is nonnegative and decreasing by \eqref{item:U2}--\eqref{item:U3}, and summation by parts gives
\[
 \sum_{r=1}^m r\Delta U(r)=\sum_{r=1}^m U(r)-mU(m+1).
\]
Consequently $\sum_{r\geq 1} r \Delta U(r)=\sum_{r\geq 1} U(r)$ and $\sum_{r\geq 1}(r+1)\Delta U(r)=\sum U(r)+U(1)$. Also $(r^2/4)\Delta U(r)\le\sum_{j=\lceil r/2\rceil}^{r}j \Delta U(j)\to0$. Finally,
\[
 \sum_{j=2}^r j^2\Delta^2U(j)
 =-r^2\Delta U(r)-(2r+1)U(r)
 +2\sum_{j=1}^r U(j)+2U(1)-U(2).
\]
The right-hand side has a finite limit, proving the remaining claim.
\end{proof}

\section{Periodic minimizing configurations}
\label{sec:variational}\label{sec:exact-beta}

Convex repulsion favors the most uniform arrangement allowed by the lattice \cite{Hubbard1978}.
We express this principle through block counts: a minimizing word of density $\rho$ has either $\lfloor n\rho\rfloor$ or $\lceil n\rho\rceil$ occupied sites in every block of length $n$.
A positive decomposition of the energy will make this condition both necessary and sufficient.

\subsection{Block counts and the energy identity}
Use the block counts $S_n(j;\omega)$ from \eqref{eq:Njn}.
For an $L$-periodic word all indices are understood modulo $L$, including when $n>L$. Write $f(k)=\binom{k}{2}=k(k-1)/2$ for $k\in\mathbb Z$, and let $F$ be its piecewise affine interpolation:
\begin{equation}\label{eq:phi-formula}
 F(y)=ky-\frac{k(k+1)}2\quad(k\le y\le k+1,\ k\in\mathbb Z).
\end{equation}

The pair count $f(S_n)$ recovers the energy with weights given by the discrete curvature of $U$.

\begin{Lem}[Discrete convex decomposition]\label{lem:hinge-U}
For $m\ge1$,
\begin{equation}\label{eq:U-hinge}
 U(m)=\sum_{n=m+1}^{\infty}(n-m)\Delta^2U(n),
\end{equation}
and
\begin{equation}\label{eq:sumU-an}
 \sum_{m\ge1}U(m)=\frac12\sum_{n\ge2}n(n-1)\Delta^2U(n).
\end{equation}
\end{Lem}
\begin{proof}
Since $U(r)\to0$ and $\Delta U(r)\to0$,
\[
 U(m)=\sum_{r=m}^{\infty}\Delta U(r)
 =\sum_{r=m}^{\infty}\sum_{n=r+1}^{\infty}\Delta^2U(n)
 =\sum_{n=m+1}^{\infty}(n-m)\Delta^2U(n).
\]
Summing over $m$ and using nonnegativity proves \eqref{eq:sumU-an}.
\end{proof}

\begin{Lem}[Block energy identity]\label{lem:block-energy}
For every $L$-periodic word,
\begin{equation}\label{eq:periodic-block-energy}
 e_U(\omega)=\sum_{n\ge2}\Delta^2U(n)
 \frac1L\sum_{j=0}^{L-1}f\bigl(S_n(j;\omega)\bigr).
\end{equation}
The corresponding local identity is
\begin{equation}\label{eq:A0-Bn}
 A_0=\sum_{n\ge2}\Delta^2U(n)B_n,
\end{equation}
with uniform convergence, where
\begin{equation}\label{eq:Bn}
 B_n(\omega)=\omega_0\sum_{r=1}^{n-1}(n-r)\omega_r.
\end{equation}
\end{Lem}
\begin{proof}
Insert \eqref{eq:U-hinge} into $A_0$. The resulting series is uniformly convergent because $\|B_n\|_\infty\le n(n-1)/2$ and \eqref{eq:sumU-an} holds. Counting pairs in each block and averaging over its starting site gives
\[
 \frac1L\sum_{j=0}^{L-1}f\bigl(S_n(j;\omega)\bigr)
 =\sum_{r=1}^{n-1}(n-r)\frac1L\sum_{j=0}^{L-1}\omega_j\omega_{j+r}
 =\frac1L\sum_{j=0}^{L-1}B_n(\sigma^j\omega).
\]
This identity also holds for $n>L$, since positions in a block are counted before their indices are reduced modulo $L$. Averaging \eqref{eq:A0-Bn} now proves \eqref{eq:periodic-block-energy}.
\end{proof}

The cyclic mean of $S_n$ is $n\rho$. We therefore need Jensen's inequality with its integer equality cases.

\begin{Lem}[Jensen inequality for integer-valued random variables]\label{lem:integer-Jensen}
Let $X$ be an integrable integer-valued random variable with mean $y$. Then
\begin{equation}\label{eq:integer-Jensen}
 \mathbb E f(X)\ge F(y),
\end{equation}
where $\mathbb E$ is the expectation with respect to a probability on $\mathbb Z$. For $y\notin\mathbb Z$, equality holds precisely when
\begin{equation}\label{eq:integer-Jensen-equality-noninteger}
 X\in\{\lfloor y\rfloor,\lceil y\rceil\}\quad\hbox{almost surely}.
\end{equation}
For $y\in\mathbb Z$, equality holds precisely when $X=y$ almost surely.
\end{Lem}
\begin{proof}
Apply Jensen's inequality to $F$, which agrees with $f$ on $\mathbb Z$. At a noninteger mean, equality requires $X$ to lie in the unique affine interval $[\lfloor y\rfloor,\lceil y\rceil]$ containing that mean. For an integer mean $m$, use
\[
 f(k)-f(m)-\left(m-\tfrac12\right)(k-m)=\tfrac12(k-m)^2.
\]
Taking expectations shows that equality requires $X=m$ almost surely.
\end{proof}

Mechanical words attain these bounds at every length.

For the mechanical words of Definition~\ref{defn:mechanical},
the floor identity gives
\begin{equation}\label{eq:block-telescope}
 S_n(j;\omega^{\rho,\vartheta})
 =\lfloor(j+n)\rho+\vartheta\rfloor-\lfloor j\rho+\vartheta\rfloor
 \in\{\lfloor n\rho\rfloor,\lceil n\rho\rceil\}.
\end{equation}

\begin{Lem}[Rational block characterization]\label{lem:rational-balanced}
Let $0<p<q$ be coprime. A binary sequence satisfies
\begin{equation}\label{eq:rational-block-condition}
 S_n(j;\omega)\in\{\lfloor np/q\rfloor,\lceil np/q\rceil\}
 \quad(j\in\mathbb Z,\ n\ge1)
\end{equation}
if and only if it belongs to $X_{p/q}^{\rm per}$.
\end{Lem}
\begin{proof}
The condition at $n=q$ gives $S_q(j;\omega)=p$. Subtracting successive identities yields $\omega_{j+q}=\omega_j$. Put $P_0=0$, $P_j=\sum_{i=0}^{j-1}\omega_i$, and $T_j=P_j-jp/q$ for $0\le j\le q$. For $0\le r<s<q$, the block condition implies $T_s-T_r=S_{s-r}(r;\omega)-(s-r)p/q\in\{\lfloor(s-r)p/q\rfloor-(s-r)p/q,
   \lceil(s-r)p/q\rceil-(s-r)p/q\}$. Then $|T_s-T_r|<1$ by \eqref{eq:rational-block-condition}, since $(s-r)p/q$ is noninteger. Also $T_q=T_0=0$. Thus $\max_{1\leq j\leq q}T_j-\min_{1\leq j\leq q}T_j<1$. Choosing $\vartheta=\max_jT_j$ gives $ 0\le\vartheta<\min_{0\le j\le q}T_j+1\le1$ because of $T_0=0$ and
\[
 P_j\le jp/q+\vartheta<P_j+1\quad(0\le j\le q),
\]
so $\omega$ is mechanical. Conversely, mechanical words satisfy the block condition by \eqref{eq:block-telescope}. As $\vartheta$ varies, the word is constant on each interval $[r/q,(r+1)/q)$. Shifting the word adds $p/q$ to $\vartheta$ modulo one, so coprimality gives a single orbit of $q$ words. Any period $d$ must satisfy $dp/q\in\mathbb Z$; hence its least period is $q$.
\end{proof}

\subsection{Minimum energy and classification}
\begin{proof}[Proof of Theorem~\ref{intro:variational}: periodic minimizers and energy]
For an $L$-periodic word of density $\rho=N/L$, the cyclic mean of $S_n(j;\omega)$ is $n\rho$. Lemmas~\ref{lem:block-energy} and \ref{lem:integer-Jensen} imply
\[
 e_U(\omega)\ge\sum_{n\ge2}\Delta^2U(n)F(n\rho).
\]
The mechanical word of density $N/L$ is $L$-periodic, and \eqref{eq:block-telescope} makes every inequality an equality.
Conversely, equality forces every nonnegative Jensen defect to vanish, because $\Delta^2U(n)>0$.
The equality cases in Lemma~\ref{lem:integer-Jensen} give \eqref{eq:rational-block-condition} at every starting site; the length-one condition is automatic.
Lemma~\ref{lem:rational-balanced} then identifies the minimizers as the $q$ translates of the mechanical word at the reduced density $p/q=N/L$.

The resulting minimum has a continuous extension to $0\le\rho\le1$:
\begin{equation}\label{eq:beta-formula}
 \beta_U(\rho)=\sum_{n\ge2}\Delta^2U(n)F(n\rho)
 =\sum_{n\ge2}\Delta^2U(n)\sum_{j=1}^{n-1}(n\rho-j)_+.
\end{equation}
Indeed, $F(y)=\sum_{j\ge1}(y-j)_+$ for $y\ge0$, and all series converge uniformly by
\[
 \sum_{n\ge2}\Delta^2U(n)F(n)=\sum_{r\ge1}U(r)<\infty.
\]
Consequently, if $N_L/L\to\rho$, then
$\min_{N(\omega)=N_L}e_U(\omega)=\beta_U(N_L/L)\to\beta_U(\rho)$.
The boundary estimate in Section~\ref{sec:physical-setting} gives the same limit for the minimum interval energy per site.
\end{proof}

\section{Rational plateau endpoints and widths}
\label{sec:minimizing-sets}

The supporting slopes of $\beta_U$ determine the rational plateaux. Appendix~\ref{app:aubry-mather} proves completeness of the staircase.

The convex dual $\alpha_U$ is defined in \eqref{eq:alpha}.
For a convex function $g$, denote its \textit{subdifferential} by
\[
 D^-g(x):=\{a:g(y)\ge g(x)+a(y-x)\ \hbox{for every }y\}.
\]
At an interior point, $D^-g(x)=[g'_-(x),g'_+(x)]$. Therefore, a density $\rho$ minimizes $\beta_U(\rho)-c\rho$ exactly when $c\in D^-\beta_U(\rho)$. The rational plateau interiors are therefore the intervals stated in Theorem~\ref{intro:variational}.

\begin{Pro}[One-sided derivatives]\label{pro:beta-regularity}
The function $\beta_U$ is Lipschitz and strictly convex on $[0,1]$. For $0<\rho<1$,
\begin{align}\label{eq:beta-onesided}
 \beta'_{U,-}(\rho)&=\sum_{n\ge2}n\Delta^2U(n)
 (\lceil n\rho\rceil-1),\\
 \beta'_{U,+}(\rho)&=\sum_{n\ge2}n\Delta^2U(n)
 \lfloor n\rho\rfloor.\notag
\end{align}
At the endpoints, $\beta'_{U,+}(0)=0$ and $\beta'_{U,-}(1)=2\sum_{r\ge1}U(r)$.
\end{Pro}
\begin{proof}
Each summand in \eqref{eq:beta-formula} is convex. Its one-sided slopes have absolute value at most $n^2\Delta^2U(n)$ on $[0,1]$. Lemma~\ref{lem:U} therefore permits passage of one-sided difference quotients through the series and gives a common Lipschitz bound. The slopes of $F$ are $F'_-(y)=\lceil y\rceil-1$ and $F'_+(y)=\lfloor y\rfloor$, proving \eqref{eq:beta-onesided}. Every open subinterval of $(0,1)$ contains some $j/n$ with $1\le j<n$. The summand $\Delta^2U(n)F(n\rho)$ has a positive slope jump there, so the sum is affine on no nontrivial interval. Hence $\beta_U$ is strictly convex. Finally the same summable slope bound allows $\rho\to0^+$ and $\rho\to 1^-$; \eqref{eq:sumU-an} gives the stated endpoint values.
\end{proof}

\begin{proof}[Proof of Theorem~\ref{intro:variational}: rational plateau widths]
Strict convexity gives a unique minimizer of $\beta_U(\rho)-c\rho$. The subgradient criterion identifies its plateau at a reduced rational $p/q\in(0,1)$ with
\[
 \bigl(c_{U,-}(p/q),c_{U,+}(p/q)\bigr)
 =\bigl(\beta'_{U,-}(p/q),\beta'_{U,+}(p/q)\bigr).
\]
The slope of $F$ jumps by one at each integer. At $p/q$, precisely the terms with $q\mid n$ contribute to the jump in \eqref{eq:beta-onesided}, giving
\[
 W_{U,p/q}=c_{U,+}(p/q)-c_{U,-}(p/q)
 =\sum_{\substack{n\ge2\\q\mid n}}n\Delta^2U(n)
 =q\sum_{k\ge1}k\Delta^2U(kq)>0,
\]
Convergence follows from Lemma~\ref{lem:U}. The width depends only on $q$; the endpoints can depend on $p$.
\end{proof}

\section{Diagonal operators on exterior powers and a uniform spectral gap}
\label{pc:section}

We use the exterior spaces, occupation basis, and operators defined
in Section \ref{sec:classical-quantum}. In addition to the diagonal action of $n_j$ and
$\Num$, exterior multiplication and its adjoint satisfy
\[
 c_ic_j+c_jc_i=0,\qquad
 c_ic_j^*+c_j^*c_i=\delta_{ij}I,\qquad n_j=c_j^*c_j.
\]
The minimum energy $E_H(N)$, neutral gap, and density plateau
retain the meanings given there.

For the diagonal operator considered here, the block identity gives
a gap independent of the number of repeated periods: every
nonminimizing word has a positive integer defect at a block length
at most $q$.

\subsection{Periodic words and the correction for repeated indices}

Assume \textnormal{(U1)--(U3)}. For $L\ge2$, write $\Lambda_L=\mathbb Z/L\mathbb Z$, with representatives $0,\ldots,L-1$, and use the kernel $K_L$ and operator $D_L$ from \eqref{pc:periodization}--\eqref{pc:diagonal-energy}. The kernel $K_L$ is symmetric and its defining series converges.
We retain the self-image correction $s_L$ from
Section~\ref{sec:classical-quantum}. Cyclic words are extended
periodically, with mean energy \eqref{eq:periodic-mean-energy}.
Write $D_L(\omega)$ for its eigenvalue on $\Psi_\omega$;
$E_{D_L}(N)$ retains the minimum-energy notation of Section \ref{sec:classical-quantum}.

\begin{Pro}[Exact periodic energy]\label{pc:energy}
For every cyclic word with $N(\omega)=N$,
\begin{equation}\label{pc:cell-identity}
 D_L(\omega)=L e_U(\omega)-Ns_L.
\end{equation}
Consequently, for every integer $0\le N\le L$,
\begin{equation}\label{pc:exact-canonical}
 E_{D_L}(N)=L\beta_U(N/L)-Ns_L.
\end{equation}
\end{Pro}
\begin{proof}
In \eqref{eq:periodic-mean-energy}, terms with $r=kL$ contribute $Ns_L$ after multiplication by $L$. The other terms regroup into the unordered pairs of \eqref{pc:diagonal-energy}; nonnegativity justifies regrouping. This proves \eqref{pc:cell-identity}. Every periodic word of density $N/L$ has mean energy at least $\beta_U(N/L)$ by Theorem~\ref{thm:periodic-minimizers}. If $N/L=p/q$ is reduced, then $q\mid L$, so repetition of a balanced $q$-periodic word attains this bound on the ring. The constant words cover $N=0,L$, proving \eqref{pc:exact-canonical}.
\end{proof}

Fix a primitive binary word $w$ of length $q$ with $p$ occupied sites, and denote its $M$-fold repetition by $w^{(M)}$, where $L=qM$ and $N=pM$. Define
\begin{equation}\label{pc:root-orbit}
 \mathcal O_{w,M}=\{\sigma^a w^{(M)}:0\le a<q\},
\end{equation}
\begin{equation}\label{pc:root-space}
 \mathcal R_{w,M}=\Span\{\Psi_{\sigma^a w^{(M)}}:0\le a<q\},
\end{equation}
and let $P_{w,M}$ be the orthogonal projection onto this span.

\subsection{A uniform gap from integer block defects}

\begin{proof}[Proof of Theorem~\ref{intro:periodic}: minimizers and gap]
By \eqref{pc:cell-identity}--\eqref{pc:exact-canonical}, a cyclic word minimizes $D_L$ at $N=pM$ precisely when its periodic extension minimizes $e_U$ at $\rho=p/q$. Theorem~\ref{thm:periodic-minimizers} identifies these words with \eqref{pc:root-orbit}. Their $q$ distinct translates give $q$ orthonormal binary-word basis vectors, so the lowest eigenspace is exactly \eqref{pc:root-space}.

The block identity, with its self-image terms canceled, gives
\begin{align}
 D_L(\omega)-E_{D_L}(N)
 &=\sum_{\ell=2}^\infty\Delta^2U(\ell)\mathscr D_\ell(\omega),
 \label{pc:integer-defects}\\
 \mathscr D_\ell(\omega)
 &:=\sum_{j=0}^{L-1}\binom{S_\ell(j;\omega)}2-LF(\ell\rho).
 \label{pc:block-defect}
\end{align}
The series converges by Lemma~\ref{lem:U}, and Jensen's inequality gives $\mathscr D_\ell\ge0$.
Moreover, for $k=\lfloor\ell p/q\rfloor$,
\[
 LF(\ell p/q)=k\ell pM-qM\frac{k(k+1)}2\in\Z,
\]
so every $\mathscr D_\ell(\omega)$ is an integer.

If $\mathscr D_2(\omega)=\cdots=\mathscr D_q(\omega)=0$, equality at length $q$ gives $S_q(j;\omega)=p$ for every $j$. Subtracting adjacent counts yields $\omega_{[j+q]}=\omega_j$. Equality at lengths $r<q$ gives $S_r(j;\omega)\in\{\lfloor r\rho\rfloor,\lceil r\rho\rceil\}$, with $r=1$ automatic. Writing $\ell=aq+r$, $0\le r<q$, shows that every block has one of the two balanced counts. Thus every $\mathscr D_\ell(\omega)$ vanishes and $\omega\in\mathcal O_{w,M}$. Outside that orbit, some $\mathscr D_\ell(\omega)\ge1$ with $2\le\ell\le q$. Consequently every binary-word basis vector outside $\mathcal R_{w,M}$ has energy at least $\delta(U,q)$ above the minimum.
Since $D_L$ is diagonal, this is the operator inequality
\begin{equation}\label{pc:uniform-gap}
 D_L|_{\mathcal H_{L,N}}-E_{D_L}(N)I
 \ge\delta(U,q)(I-P_{w,M}),
 \qquad \delta(U,q):=\min_{2\le\ell\le q}\Delta^2U(\ell)>0.
\end{equation}
The lower bound is independent of $M$.
\end{proof}

\subsection{Finite cyclic kernels without a summability assumption}

The exact diagonal kernels below need not arise from a summable infinite-lattice potential.
The following finite block decomposition suffices for their minimization; the possibly negative constant $v_*$ must be retained in comparisons between exterior degrees.

\begin{Lem}[Minimum energy for a finite cyclic kernel]\label{pc:cyclic-kernel}
Let $L\ge2$, $J_L=\lfloor L/2\rfloor$, $v_*\in\R$, and $a_\ell\ge0$ for $2\le\ell\le J_L$. Suppose
\begin{equation}\label{pc:cyclic-hinge}
 K(r)=v_*+\sum_{\ell=2}^{J_L} a_\ell
 (\ell-\min\{r,L-r\})_+,
 \qquad 1\le r<L.
\end{equation}
For $D_K=\sum_{0\le i<j<L}K(j-i)n_in_j$,
\begin{equation}\label{pc:cyclic-minimum}
 E_{D_K}(N)=v_*\binom N2+L\sum_{\ell=2}^{J_L} a_\ell F(\ell N/L),
 \qquad 0\le N\le L.
\end{equation}
Here $F$ is the interpolation in \eqref{eq:phi-formula}; the sums are empty when $J_L=1$, and $v_*$ may be negative.
\end{Lem}
\begin{proof}
For a cyclic word $\omega$ with $N(\omega)=N$, use the block counts in \eqref{eq:Njn}. A pair at cyclic distance $d$ belongs to $(\ell-d)_+$ windows of length $\ell\le L/2$. Counting pairs gives
\begin{equation}\label{pc:cyclic-block-identity}
 D_K(\omega)=v_*\binom N2+\sum_{\ell=2}^{J_L} a_\ell
 \sum_{j=0}^{L-1}\binom{S_\ell(j;\omega)}2.
\end{equation}
The average of $S_\ell(j;\omega)$ is $\ell N/L$, so Lemma~\ref{lem:integer-Jensen} bounds each inner sum below by $LF(\ell N/L)$. Formula~\eqref{eq:block-telescope}, at density $N/L$, attains all these bounds on an $L$-periodic mechanical word; the constant words cover $N=0,L$. This proves \eqref{pc:cyclic-minimum} using only finite sums.
\end{proof}

\section{Magnetic boundary parameters and Chern classes}
\label{mg:section}

We use the cohomology parameter space introduced in
Section~\ref{sec:classical-quantum}.
We now study the lowest eigenspaces as a family over this quotient.
There are two distinct bundles. The magnetic line bundle
$\mathscr L_\theta$ is over the \emph{spatial} torus; adding
$A_{\mathrm{flat}}$ changes its connection but leaves its curvature
and degree $L$ unchanged. The bundle $\mathcal E_L$ of lowest
eigenspaces is over the \emph{parameter} torus, with fiber
$\mathcal E_L(\theta)\subset L^2(\mathscr L_\theta)$.
Its first Chern class belongs to
$H^2(\mathbb T^2_{\rm tw};\mathbb Z)$, not to the cohomology
of the spatial torus.

An integral-period change of $A_{\mathrm{flat}}$ represents the
same parameter, but the corresponding identification can act
nontrivially on an eigenfunction basis. These identifications
are the transition maps of $\mathcal E_L$. Thus closed spatial
one-forms supply the parameters, while the failure of local
spectral bases to glue to a global basis is measured by a Chern
class on their parameter space. We compute this class from
explicit transition matrices and then show that it is preserved
under gapped continuation.

\subsection{The magnetic operator and its lowest eigenspaces}

Fix an integer $L\geq1$ and $a,b>0$ with $ab=2\pi L$. The rectangle $[0,a]\times[0,b]$ represents the spatial torus. Use the magnetic expression \eqref{mg:operator} with boundary
conditions \eqref{mg:bcx}--\eqref{mg:bcy}.
The condition $ab=2\pi L$ makes the two boundary identifications
compatible. They define a Hermitian line bundle $\mathscr L_\theta$
over the spatial torus. The differential expression $H_0$ is fixed;
$H_{0,\theta}$ denotes its self-adjoint realization on
$L^2(\mathbb T^2_{\rm space},\mathscr L_\theta)$, obtained from the
quadratic form
\[
 q_\theta(\psi)=\frac12\int_{[0,a]\times[0,b]}
 \left(|(-i\partial_x+y)\psi|^2+|\partial_y\psi|^2\right)\,dx\,dy
\]
on $H^1$ sections satisfying the stated identifications.

On the parameter torus \eqref{mg:twisttorus}, let $Q=[0,2\pi]^2$
be the fundamental square. With the orientation fixed in
Section~\ref{sec:classical-quantum}, its positive generator is
\begin{equation}\label{mg:eta}
 \mathfrak u=\left[\frac{d\theta_x\wedge d\theta_y}{(2\pi)^2}\right]
 \in H^2(\T^2_{\rm tw};\Z),
 \qquad \langle\mathfrak u,[\T^2_{\rm tw}]\rangle=1.
\end{equation}

\begin{Lem}[Explicit frame of the lowest eigenspaces]\label{mg:orbitals}
The lowest eigenvalue of $H_{0,\theta}$ is $1/2$, of multiplicity $L$. An orthonormal basis of this eigenspace is
\begin{equation}\label{mg:gaussians}
 \phi_m^\theta(x,y)=\frac{\pi^{-1/4}}{\sqrt a}
 \sum_{n\in\Z}e^{-in\theta_y}
 e^{ik_{m+nL}(\theta_x)x}
 e^{-(y+k_{m+nL}(\theta_x))^2/2},
 \qquad 0\leq m<L,
\end{equation}
where $k_j(\theta_x)=(2\pi j-\theta_x)/a$. These vectors form a smooth frame over the square $Q$.
\end{Lem}
\begin{proof}
Gaussian decay gives smooth convergence. The $x$ boundary condition is immediate; $k_{j+L}=k_j+b$ and reindexing give the $y$ condition. Each summand is a lowest eigenfunction of the harmonic oscillator of energy $1/2$. Fourier orthogonality in $x$ makes distinct residues orthogonal, and
\[
 \|\phi_m^\theta\|^2
 =\pi^{-1/2}\sum_{n\in\Z}\int_0^b e^{-(y+k_m+nb)^2}\,dy=1,
\]
since the shifted intervals partition $\R$.

To prove completeness, we describe the Fourier decomposition explicitly.
For fixed $y$, the first boundary condition makes
$e^{i\theta_x x/a}\psi(x,y)$ periodic in $x$. Hence
\[
 \psi(x,y)=\frac1{\sqrt a}\sum_{j\in\mathbb Z}
 F_j(y)e^{ik_jx},\qquad
 F_j(y)=\frac1{\sqrt a}\int_0^a\psi(x,y)e^{-ik_jx}\,dx.
\]
The series converges in $L^2$, and Parseval's identity gives
\[
 \|\psi\|^2=\sum_{j\in\mathbb Z}\int_0^b|F_j(y)|^2\,dy.
\]
Here and below $k_j=k_j(\theta_x)$. Since $ab=2\pi L$,
we have $k_{j+L}=k_j+b$. Substitution into the second boundary
condition yields
\[
 \sum_jF_j(y+b)e^{ik_jx}
 =e^{i\theta_y}\sum_jF_j(y)e^{i(k_j-b)x}
 =e^{i\theta_y}\sum_jF_{j+L}(y)e^{ik_jx}.
\]
Uniqueness of Fourier coefficients therefore implies
\[
 F_j(y+b)=e^{i\theta_y}F_{j+L}(y).
\]
For form-domain sections this identity is used in the trace sense
at $y=0,b$; the identity for other $y$ refers to their extension
by the boundary identifications.

For each $m\in\{0,\ldots,L-1\}$, the intervals
\[
 I_{m,n}=[k_m+nb,k_m+(n+1)b],\qquad n\in\mathbb Z,
\]
cover $\mathbb R$ with disjoint interiors. Define on each interval
\[
 g_m(t)=e^{in\theta_y}F_{m+nL}(t-k_m-nb),
 \qquad t\in I_{m,n}.
\]
The traces at adjacent endpoints agree, because
\[
 e^{in\theta_y}F_{m+nL}(b)
 =e^{i(n+1)\theta_y}F_{m+(n+1)L}(0).
\]
Parseval's identity now becomes
\[
 \|\psi\|^2=\sum_{m=0}^{L-1}\int_{\mathbb R}|g_m(t)|^2\,dt.
\]
Conversely, given $g_m\in L^2(\mathbb R)$, set
$F_{m+nL}(y)=e^{-in\theta_y}g_m(y+k_m+nb)$ on $[0,b]$
and use the Fourier series above. This is the inverse map, so the
construction is a unitary identification with
$\bigoplus_{m=0}^{L-1}L^2(\mathbb R)$.

On the form domain, Fourier orthogonality also gives
\begin{align*}
 q_\theta(\psi)
 &=\frac12\sum_{j\in\mathbb Z}\int_0^b
 \bigl(|F_j'(y)|^2+(y+k_j)^2|F_j(y)|^2\bigr)\,dy\\
 &=\frac12\sum_{m=0}^{L-1}\int_{\mathbb R}
 \bigl(|g_m'(t)|^2+t^2|g_m(t)|^2\bigr)\,dt.
\end{align*}
The matching traces ensure that the joined functions lie in
$H^1(\mathbb R)$; conversely, functions in $H^1(\mathbb R)$ with
$tg_m\in L^2(\mathbb R)$ reconstruct precisely the form domain.
Thus the associated self-adjoint operator is unitarily equivalent
to $L$ copies of $\frac12(-\partial_t^2+t^2)$.
For this one-dimensional form,
\[
 \frac12\int_{\mathbb R}(|g'|^2+t^2|g|^2)\,dt
 =\frac12\|g\|^2+\frac12\|g'+tg\|^2,
\]
first by integration by parts on smooth compactly supported
functions and then by form closure. Equality holds precisely when
$g'+tg=0$, that is, when $g$ is a multiple of $e^{-t^2/2}$.
Choosing the normalized Gaussian $\pi^{-1/4}e^{-t^2/2}$ in the
$m$th summand and zero in the others reconstructs exactly
$\phi_m^\theta$ in \eqref{mg:gaussians}. This proves the lower
bound, multiplicity $L$, and completeness of the stated basis.
\end{proof}

The spaces $\mathcal E_L(\theta)=\ker(H_{0,\theta}-\tfrac12)$
therefore give the rank-$L$ bundle over $\T^2_{\rm tw}$ introduced
in Section \ref{sec:classical-quantum}. On the covering space $\mathbb R^2$ of parameters,
one may choose representatives $A_{\mathrm{flat}}$ and smooth
local bases. To descend to the quotient, one must identify the
bases at $\theta$ and $\theta+2\pi k$, $k\in\mathbb Z^2$.
The following calculation makes these identifications explicit.

\subsection{Transition matrices of the lowest-eigenspace bundle}

In the coordinate basis of Section~\ref{sec:classical-quantum}, set
\begin{equation}\label{mg:clockshift}
 S_Le_m=e_{m+1},\qquad
 C_Le_m=e^{2\pi im/L}e_m,
 \qquad S_LC_L=e^{-2\pi i/L}C_LS_L,
\end{equation}
where indices are taken modulo $L$. A pair of matrices $u_x(\theta_y),u_y(\theta_x)$ describes the bundle by the coordinate identifications
\begin{align}
 (0,\theta_y,v)&\sim(2\pi,\theta_y,u_x(\theta_y)v),\label{mg:quotientx}\\
 (\theta_x,0,v)&\sim(\theta_x,2\pi,u_y(\theta_x)v).
 \label{mg:quotienty}
\end{align}
Consistency at the corner requires $u_x(2\pi)u_y(0)=u_y(2\pi)u_x(0)$.

\begin{Pro}[Transition maps obtained from the magnetic eigenfunctions]
\label{mg:orbitalclutch}
There is an eigenfunction frame over $Q$, obtained from \eqref{mg:gaussians} by multiplying each basis eigenfunction by a scalar of modulus one, in which the exact coordinate identifications are
\begin{equation}\label{mg:lllclutch}
 u_x(\theta_y)=e^{i\theta_y/L}S_L,
 \qquad u_y(\theta_x)=C_L.
\end{equation}
\end{Pro}
\begin{proof}
Rescale the basis in \eqref{mg:gaussians} by setting
\[
 \widetilde\phi_m^\theta
 =e^{-im\theta_y/L}\phi_m^\theta,\qquad 0\le m<L.
\]
This choice distributes the boundary phase equally among the $L$
cyclic indices. Indeed, reindexing the series in
\eqref{mg:gaussians} gives
\[
 \phi_m^{(\theta_x+2\pi,\theta_y)}
 =\begin{cases}
 \phi_{m-1}^\theta,&1\le m<L,\\
 e^{-i\theta_y}\phi_{L-1}^\theta,&m=0,
 \end{cases}
 \qquad
 \phi_m^{(\theta_x,\theta_y+2\pi)}=\phi_m^\theta.
\]
After the rescaling, these become the uniform identities
\[
 \widetilde\phi_m^{(\theta_x+2\pi,\theta_y)}
 =e^{-i\theta_y/L}\widetilde\phi_{m-1}^\theta,
 \qquad
 \widetilde\phi_m^{(\theta_x,\theta_y+2\pi)}
 =e^{-2\pi im/L}\widetilde\phi_m^\theta,
\]
where $m-1$ is read modulo $L$.

To read off the coordinate identifications, represent the same
vector on opposite edges as $\sum_m v_m\widetilde\phi_m$.
The first identity gives
$v_m^{\rm right}=e^{i\theta_y/L}v_{m-1}^{\rm left}$,
and the second gives
$v_m^{\rm top}=e^{2\pi im/L}v_m^{\rm bottom}$.
By \eqref{mg:clockshift}, these are precisely
$u_x=e^{i\theta_y/L}S_L$ and $u_y=C_L$.
Finally, $S_LC_L=e^{-2\pi i/L}C_LS_L$ implies
$u_x(2\pi)u_y(0)=u_y(2\pi)u_x(0)$, as required.
\end{proof}

\subsection{Subbundles generated by periodic binary words}

For each degree $N$, take the exterior-power bundle
\begin{equation}\label{mg:manybody}
 \mathcal H_{L,N}=\bigwedge^N\mathcal E_L,
 \qquad \operatorname{rank}\mathcal H_{L,N}=\binom LN.
\end{equation}
In the occupation basis of Section \ref{sec:classical-quantum}, the transitions are
$\bigwedge^Nu_x$ and $\bigwedge^Nu_y$.

Fix integers $0<p<q$, let $M\geq1$ be an integer, and set $L=qM$, $N=pM$. Let $w=(w_0,\ldots,w_{q-1})$ be a primitive binary word, meaning that its least period is $q$, with $\sum_rw_r=p$. Its $M$-fold repetition has occupied sets
\begin{equation}\label{mg:rootsets}
 R_0=\{r+qk:w_r=1,\ 0\leq r<q,\ 0\leq k<M\},
 \qquad R_a=R_0+a\pmod L,
\end{equation}
where $0\leq a<q$. Write $\Psi_a$ for the occupation vector of $R_a$. Under the local
exterior-power identification, the span \eqref{pc:root-space} has fiber
\begin{equation}\label{mg:rootspace}
 \Psi_a=\bigwedge_{m\in R_a}^{\uparrow}e_m,
 \qquad
 \mathcal R_{w,M}(\theta)=\operatorname{span}\{\Psi_0,\ldots,\Psi_{q-1}\}
 \subset\mathcal H_{L,N}(\theta).
\end{equation}
With the Chern-number convention of Section \ref{sec:classical-quantum}, the corresponding
cohomology class is

\begin{equation}\label{mg:chernweil}
 c_1(E)=\left[\frac{i}{2\pi}\operatorname{tr}F\right].
\end{equation}
Here $F$ is the curvature of a unitary connection on the
\emph{parameter-space} bundle $E$. This is different from the
magnetic curvature $dA$ on the spatial torus. In particular,
$dA_{\mathrm{flat}}=0$ on each spatial torus does not imply
$F=0$: the latter measures how eigenspaces are glued as the
cohomology parameter varies. The class $\mathfrak u$ in \eqref{mg:eta} has integral one.
Thus $c_1(E)=p\mathfrak u$ means that the first Chern number of
$E$ is $p$; $\mathfrak u$ is a cohomology class, not an operator.

\begin{Pro}[The bundle of a primitive occupation word]
\label{mg:primitive-bundle}\label{mg:word-bundle}
Let $0<p<q$ be integers, $M\ge1$, $L=qM$, and $N=pM$.
For every primitive binary word $w$ of length $q$ containing $p$ occupied sites, the spaces \eqref{mg:rootspace} form a rank-$q$ Hermitian subbundle of $\bigwedge^N\mathcal E_L$.
For every rectangle of area $2\pi L$, and with the orientation and convention \eqref{mg:eta}--\eqref{mg:chernweil},
\[
 c_1(\mathcal R_{w,M})=p\mathfrak u,
 \qquad
 \bigl\langle c_1(\mathcal R_{w,M}),[\T^2_{\rm tw}]\bigr\rangle=p.
\]
No coprimality or interaction assumption is needed.
\end{Pro}

\begin{proof}
\emph{Step 1: the transition maps preserve the proposed fibers.}
Primitivity makes the $q$ occupied sets $R_a$ distinct, so their
occupation vectors $\Psi_a$ are orthonormal. It remains to check
that their span is preserved by both boundary identifications.
For the $x$ identification,
\begin{equation}\label{mg:firsttransition}
 \bigwedge^Nu_x=e^{iN\theta_y/L}\bigwedge^NS_L
 =e^{ip\theta_y/q}\bigwedge^NS_L.
\end{equation}
The operator $\bigwedge^NS_L$ shifts every occupied index by one.
Since $R_a+1=R_{a+1}$ modulo $L$, its action is
\[
 \left(\bigwedge^NS_L\right)\Psi_a=\varepsilon_a\Psi_{a+1},
 \qquad \varepsilon_a\in\{1,-1\},
\]
where $a+1$ is read modulo $q$. The sign only records the
permutation needed to restore increasing order in the exterior
product. Thus the restriction is a cyclic shift of the $q$ configurations,
with possible signs.

The $y$ identification preserves each basis line:
\[
 \left(\bigwedge^NC_L\right)\Psi_a=d_a\Psi_a,
 \qquad d_a=\exp\left(\frac{2\pi i}{L}\sum_{m\in R_a}m\right).
\]
Shifting all $N$ occupied indices by one changes their sum by
$N$ modulo $L$. Consequently
\[
 d_a=d_0e^{2\pi ipa/q}.
\]
Both transitions therefore preserve the span, which defines the
claimed rank-$q$ Hermitian subbundle.

\emph{Step 2: choose a basis with simple transition matrices.}
After $q$ shifts, each of the $p$ occupied residue classes modulo
$q$ undergoes an $M$-cycle. Each such cycle has sign $(-1)^{M-1}$.
It follows that
\begin{equation}\label{mg:fermionsign}
 \prod_{a=0}^{q-1}\varepsilon_a=(-1)^{p(M-1)}.
\end{equation}
Choose a unit complex number $\lambda$ with
$\lambda^q=(-1)^{p(M-1)}$. Set $z_0=1$ and recursively
$z_{a+1}=z_a\varepsilon_a/\lambda$. The product identity gives
$z_q=z_0$, and the rephased basis satisfies
\[
 \left(\bigwedge^NS_L\right)(z_a\Psi_a)
 =\lambda z_{a+1}\Psi_{a+1}.
\]
Thus the shift matrix becomes $\lambda S_q$, while the diagonal
matrix remains $d_0C_q^p$. Here $S_q$ shifts the $q$ basis labels
cyclically, and $C_q$ multiplies label $a$ by $e^{2\pi ia/q}$,
as in \eqref{mg:clockshift} with $L$ replaced by $q$.

Write $\lambda=e^{i\ell}$ and $d_0=e^{i\gamma}$.
On the cut square $Q$, multiply all coordinate vectors by
$\exp[-i(\ell\theta_x+\gamma\theta_y)/(2\pi)]$.
Its right-to-left boundary ratio is $e^{-i\ell}$ and its
top-to-bottom ratio is $e^{-i\gamma}$, so it removes exactly the
two constant phases. The resulting coordinate transitions are
\begin{equation}\label{mg:rootclutch}
 g_x(\theta_y)=e^{ip\theta_y/q}S_q,
 \qquad g_y(\theta_x)=C_q^p.
\end{equation}
This change of coordinates is made on $Q$, not asserted to be
periodic on the parameter torus. The corner compatibility follows
from $S_qC_q^p=e^{-2\pi ip/q}C_q^pS_q$.

\emph{Step 3: compute the first Chern class.}
The matrices in \eqref{mg:rootclutch} implement the two integral
identifications in the quotient of spatial $H^1$ described above.
Taking determinants gives the transition functions of the line
bundle $\det\mathcal R_{w,M}$, whose first Chern class equals
that of $\mathcal R_{w,M}$. The $q$-fold determinant of the scalar factor in $g_x$ is
$e^{ip\theta_y}$. This suggests the integer $p$; to verify its
sign with our coordinate convention, use the connection form
\begin{equation}\label{mg:connection}
 A^0=-i\frac pq\frac{\theta_x}{2\pi}\,d\theta_y\,I_q.
\end{equation}
At the left edge it is zero, and at the right edge it is
$-i(p/q)d\theta_y I_q=-(dg_x)g_x^{-1}$. Hence
$A_{\rm right}=g_xA_{\rm left}g_x^{-1}-(dg_x)g_x^{-1}$.
At the top and bottom the forms agree, and the constant matrix
$g_y$ commutes with them. Thus $d+A^0$ defines a unitary connection
on the subbundle. Its curvature is
\[
 F^0=dA^0+A^0\wedge A^0
 =-i\frac{p}{2\pi q}\,d\theta_x\wedge d\theta_y\,I_q.
\]
Taking the trace multiplies this scalar by $q$, so
\[
 \frac{i}{2\pi}\operatorname{tr}F^0
 =p\frac{d\theta_x\wedge d\theta_y}{(2\pi)^2}.
\]
The two angles each range over $[0,2\pi]$; the integral of this
form is therefore $p$. Since the second integral cohomology of
the torus has no torsion, this also identifies the integral class:
\begin{equation}\label{mg:mainconclusion}
 \operatorname{rank}\mathcal R_{w,M}=q,\qquad
 c_1(\mathcal R_{w,M})=p\mathfrak u,\qquad
 \langle c_1(\mathcal R_{w,M}),[\T^2_{\rm tw}]\rangle=p.
\end{equation}
\end{proof}

For the balanced word selected by the classical interaction, Proposition~\ref{mg:primitive-bundle} completes the proof of the Chern-class assertion in Theorem~\ref{intro:periodic}.

\begin{Lem}[Boundary loops and translation of periodic configurations]
\label{mg:shift-holonomy}
Use the rephased frame of \eqref{mg:rootclutch} and the connection
$d+A^0$ of \eqref{mg:connection} on $\mathcal R_{w,M}$.
The $x$-edge identification translates the occupation configuration
by one site, up to its common scalar phase; the $y$-edge
identification preserves each configuration and assigns successive
translates phases with ratio $e^{2\pi ip/q}$.
For the positively oriented coordinate loops based at $(0,0)$,
\[
 \operatorname{Hol}_x(d+A^0)=S_q^{-1},\qquad
 \operatorname{Hol}_y(d+A^0)=C_q^{-p}.
\]
Moreover, the determinant of the $x$-edge identification has winding
number $p$ as $\theta_y$ traverses its circle, recovering the Chern
number in \eqref{mg:mainconclusion}.
\end{Lem}
\begin{proof}
By \eqref{mg:rootsets}, shifting every occupied index by one sends
$R_a$ to $R_{a+1}$. The exterior-product signs are absorbed by
the rephasing used in \eqref{mg:rootclutch}, whose matrix $S_q$
therefore represents this cyclic permutation. The other matrix,
$C_q^p$, has successive diagonal entries with ratio
$e^{2\pi ip/q}=e^{2\pi iN/L}$, as computed above.

Along the bottom and left edges of $Q$, the pullback of $A^0$
vanishes, so parallel sections have constant coordinates there.
At the endpoint of either loop, identification back to the initial
fiber uses the inverse of the corresponding transition matrix.
Evaluating \eqref{mg:rootclutch} at the base point proves the
holonomy formulas. Finally,
\[
 \det g_x(\theta_y)=e^{ip\theta_y}\det S_q,
\]
whose winding number is $p$, with the orientation of
\eqref{mg:eta}.
\end{proof}

These holonomies refer to the specified connection, not necessarily
the Berry connection of a chosen operator family. Since $A^0$ has
nonzero curvature, they do not define a homotopy-invariant monodromy
representation of $\pi_1(\mathbb T^2_{\rm tw})$. The projective
commutation of $S_q$ and $C_q^p$ is compatible with the corner
identity because the scalar factor in $g_x(\theta_y)$ also varies.

\subsection{Unitary identifications and continuation of spectral subspaces}
The next lemma gives explicit unitary identifications between the
section spaces at different boundary parameters. Their exterior
powers will identify the corresponding finite-dimensional operators.

\begin{Lem}[Unitary identification between boundary conditions]\label{mg:magnetic-transport}
For $ab=2\pi L$ and $\theta=(\theta_x,\theta_y)\in Q$, define
\begin{equation}\label{vb:magnetic-translation}
 (T_\theta\psi)(x,y)
 =e^{-i\theta_x x/a}
 \psi\left(x-\frac{\theta_y}{b},y-\frac{\theta_x}{a}\right),
\end{equation}
where a zero-parameter magnetic section is extended to the plane using its boundary conditions.
Then $T_\theta$ is unitary from the zero-parameter one-variable section space onto the space with boundary parameter $\theta$, intertwines the magnetic operators, and satisfies
\begin{equation}\label{vb:transport-orbitals}
 T_\theta\phi_j^0=e^{-ij\theta_y/L}\phi_j^\theta,
 \qquad 0\le j<L.
\end{equation}
The induced transport
$\mathcal U_\theta=\bigoplus_{N=0}^L\bigwedge^N T_\theta$
preserves exterior degree.
\end{Lem}
\begin{proof}
Substituting the zero-parameter boundary conditions into \eqref{vb:magnetic-translation} gives \eqref{mg:bcx}--\eqref{mg:bcy} at parameter $\theta$.
Since $|\psi|^2$ is periodic in both spatial directions, translation preserves its integral over a fundamental rectangle.
The inverse translation with the inverse phase proves surjectivity, hence unitarity.
For $D_x=-i\partial_x+y$ and $D_y=-i\partial_y$, direct differentiation gives
$D_xT_\theta=T_\theta D_x$ and $D_yT_\theta=T_\theta D_y$.
These identities extend from smooth magnetic sections to the quadratic-form domains and intertwine the self-adjoint operators and their lowest eigenspaces.

In the eigenfunction series \eqref{mg:gaussians}, the mode indexed by $j+nL$ acquires the factor
\[
 \exp\!\left[-i\frac{2\pi(j+nL)}a\frac{\theta_y}b\right]
 =e^{-ij\theta_y/L}e^{-in\theta_y},
\]
while its momentum and Gaussian center change from $k_{j+nL}(0)$ to
$k_{j+nL}(0)-\theta_x/a$.
This proves \eqref{vb:transport-orbitals}, which is precisely the rephased frame used in \eqref{mg:lllclutch}.
The exterior powers preserve their fixed-degree subspaces by construction.
\end{proof}

The identifications are defined on the cut square $Q$; their edge
compatibility is governed by the transition matrices above. We next
prove invariance of the spectral bundle under a gapped deformation.

\begin{Pro}[Continuation of an isolated spectral bundle]
\label{sp:continuation}
Let $\mathcal H\to X$ be a finite-rank Hermitian bundle over a compact space and let $H_s$, $s\in[0,1]$, be a jointly continuous family of self-adjoint endomorphisms. If its ordered eigenvalues satisfy $\lambda_{r+1,s}(x)>\lambda_{r,s}(x)$ for every $(s,x)$, the lowest rank-$r$ spectral bundles at $s=0$ and $s=1$ are isomorphic and have the same Chern classes.
\end{Pro}
\begin{proof}
Locally choose a contour $\Gamma$ enclosing exactly the lowest $r$ eigenvalues. The Riesz formula
\[
 P_s(x)=\frac{1}{2\pi i}\int_\Gamma(z-H_s(x))^{-1}\,dz
\]
gives a continuous orthogonal projection; the local formulas agree on overlaps \cite{Kato1995}. Compactness gives a subdivision $0=s_0<\cdots<s_m=1$ such that $\sup_x\|P_{s_{j+1}}(x)-P_{s_j}(x)\|<1$. To see injectivity, let $v\in\operatorname{Ran}P_{s_j}(x)$ and suppose $P_{s_{j+1}}(x)v=0$.
Then $v=(P_{s_j}(x)-P_{s_{j+1}}(x))v$, which contradicts the strict norm bound unless $v=0$.
Both ranges have dimension $r$, so the restricted projection is a fiberwise isomorphism.
It varies continuously with $x$ and therefore defines a bundle isomorphism. Composing these continuous bundle isomorphisms proves the claim.
\end{proof}

\section{The \texorpdfstring{$V_1$}{V1} operator: minimizing degrees, spectral gaps, and Chern number}
\label{ml:section}

We study the nonnegative matrix $H_{L,\kappa}$ defined by
\eqref{ml:eq:QH}; its identification with the physical $V_1$ model
is taken from the literature. The proof has three algebraic steps.
A positive comparison operator is a small local perturbation of a
diagonal reference matrix. Spectral stability continues the
three-dimensional lowest spectral subspace; commutation with
$\Num$ and three known kernel vectors identify its exterior degree
and its eigenvalue. The resulting inequality on $\mathcal F_L$
then determines which degree minimizes $H_{L,\kappa}-c\Num$.
The Chern number follows by continuation of the isolated eigenbundle.
Throughout, $\kappa$ is fixed first and $L\in6\mathbb N$ is arbitrary.

\subsection{The operator and its diagonal decomposition}

Fix $L\in6\mathbb N$ and use the coordinate projections of
\eqref{intro:occupation-space}, initially at zero boundary parameters.
The operators $Q_R,H_{L,\kappa}$ are defined in \eqref{ml:eq:QH}. All these sums converge in operator norm. Using
$f_{L,\kappa}$ from \eqref{intro:periodic-amplitude}, the equivalent expression is
\begin{equation}\label{ml:eq:f-periodic}
 Q_R=\frac12\sum_{x\in(\Z+R)/L\Z}
 f_{L,\kappa}(x)c_{[R+x]}c_{[R-x]}.
\end{equation}
The factor $1/2$ compensates for the two orderings of a pair.
For comparison with the convention cited in Section \ref{sec:V1}, the
geometric parameter $\alpha$ in \cite{WY} corresponds to
our $\kappa$ after exchanging the coordinate axes and applying the
corresponding unitary change of section coordinates. The common positive scalar in the operator is absorbed into
our normalization, which is fixed by \eqref{ml:eq:QH}. We use the
periodized coefficients in \eqref{ml:eq:f-periodic} without truncation.
Multiplication by a scalar $g>0$ multiplies all gaps and
plateau widths by $g$.

We use the minimum energies $E_{H_{L,\kappa}}(N)$ and convex-dual
functions $\alpha_{H_{L,\kappa},L}$ defined in Section \ref{sec:classical-quantum}.
We use the diagonal decomposition \eqref{ml:eq:Hubbard-perturbation}.
Periodic amplitudes must be summed before taking their squared moduli:
different relative coordinates may represent the same unordered pair.

\begin{Lem}[Exact diagonal pair coefficient]\label{ml:lem:exact-diagonal-pair}
The operator $D_{L,\kappa}$ defined in Section~\ref{sec:V1}
is the diagonal part of $H_{L,\kappa}$ in the ordered
exterior-product basis, with the pair coefficient
\eqref{ml:eq:finite-V1-potential}.
\end{Lem}

\begin{proof}[Proof of Lemma~\ref{ml:lem:exact-diagonal-pair}]
Fix the unordered pair $\{i,j\}$, with $j-i=r$. It occurs in precisely the two summands indexed by $s=(i+j)/2$ and $s+L/2$, modulo $L$. In the first, the terms with $x\equiv r/2\pmod L$ have one ordering of the two annihilators, and those with $x\equiv-r/2\pmod L$ have the opposite ordering. Their combined coefficient is $f_{L,\kappa}(r/2)$. In the second, the combined coefficient is, up to its irrelevant sign, $f_{L,\kappa}((L-r)/2)$. On an binary-word basis vector, removing different unordered occupied pairs gives orthogonal vectors. Thus the diagonal expectation of $Q_s^*Q_s$ is the sum of the squared combined coefficients for occupied pairs. Summing over the two possible values of $s$ proves \eqref{ml:eq:finite-V1-potential}, with no additional factor of two.
\end{proof}

The potential $V_\kappa$ in \eqref{ml:eq:classical-U} is the
pointwise limit of these finite cyclic coefficients. Its limiting
variational energy and plateau widths are established in
Proposition~\ref{ml:prop:limits}.

\subsection{Positive comparison and the local perturbation bound}

Use $a_\kappa,\varepsilon_\kappa$ from Section~\ref{sec:V1},
and put $N_0=L/3$.
As in Section \ref{sec:V1}, the leading diagonal terms are
$a_\kappa A_L+\varepsilon_\kappa B_L$; their coefficients arise
from $x=y=1/2$ and $x=y=1$ in the squares $Q_R^*Q_R$.
Since $a_\kappa/\varepsilon_\kappa=\tfrac14e^{3\kappa^2/2}$ is
unbounded, we will reweight the positive squares before applying a
fixed spectral-stability theorem.

We use an interaction norm, rather than the norm of the whole operator, which can grow with $L$. For a decomposition $V=\sum_X\Phi(X)$ on the ring, each $\Phi(X)$ acts only on the tensor factors indexed by $X$. Define
\begin{equation}\label{ml:eq:interaction-norm}
 \|\Phi\|_*=\sup_j\sum_{X\ni j}
 |X|^3 e^{\diam(X)}\|\Phi(X)\|,
\end{equation}
where diameter uses the ring metric.
This norm controls the interaction at a fixed site, with sufficient
spatial decay for the spectral stability theorem.
To apply the spectral stability theorem, we identify the occupation
space with $(\mathbb C^2)^{\otimes L}$. Here a local operator means
an operator acting as the identity outside a specified set of tensor
factors. Fermionic signs require a choice of phases for this
identification. We use the standard choice from
\cite[\S1.2, Eq.~(1.3) and the discussion following Eq.~(1.5)]{WY}.
It preserves every occupation number. The elementary matrix
verification below explains why the full periodic interaction then
has the required local supports and translation symmetry.

For $\kappa\ge3$, $2\varepsilon_\kappa<a_\kappa$. Reduce the weight of the half-integer-center squares and define
\begin{equation}\label{ml:eq:comparison}
 K_{L,\kappa}
 =\sum_{R\in\Z/L\Z}Q_R^*Q_R
 +\frac{2\varepsilon_\kappa}{a_\kappa}
 \sum_{R\in(\Z+1/2)/L\Z}Q_R^*Q_R.
\end{equation}
Since all coefficients in this weighted sum are positive,
\begin{equation}\label{ml:eq:comparison-order}
 0\le K_{L,\kappa}\le H_{L,\kappa},
 \qquad \ker K_{L,\kappa}=\ker H_{L,\kappa}.
\end{equation}
In particular, the reweighting preserves the common kernel of the operators $Q_R$.

The leading terms in \eqref{ml:eq:comparison} give
\begin{equation}\label{ml:eq:reference-plus-tail}
 \varepsilon_\kappa^{-1}K_{L,\kappa}=2A_L+B_L+\mathcal R_{L,\kappa}.
\end{equation}
The reweighting gives a reference operator with coefficients bounded independently of $\kappa$ after division by $\varepsilon_\kappa$. Equality of the kernels will later allow the resulting bound to be transferred back to the original operator.

\begin{Lem}\label{ml:lem:tail}
After that unitary change, $\mathcal R_{L,\kappa}$ is translation invariant, commutes with $\Num$, and admits a decomposition
\[
 \mathcal R_{L,\kappa}
 =\sum_{\varnothing\ne X\subset\Z/L\Z}\Phi_{L,\kappa}(X),
 \qquad \Phi_{L,\kappa}(X)^*=\Phi_{L,\kappa}(X),
\]
where $\Phi_{L,\kappa}(X)$ acts only on the tensor factors indexed by $X$.
Thus $\mathcal R_{L,\kappa}$ is a single operator on
$(\mathbb C^2)^{\otimes L}$, whereas $\Phi_{L,\kappa}$ denotes the
family of local operators $X\mapsto\Phi_{L,\kappa}(X)$ whose sum is
that operator. Such a decomposition need not be unique; the proof
constructs one satisfying
\begin{equation}\label{ml:eq:small-tail}
 \|\Phi_{L,\kappa}\|_*\le C e^{-2\kappa^2}
 \qquad(\kappa\ge1),
\end{equation}
where $C$ is an explicit finite constant independent of $L$.
The norm in this estimate is the weighted norm of the family defined
in \eqref{ml:eq:interaction-norm}, not the operator norm of its sum.
In particular, it measures the total contribution near each site,
including a penalty for the size and diameter of the supporting set.
It implies $\|\mathcal R_{L,\kappa}\|\le
L\|\Phi_{L,\kappa}\|_*$; the factor $L$ explains why we use the
local bound when applying spectral stability uniformly in $L$.
\end{Lem}

\begin{proof}[Proof of Lemma~\ref{ml:lem:tail}]
The purpose of the first step is to make locality explicit in a
tensor-product basis. Small matrix coefficients alone do not give
the required estimate: the norm in \eqref{ml:eq:interaction-norm}
also weights the size and diameter of each support. We therefore
represent the exterior contractions by explicit matrices, with a
choice of basis phases adapted to the periodic boundary. This
separates the argument into a support estimate and a Gaussian
coefficient estimate. We use the Jordan--Wigner representation in
\cite[\S1.2, Eq.~(1.3)]{WY}; its periodic boundary correction is
discussed there after Eq.~(1.5). The identities below record the
form needed for all separations in the present operator.

\smallskip
\noindent\emph{Step 1: A local representation.}
We first relate the coordinates used here to Section~\ref{mg:section}.
On the cut parameter square, the coordinate vector $e_j\in\mathbb C^L$
represents the orbital
$\widetilde\phi_j^\theta=e^{-ij\theta_y/L}\phi_j^\theta
=T_\theta\phi_j^0$ by \eqref{vb:transport-orbitals}.
Consequently the occupation vector for a binary word is
\[
 \Psi_\omega(\theta)
 =\bigwedge_{j:\omega_j=1}^{\uparrow}\widetilde\phi_j^\theta,
 \qquad
 \Psi_\omega(\theta)=\mathcal U_\theta\Psi_\omega(0).
\]
These are the same exterior-product coordinates as in Section \ref{sec:classical-quantum} and \eqref{mg:rootspace}; the latter selects only the $q$ translated periodic words, whereas here we use all binary words and all degrees.
The estimates below are carried out at $\theta=0$.

The tensor-product basis is this occupation basis with a fixed
choice of signs. More precisely, we identify
\[
 |\omega_0\rangle\otimes\cdots\otimes|\omega_{L-1}\rangle
 \longmapsto
 (-1)^{\sum_{j=0}^{L-1}j\omega_j}\Psi_\omega(0).
\]
This diagonal unitary changes neither occupation numbers nor the
span of any collection of occupation vectors. With the unsigned
identification, exterior contraction has the usual preceding-site
factors $\operatorname{diag}(1,-1)$. The displayed sign change
multiplies its $j$th matrix by $(-1)^j$, giving exactly
$c_j=(\prod_{i<j}Z_i)a_j$ with the convention for $Z$ below.
Thus the tensor factors label the same orbitals as the coordinate
basis in Section~\ref{mg:section}; they do not introduce a different
one-particle space. The further parity-dependent unitary
$U_\epsilon$ introduced below handles the periodic boundary signs
and is distinct from the geometric transport $\mathcal U_\theta$.

On $(\mathbb C^2)^{\otimes L}$, let
$a=\left(\begin{smallmatrix}0&1\\0&0\end{smallmatrix}\right)$ and
$Z=\operatorname{diag}(-1,1)$. Recall that $Q_R$ is a weighted sum of the pair contractions
$c_{[R+x]}c_{[R-x]}$. We now represent each such contraction by
a matrix supported on the interval joining its two indices.
For $p=R-x$, $q=R+x$, set
\[
 A_{R,x}=-a_{[p]}\left(\prod_{j=p+1}^{q-1}Z_{[j]}\right)a_{[q]}.
\]
To interpret this formula, write the binary-word basis as
$|\omega_0\rangle\otimes\cdots\otimes|\omega_{L-1}\rangle$,
where $a|1\rangle=|0\rangle$ and $a|0\rangle=0$.
The two endpoint matrices replace the entries at $p$ and $q$
by zero, and the result vanishes unless both entries were one.
The intervening diagonal matrices leave the binary entries
unchanged and supply the sign required by exterior contraction:
$Z|\omega_j\rangle=(-1)^{1-\omega_j}|\omega_j\rangle$.
In particular, for $0\le p<q<L$, cancellation of the common
initial products in $c_qc_p$ gives exactly
$-a_p(\prod_{p<j<q}Z_j)a_q$. Thus $A_{R,x}$ is the interval
representation of a pair of contractions, with $R$ its midpoint
and $2x$ its length. Its support includes the intervening indices,
not just the two endpoints.

This operator has norm at most one and is supported on the image of
$[R-x,R+x]$ in the ring. Products retain repeated factors when the
interval winds around the ring; equal endpoint sites give zero.
These are the only properties of the matrices needed for the estimates.

We record briefly why this representation also handles the boundary.
Use $c_j=(\prod_{i<j}Z_i)a_j$. On the subspace
$\prod_j Z_j=\epsilon$, choose $e^{i\theta_\epsilon}=-\epsilon$
and define the diagonal unitary
$U_\epsilon|\omega\rangle=
 e^{-i\theta_\epsilon\sum_jj\omega_j/L}|\omega\rangle$.
Multiplication of the diagonal signs gives, with
$k=\lfloor p/L\rfloor+\lfloor q/L\rfloor$,
\[
 U_\epsilon c_{[q]}c_{[p]}U_\epsilon^*
 =(-\epsilon)^k e^{i\theta_\epsilon([p]+[q])/L}A_{R,x}
 =e^{2i\theta_\epsilon R/L}A_{R,x}.
\]
In particular, substituting this identity into the definition
\eqref{ml:eq:QH} gives, on the parity sector $\epsilon$,
\[
 U_\epsilon Q_RU_\epsilon^*
 =e^{2i\theta_\epsilon R/L}
 \sum_{\substack{x>0\\x-R\in\mathbb Z}}
 x e^{-\kappa^2x^2}A_{R,x}.
\]
Here $Q_R$ preserves parity, since it lowers the exterior degree
by two. Thus $A_{R,x}$ is the local matrix for the individual
pair-removal term in $Q_R$, rather than a new interaction.
The phase is independent of $x$ and cancels in $Q_R^*Q_R$.
Taking the direct sum over $\epsilon=\pm1$ gives the required
unitary change, preserving every $n_j$. This is the standard
periodic correction discussed in \cite[\S1.2]{WY}; the displayed
identity applies to all separations, including winding intervals.

\smallskip
\noindent\emph{Step 2: Expand the remainder into local terms.}
The leading relative separation is $u=1/2$ for half-integer centers
and $u=1$ for integer centers. For $x\in u+\mathbb N_0$, put
\[
 r_u(x)=\frac xu e^{-\kappa^2(x^2-u^2)},
 \qquad \omega_{1/2}=2,\quad\omega_1=1.
\]
For a center $R\in(u+\Z)/L\Z$, the preceding identity reads
\[
 U_\epsilon Q_RU_\epsilon^*
 =e^{2i\theta_\epsilon R/L}\,
 u e^{-\kappa^2u^2}
 \sum_{x\in u+\mathbb N_0}r_u(x)A_{R,x}.
\]
Thus $r_u(x)$ is the coefficient of the separation $x$ relative
to the leading coefficient at $x=u$. After applying the weights
in \eqref{ml:eq:comparison} and dividing by $\varepsilon_\kappa$,
the contribution of this center to $K_{L,\kappa}/\varepsilon_\kappa$
is the normalized square
\[
 \omega_u\left(\sum_{x\in u+\mathbb N_0}r_u(x)A_{R,x}\right)^*
 \left(\sum_{y\in u+\mathbb N_0}r_u(y)A_{R,y}\right).
\]
Since $r_u(u)=1$, its $(u,u)$ term is precisely the corresponding
term of $2A_L+B_L$. Removing these leading terms, therefore, gives
\[
 \mathcal R_{L,\kappa}
 =\sum_{u\in\{1/2,1\}}\sum_{R\in(u+\Z)/L\Z}
 \sum_{\substack{x,y\in u+\mathbb N_0\\(x,y)\ne(u,u)}}
 \omega_u r_u(x)r_u(y)B_{R;x,y},
 \qquad
 B_{R;x,y}=\frac{A_{R,x}^*A_{R,y}+A_{R,y}^*A_{R,x}}2.
\]
Here the sum is over ordered pairs $(x,y)$. Symmetrization does not
change the sum, because its coefficients and index set are symmetric
in $x,y$. Each $B_{R;x,y}$ is self-adjoint and has norm at most one.

\smallskip
\noindent\emph{Step 3: Count the terms meeting a fixed site.}
Fix $u,x,y$, and put $m=\max(x,y)$ and $\ell=2m+1$.
Both pair matrices, and hence $B_{R;x,y}$, are supported in
\[
 X_{R,m}=\{[R-m],[R-m+1],\ldots,[R+m]\}.
\]
Thus $|X_{R,m}|\le\ell$ and
$\diam(X_{R,m})\le2m$ in the ring metric. These bounds remain true across the boundary cut and for repeated
windings. If $\ell\ge L$, the image has $L\le\ell$ sites;
in all cases, the ring distance never exceeds the distance between
chosen integer lifts. For any fixed site $j$, there are at
most $\ell$ centers $R\in(u+\Z)/L\Z$ with $j\in X_{R,m}$.
Indeed, such a center satisfies $R\equiv j-d\pmod L$ for one of the
$\ell$ offsets $d\in\{-m,-m+1,\ldots,m\}$. Repeated residues only
reduce this number, so the same counting applies to winding intervals.

Define $\Phi_{L,\kappa}(X)$ by adding all the above self-adjoint terms
assigned to $X_{R,m}=X$; explicitly,
\[
 \Phi_{L,\kappa}(X)=
 \sum_{u\in\{1/2,1\}}\sum_{R\in(u+\Z)/L\Z}
 \sum_{\substack{x,y\in u+\mathbb N_0,\ (x,y)\ne(u,u)\\
 X_{R,\max(x,y)}=X}}
 \omega_u r_u(x)r_u(y)B_{R;x,y}.
\]
Summing over $X$ recovers exactly the expansion of
$\mathcal R_{L,\kappa}$ in Step~2: each term is assigned to one
specified supporting set, even if its minimal support is smaller.
The triangle inequality, followed by the
center count, gives
\[
 \begin{aligned}
 \sum_{X\ni j}|X|^3e^{\diam(X)}\|\Phi_{L,\kappa}(X)\|
 &\le\sum_{u\in\{1/2,1\}}\omega_u
 \sum_{\substack{x,y\in u+\mathbb N_0\\(x,y)\ne(u,u)}}
 r_u(x)r_u(y)
 \sum_{R:\,j\in X_{R,m}}\ell^3e^{2m}\\
 &\le\sum_{u\in\{1/2,1\}}\omega_u
 \sum_{\substack{x,y\in u+\mathbb N_0\\(x,y)\ne(u,u)}}
 \ell^4e^{2m}r_u(x)r_u(y).
 \end{aligned}
\]
The estimates may first be applied to finite sums; Gaussian decay
makes the displayed majorant summable and justifies the limiting
decomposition. Taking the supremum over $j$ gives
\begin{equation}\label{ml:eq:double-sum}
 \|\Phi_{L,\kappa}\|_*\le
 \sum_{u\in\{1/2,1\}}\omega_u
 \sum_{\substack{x,y\in u+\mathbb N_0\\(x,y)\ne(u,u)}}
 (2\max(x,y)+1)^4e^{2\max(x,y)}r_u(x)r_u(y).
\end{equation}
The fourth power consists of three powers from the interaction norm
and one from counting centers.

\smallskip
\noindent\emph{Step 4: Extract the small Gaussian factor.}
For every nonleading pair, $x^2+y^2-2u^2\ge2$. For $\kappa\ge1$,
\[
 r_u(x)r_u(y)\le e^{-2\kappa^2}\frac{xy}{u^2}
 e^{-(x^2+y^2-2u^2-2)}.
\]
An explicit admissible constant $C$ is $e^2$ times the right-hand side of \eqref{ml:eq:double-sum} evaluated at $\kappa=1$.
This Gaussian double series converges. Translation invariance is immediate from the sum over centers, and every term contains equally many exterior multiplications and contractions.
\end{proof}

\subsubsection{Exact kernel vectors}
\label{ml:sec:zero}

The proof uses three known kernel vectors, but not their explicit
formulas. We state precisely the input needed from the torus
Laughlin construction \cite{HaldaneRezayi1985}.

\begin{Lem}[Kernel vectors at and below degree $L/3$]\label{ml:lem:zero}
For every $\kappa>0$ and $L\in6\mathbb N$,
\[
 \dim\ker\!\left(H_{L,\kappa}|_{\mathcal H_{L,N_0}}\right)\ge3.
\]
For every $0\le N\le N_0$, $E_{H_{L,\kappa}}(N)=0$.
\end{Lem}

\subsection{Applying the spectral stability theorem}
\label{ml:sec:stability}

We seek a reference operator whose lowest eigenspace consists of
the three translates of $(100)^{L/3}$ and is separated from the
rest of the spectrum uniformly in $L$. The leading diagonal operator
$2A_L+B_L$ alone does not have this property: every binary word
whose occupied sites are separated by at least three positions has
zero energy, including the empty word and many words with fewer
than $N_0$ occupied sites. We therefore subtract $t\Num$ to
favor larger exterior degree, while retaining the penalties for
occupied sites at distance one or two. Adding the scalar $tN_0$
normalizes the energy of the three desired words to zero. The local decomposition below will be used in the proof of
Lemma~\ref{ml:lem:stability} to verify its kernel and gap when
$0<t<3$.

For $0<t<3$, put $S_j=n_j+n_{j+1}+n_{j+2}$ and define
\begin{equation}\label{ml:eq:reference}
 \mathcal L_{L,t}=2A_L+B_L-t\Num+\frac{t L}{3}
 =\sum_j h_t(S_j),
 \qquad h_t(s)=\binom{s}{2}-\frac t3(s-1).
\end{equation}
To verify the equality, use $n_j^2=n_j$ and the commutativity of
the occupation projections to obtain
\[
 \binom{S_j}{2}=\frac{S_j^2-S_j}{2}
 =n_jn_{j+1}+n_jn_{j+2}+n_{j+1}n_{j+2}.
\]
Summing over the cyclic indices, each adjacent pair occurs twice
and each pair at distance two occurs once. Also, each $n_j$
occurs in exactly three of the sums $S_j$. Hence
\[
 \sum_j\binom{S_j}{2}=2A_L+B_L,
 \qquad \sum_j S_j=3\Num,
\]
which gives the asserted decomposition after substituting the
polynomial $h_t$.

On the fixed-degree space $\mathcal H_{L,N_0}$, the added term
$-t(\Num-N_0)$ vanishes. Its purpose is to separate this space's
three reference minimizers from all other degrees, so that spectral
stability can be applied on the full exterior algebra. Indeed,
\eqref{ml:eq:reference-plus-tail} gives
\[
 \mathcal L_{L,t}+\mathcal R_{L,\kappa}
 =\varepsilon_\kappa^{-1}K_{L,\kappa}-t(\Num-N_0).
\]
Thus the same perturbation estimate will control comparisons
between different degrees. In the main proof, spectral stability
and the exact kernel vectors will show that this operator is
nonnegative. Together with $H_{L,\kappa}\ge K_{L,\kappa}$ and
$H_{L,\kappa}\ge0$, this yields the density-$1/3$ plateau for
$0<c<t\varepsilon_\kappa$. We retain the parameter $t$, rather
than fixing one convenient value, because taking $t$ arbitrarily
close to $3$ gives the desired asymptotic lower bound for the plateau
width; the required lower bound on $\kappa$ may depend on $t$.

The following lemma gives the spectral stability needed below; its
proof verifies the local reference conditions directly.

\begin{Lem}[Spectral stability of the reference operator]
\label{ml:lem:stability}
For each $t\in(0,3)$ there exist $\epsilon_t>0$ and $g_t>0$ such that any cyclically translation-invariant self-adjoint interaction $\Phi_L$ with $\sup_L\|\Phi_L\|_*<\epsilon_t$ has the following property. For every $s\in[0,1]$, the lowest three eigenvalues of
\[
 \mathcal L_{L,t}+s\sum_X\Phi_L(X)
\]
form a spectral cluster separated from the remaining spectrum by at least $g_t$, uniformly in $L\in6\mathbb N$. Its rank-three spectral projection depends continuously on $s$.
\end{Lem}

With the reference and interaction norm fixed, Lemma~\ref{ml:lem:tail} allows a choice $\kappa_t\ge3$, independent of $L$, satisfying
\begin{equation}\label{ml:eq:verified-smallness}
 C e^{-2\kappa_t^2}<\epsilon_t.
\end{equation}

\subsubsection{Limits and classical widths}
\label{ml:sec:classical}
\begin{Pro}[Limits of minimum eigenvalues per site and the diagonal comparison]\label{ml:prop:limits}
For each fixed $\kappa>0$, the limit $\alpha_\kappa(c)=\lim_{L\to\infty,\,6\mid L}\alpha_{H_{L,\kappa},L}(c)$ exists for every $c$. The corresponding diagonal limit also exists. For $\kappa\geq3$, the limiting diagonal energy at fixed density is $\beta_{V_\kappa}$ for the potential $V_\kappa$ in \eqref{ml:eq:classical-U}. For $\kappa\ge3$, every rational $0<p/q<1$ has classical width
\begin{equation}\label{ml:eq:classical-width}
 W^{\mathrm{cl}}_{p/q}(\kappa)
 =q\sum_{k\ge1}k\bigl[V_\kappa(kq-1)
 -2V_\kappa(kq)+V_\kappa(kq+1)\bigr]>0.
\end{equation}
In particular,
\begin{equation}\label{ml:eq:classical-one-third}
 W^{\mathrm{cl}}_{1/3}(\kappa)
 =3e^{-2\kappa^2}
 -\frac{27}{2}e^{-9\kappa^2/2}
 +O(e^{-8\kappa^2}).
\end{equation}

\end{Pro}

\subsubsection{The lowest-eigenspace bundle over boundary parameters}
\label{vb:section}

For the geometric part of the main proof, we relate the three-dimensional lowest eigenspace to the periodic-state bundle of Section~\ref{mg:section}. Unitary identification between boundary conditions extends the estimates over all boundary parameters, and a compatible gapped path identifies the Chern class.

Retain the geometry of Section~\ref{sec:classical-quantum} and
the normalization of Section~\ref{sec:V1}
and the notation $N_0$ introduced above.
Use the magnetic boundary conventions and the oriented parameter torus
of Section~\ref{mg:section}. Let $\widehat H_{L,\kappa}(\theta)$
be the full torus $V_1$ interaction on
$\bigwedge^N\mathcal E_L(\theta)$, with energy normalization fixed
by \eqref{ml:eq:QH}. We use its standard projected
Trugman--Kivelson representation from \cite[\S1.1]{WY}, rather
than deriving it from continuous coordinates. The two properties
needed below are invariance under the unitary transport applied simultaneously to all variables and zero interaction energy for states vanishing
to third order at every pair coincidence. Both are properties of
the projected contact interaction; the latter follows from its
first-relative-derivative quadratic form. At zero boundary parameters its matrix
is \eqref{ml:eq:QH}, with the periodized coefficients
\eqref{ml:eq:f-periodic}.

\begin{Lem}[Unitary equivalence for the boundary-parameter family]
\label{vb:transport}
Let $T_\theta$ and $\mathcal U_\theta$ be the unitary transport
and its direct sum of exterior powers from
Lemma~\ref{mg:magnetic-transport}.  For the full projected contact
interaction,
\begin{equation}\label{vb:unitary-equivalence}
 \widehat H_{L,\kappa}(\theta)
 =\mathcal U_\theta H_{L,\kappa}\mathcal U_\theta^{-1},
 \qquad
 \mathcal U_\theta\Num=\Num\mathcal U_\theta.
\end{equation}
In the rephased eigenfunction frame of \eqref{mg:lllclutch}, its
coefficient matrix is $H_{L,\kappa}$, independent of $\theta$.
\end{Lem}

We will apply Proposition~\ref{sp:continuation} to identify the interacting
bundle with the periodic-state bundle.

The transport $T_\theta$ is defined on the cut parameter square. Its edge identifications are those of the bundle of lowest eigenspaces, so the parameter-independent spectrum is compatible with the nonzero Chern class. The deformation in the main proof uses a gap at fixed exterior degree; the plateau statement follows separately from the unitary equivalence on the full exterior algebra and the energy estimates.

\subsection{From spectral stability to the global plateau and Chern number}

\begin{proof}[Proof of Theorem~\ref{intro:quantum}: energy and gaps]
Fix $t\in(0,3)$ and choose $\kappa\ge\kappa_t$ so that \eqref{ml:eq:verified-smallness} holds.

\smallskip
\noindent\emph{Step 1: Continue the three-dimensional subspace and fix its exterior degree.}
Work in the degree-preserving tensor-product coordinates used in
Lemma~\ref{ml:lem:tail}; this unitary change preserves the spectra.
Lemmas~\ref{ml:lem:tail} and~\ref{ml:lem:stability} apply to
\[
 \mathcal L_{L,t}+s\mathcal R_{L,\kappa},\qquad 0\le s\le1,
\]
and denote its projection onto the lowest spectral cluster by $P(s)$. The interpolation commutes with $\Num$, hence so does $P(s)$. If $P_N$ is the projection onto number $N$, then
\[
 s\longmapsto\Tr(P_NP(s))
\]
is continuous and integer valued.
Indeed, $P_N$ and $P(s)$ commute, so their product is the orthogonal projection onto $\mathcal H_{L,N}\cap\operatorname{Ran}P(s)$, and its trace is the dimension of this intersection.
At $s=0$, that dimension is three for $N=N_0$ and zero otherwise.
It remains constant on $[0,1]$, proving that $\operatorname{Ran}P(s)\subset\mathcal H_{L,N_0}$ throughout the deformation.

\smallskip
\noindent\emph{Step 2: Identify the continued eigenvalue and obtain a global operator bound.}

At $s=1$ the operator is
\begin{equation}\label{ml:eq:endpoint-reference}
 \mathcal L_{L,t}+\mathcal R_{L,\kappa}
 =\varepsilon_\kappa^{-1}K_{L,\kappa}-t(\Num-N_0).
\end{equation}
On $\mathcal H_{L,N_0}$ this operator is $\varepsilon_\kappa^{-1}K_{L,\kappa}$ and is therefore nonnegative.
Its three-dimensional lowest spectral subspace belongs to this fixed-degree subspace by Step 1, so the operator has no negative eigenvalues anywhere in the exterior algebra.
Lemma~\ref{ml:lem:zero} and the common-kernel identity supply three independent kernel vectors.
The lowest three eigenvalues are consequently all zero, and the isolation of this cluster excludes any further zero eigenvalue.
Thus the continued subspace is exactly the kernel of the shifted
operator in \eqref{ml:eq:endpoint-reference}.  Transport this projection
and the operator inequality back to the original exterior-product coordinates.
Let $P$ denote the projection onto
$\ker(H_{L,\kappa}|_{\mathcal H_{L,N_0}})$, extended by zero on every
other fixed-degree subspace.  It has rank three.  It is not the
projection onto the kernel on the full exterior algebra of $H_{L,\kappa}$, which
also contains states with smaller exterior degree.
The separation from the remaining spectrum implies
\begin{equation}\label{ml:eq:strong-comparison}
 \varepsilon_\kappa^{-1}K_{L,\kappa}-t(\Num-N_0)
 \ge g_t(I-P).
\end{equation}
Combining this with $H_{L,\kappa}\ge K_{L,\kappa}$ gives
\begin{equation}\label{ml:eq:strong-full}
 H_{L,\kappa}\ge t \varepsilon_\kappa(\Num-N_0)+g_t \varepsilon_\kappa(I-P).
\end{equation}
On $N=N_0$, the common kernel of $K_{L,\kappa}$ and
$H_{L,\kappa}$ is exactly $\operatorname{ran}P$. This yields
\begin{equation}\label{ml:eq:neutral-gap}
 \Delta_{H_{L,\kappa}}(N_0)\ge g_t\varepsilon_\kappa.
\end{equation}

\smallskip
\noindent\emph{Step 3: Compare all exterior degrees.}
At degree $N_0$, the minimum eigenvalue of $H_{L,\kappa}-c\Num$ is $-cN_0$.

For $N<N_0$, the identity $H_{L,\kappa}=\sum_RQ_R^*Q_R\ge0$
gives, when $c>0$,
\[
 E_{H_{L,\kappa}}(N)-c N\ge-c N_0+c(N_0-N)>-c N_0.
\]
For $N>N_0$ and $c<t \varepsilon_\kappa$, the operator lower bound gives
\[
 E_{H_{L,\kappa}}(N)-c N
 \ge-c N_0+(t \varepsilon_\kappa-c)(N-N_0)>-c N_0.
\]
The three kernel vectors at $N=N_0$ attain $-c N_0$. Hence, for $0<c<t\varepsilon_\kappa$,
\begin{equation}\label{ml:eq:main-plateau}
 \min\spec(H_{L,\kappa}-c\Num)=-c N_0,
 \qquad \alpha_{H_{L,\kappa},L}(c)=\frac{c}{3}.
\end{equation}
More precisely, \eqref{ml:eq:strong-full} gives separation at least $g_t \varepsilon_\kappa$ for $N>N_0$ and for the orthogonal complement of $P$ inside $N=N_0$. The subspaces $N<N_0$ are separated by at least $c$. Therefore
\begin{equation}\label{ml:eq:grand-gap}
 \gap(H_{L,\kappa}-c\Num)\ge\min\{c,g_t\varepsilon_\kappa\}.
\end{equation}

\smallskip
\noindent\emph{Step 4: Exclude plateaux below density $1/3$.}

For the final assertion, no lower bound on $\kappa$ is needed. At every $\kappa>0$, positivity and Lemma~\ref{ml:lem:zero} give $E_{H_{L,\kappa}}(N)=0$ for all $0\le N\le N_0$. When $c<0$, the lowest eigenspace is the one-dimensional degree-zero subspace. When $c>0$, the minimum eigenvalue of $H_{L,\kappa}-c\Num$ at degree $N_0$ is strictly smaller than its minimum in every degree $N<N_0$. At $c=0$, kernel vectors coexist at all exterior degrees $0\le N\le N_0$, which does not define an open plateau. This proves the energy and gap assertions of Theorem~\ref{intro:quantum}.
\end{proof}

\begin{proof}[Proof of Theorem~\ref{intro:quantum}: the limiting width]
By Proposition~\ref{ml:prop:limits}, the limits used below exist. For $\kappa\ge\kappa_t$, the finite-volume energy estimate proved above gives $\alpha_\kappa(c)=c/3$ on $(0,t\varepsilon_\kappa)$. For every $\kappa>0$, no minimizing exterior degree at finite $L$ has density in $(0,1/3)$ when $c\ne0$. On $c>0$, all the slopes of $\alpha_{H_{L,\kappa},L}$ are at least $1/3$; their limit has the same one-sided slope bound. On $c<0$ the limit is zero. Hence no fractional plateau below $1/3$ reappears in the limit.

For each fixed $t<3$, the width $W^{\mathrm{full}}_{1/3}(\kappa)$ defined in Section~\ref{sec:V1} is at least $t\varepsilon_\kappa$ for all sufficiently large $\kappa$. Since $W^{\mathrm{cl}}_{1/3}(\kappa)/\varepsilon_\kappa\to3$, one obtains
\[
 \liminf_{\kappa\to\infty}
 \frac{W^{\mathrm{full}}_{1/3}(\kappa)}
 {W^{\mathrm{cl}}_{1/3}(\kappa)}\ge t/3.
\]
Letting $t\uparrow3$ proves the width comparison in Theorem~\ref{intro:quantum}.

\end{proof}

\begin{proof}[Proof of Theorem~\ref{intro:quantum}: the bundle over boundary parameters]
Lemma~\ref{vb:transport} and the energy and gap estimates proved above give the spectral assertions. To compute the Chern class, we verify a continuous deformation with an isolated lowest spectral subspace over the entire parameter torus in the fixed-fixed-degree subspace $N=N_0$.

\smallskip
\noindent\emph{Step 1: Check that the deformation is defined on the parameter torus.}
On the full exterior algebra $\bigoplus_{N=0}^L\bigwedge^N\mathbb C^L$,
define the unitary operators
\begin{equation}\label{eq:exterior-clockshift}
 \mathsf S:=\bigoplus_{N=0}^L\bigwedge^NS_L,
 \qquad
 \mathsf C:=\bigoplus_{N=0}^L\bigwedge^NC_L.
\end{equation}
Here $S_L,C_L$ are defined in \eqref{mg:clockshift}, and the
zeroth exterior powers act as the identity on $\mathbb C$.
Thus, on degree $N$, these operators are precisely
$\bigwedge^NS_L$ and $\bigwedge^NC_L$. Explicitly,
\[
 \begin{aligned}
 \mathsf S(e_{j_1}\wedge\cdots\wedge e_{j_N})
 &=e_{j_1+1}\wedge\cdots\wedge e_{j_N+1},\\
 \mathsf C(e_{j_1}\wedge\cdots\wedge e_{j_N})
 &=e^{2\pi i(j_1+\cdots+j_N)/L}
   e_{j_1}\wedge\cdots\wedge e_{j_N}.
 \end{aligned}
\]
Indices are taken modulo $L$, with the usual exterior-product
sign when the shifted indices are reordered. The direct-sum
notation is useful because $Q_R$ lowers degree by two: on degree
$N\ge2$, for example, $\mathsf S Q_R\mathsf S^*$ means
$(\bigwedge^{N-2}S_L)Q_R(\bigwedge^NS_L)^*$. All coefficient matrices here are in the original exterior-product eigenfunction frame. The tensor-product change of coordinates used to prove stability is transported back before making these identifications; it preserves exterior degree. In particular, $\mathcal R_{L,\kappa}$ is expressed in this frame by \eqref{ml:eq:reference-plus-tail}. Conjugation of the exterior contractions gives
\begin{equation}\label{vb:pair-covariance}
 \mathsf S Q_R\mathsf S^*=Q_{R+1},\qquad
 \mathsf C Q_R\mathsf C^*=e^{-4\pi iR/L}Q_R.
\end{equation}
For the second identity, every term has sum of the two basis indices $2R$ modulo $L$, including the winding terms in \eqref{ml:eq:QH}. Consequently $H_{L,\kappa}$ commutes with both operators. The same is true of $K_{L,\kappa}$: translation preserves the integer and half-integer center classes, and their separate weights do not affect the diagonal-symmetry identity. The reference $\mathcal L_{L,t}$ also commutes with both operators, because it is a translation-invariant polynomial in the occupations and $\Num$. Equation~\eqref{ml:eq:reference-plus-tail} then gives these commutation relations for $\mathcal R_{L,\kappa}$ as well.

In the fixed-fixed-degree subspace, the transition matrices from \eqref{mg:lllclutch} are $e^{iN_0\theta_y/L}\bigwedge^{N_0}S_L$ and $\bigwedge^{N_0}C_L$. Every matrix just listed commutes with these transitions. Hence constant coefficient matrices and the following paths define continuous self-adjoint bundle endomorphisms; they are compatible at both edges of the parameter square.

\smallskip
\noindent\emph{Step 2: Deform the classical subspace to the interacting kernel.}

The first path is
\begin{equation}\label{vb:first-path}
 A_s=\left.(\mathcal L_{L,t}
 +s\mathcal R_{L,\kappa})\right|_{N=N_0},
 \qquad 0\leq s\leq1.
\end{equation}
Its endpoints are the diagonal reference with periodic minimizing configurations $(100)^{L/3}$ and its two translates, and $K_{L,\kappa}/\varepsilon_\kappa$ in this subspace. The energy proof has already shown that the isolated rank-three spectral subspace of Lemma~\ref{ml:lem:stability} stays in $N=N_0$ throughout this path. Restriction therefore retains its external gap. The coefficient matrices are independent of $\theta$, so the gap bound holds uniformly over the boundary-parameter torus.

At $s=1$, Lemma~\ref{ml:lem:zero} and
\eqref{ml:eq:strong-comparison} identify that subspace with the
common three-dimensional kernel in $N=N_0$.  Use the restriction of
the projection $P$ from the energy proof. The second path is
\begin{equation}\label{vb:second-path}
 B_u=\left.\varepsilon_\kappa^{-1}
 \bigl(K_{L,\kappa}+u(H_{L,\kappa}-K_{L,\kappa})\bigr)
 \right|_{N=N_0},
 \qquad 0\leq u\leq1.
\end{equation}
By \eqref{ml:eq:comparison-order} and \eqref{ml:eq:strong-comparison}, $B_uP=0$ and $B_u\geq g_t(I-P)$ for every $u$. Thus both paths preserve an isolated subspace of rank three over the entire parameter torus.

\smallskip
\noindent\emph{Step 3: Transfer the Chern class.}
Proposition~\ref{sp:continuation}, applied to these two paths, identifies the full lowest-eigenspace bundle with the periodic-state bundle of the word $100$. Writing $\mathcal G_{L,\kappa}(\theta)= \ker(\widehat H_{L,\kappa}(\theta)|_{\mathcal H_{L,N_0}(\theta)})$, Proposition~\ref{mg:primitive-bundle}, with $p=1$ and $q=3$, gives
\begin{equation}\label{vb:full-chern-number}
 c_1(\mathcal G_{L,\kappa})=\mathfrak u.
\end{equation}
Lemma~\ref{vb:transport} also transports all energy and gap bounds, and the absence of lower-density plateaux, to every boundary parameter.
\end{proof}

\subsection{Proofs of the auxiliary results}

\begin{proof}[Proof of Lemma~\ref{ml:lem:zero}]
The construction in \cite{HaldaneRezayi1985} supplies three linearly
independent antisymmetric sections in
$\bigwedge^{N_0}\mathcal E_L(0)$ when $L=3N_0$. Their holomorphic
representatives vanish to order at least three whenever two
coordinates coincide. This holds for every rectangular aspect
ratio. The standard representation of the $V_1$ operator by the
squared first relative derivative at coincidence
\cite[\S1.1]{WY} gives zero quadratic form on these sections.
Since $H_{L,\kappa}=\sum_RQ_R^*Q_R$, they belong to every
$\ker Q_R$, proving the first assertion and, by
\eqref{ml:eq:comparison-order}, also belong to $\ker K_{L,\kappa}$.

To decrease the exterior degree, note that $[Q_R,c_j]=0$:
each $Q_R$ is a sum of products of two exterior contractions.
If $\psi\ne0$ has degree $N>0$ and belongs to $\ker H_{L,\kappa}$,
then
\[
 \sum_j\|c_j\psi\|^2=N\|\psi\|^2>0.
\]
Some $c_j\psi$ is therefore nonzero. It has degree $N-1$ and
$Q_Rc_j\psi=c_jQ_R\psi=0$ for every $R$.
Iteration gives a nonzero kernel vector in every degree below $N_0$.
The operator is nonnegative, so its minimum eigenvalue in each of
these degrees is zero.
\end{proof}

We use the following specialization of the stability theorem in
\cite[Theorems 6.8 and 6.12, with the translation-averaging
modification in \S8.4.2]{NSY}. We state the needed hypotheses
and conclusion explicitly, and take this proposition as an
external result.

\begin{Pro}[Uniform stability on periodic chains]
\label{ml:prop:nsy-stability}
Consider rings of lengths divisible by a fixed integer $q$,
with a fixed finite-dimensional space at each site. Let $H_L^0$
be the periodic restrictions of a fixed translation-invariant,
finite-range interaction consisting of bounded positive operators
with a nonzero common kernel. Assume:
\begin{enumerate}
\item The ring kernel has a fixed dimension $d$. The ring operators
and their restrictions to intervals have a gap at least
$\gamma_0>0$ above their kernels whenever they have nonzero spectrum.
\item There are constants $b>0$ and $0<c<1$, independent of $L$,
such that for every sufficiently large interval $I$ with
$|I|\le cL$ and every operator $A$ supported at distance at least
$b$ from its boundary,
\[
 P_I\left(\frac1q\sum_{r=0}^{q-1}\mathcal T^r(A)\right)P_I
 =\omega^{\rm av}(A)P_I.
\]
Here $P_I$ is the local kernel projection, $\mathcal T$ is
translation by one site, all displayed translates are supported
in $I$, and $\omega^{\rm av}$ is a translation-invariant state
independent of $I,L$.
\end{enumerate}
Then there is $\epsilon>0$, independent of $L$, such that every
self-adjoint translation-invariant interaction $\Phi_L$ with
$\sup_L\|\Phi_L\|_*<\epsilon$, in the norm
\eqref{ml:eq:interaction-norm}, has the following property:
for all $s\in[0,1]$, the lowest $d$ eigenvalues of
\[
 H_L^0+s\sum_X\Phi_L(X)
\]
are separated from the remaining spectrum by at least
$\gamma_0/2$. Their spectral projection is continuous in $s$
and has rank $d$. The eigenvalues within this cluster need
not be equal.
\end{Pro}

The norm used here implies the exponential interaction bounds
in the cited theorem. Indeed, distributing each nonempty term
equally over its sites gives terms supported in radius-$m$
balls,
\[
 \Psi_L(z,m)=
 \sum_{\substack{X\ni z\\\max\{2,\diam(X)\}=m}}
 \frac{\Phi_L(X)}{|X|},
 \qquad
 \|\Psi_L(z,m)\|\le e^2\|\Phi_L\|_*e^{-m}.
\]
If two sites are at distance $r$, the weighted sum of norms
over balls containing both is bounded by
\[
 e^2\|\Phi_L\|_*
 \sum_{m\ge\max\{2,\lceil r/2\rceil\}}(2m+1)^2e^{-m}
 \le C\|\Phi_L\|_*\frac{e^{-r/4}}{(1+r)^2},
\]
with $C$ independent of $L$. The pointwise bound also has this
exponential form. Thus the proposition uses the same perturbation
class as the reference. Condition (2) has zero error beyond a
fixed buffer and holds on intervals of length proportional to
$L$; together with the uniform local gaps, this is the
specialization of its uniformity hypotheses needed here.

\begin{proof}[Proof of Lemma~\ref{ml:lem:stability}]
We apply Proposition~\ref{ml:prop:nsy-stability} to
$H_L^0=\mathcal L_{L,t}$ with $q=d=3$.

\smallskip
\noindent\emph{Local kernels and gaps.}
The decomposition \eqref{ml:eq:reference} consists of commuting
positive terms, since
\[
 h_t(0)=t/3,\quad h_t(1)=0,\quad
 h_t(2)=1-t/3,\quad h_t(3)=3-2t/3.
\]
Every positive eigenvalue of a local term is therefore at least
\begin{equation}\label{ml:eq:reference-gap}
 \delta_t=\min\{t/3,1-t/3\}>0.
\end{equation}
The same lower bound holds for any sum of these terms.
A binary word has zero energy precisely when
$n_j+n_{j+1}+n_{j+2}=1$ for every included triple.
Subtracting consecutive equations gives $n_j=n_{j+3}$.
Thus on any interval of at least three sites the kernel is
spanned by the restrictions of the three translates of the
periodic word $100$; on a ring of length divisible by three
these give exactly three kernel vectors. This proves condition (1).

\smallskip
\noindent\emph{Scalar compression after translation averaging.}
Let $w_0,w_1,w_2$ denote these three local words on $I$, and
suppose that $A$ is supported at least six sites from its
boundary. For $a\ne a'$, the words $w_a,w_{a'}$ differ on a
site outside the support of each $\mathcal T^r(A)$, $0\le r\le2$.
Hence all off-diagonal matrix elements of the compressed
average vanish. Translation cycles the three words, so its
three diagonal entries are all equal to
\[
 \omega^{\rm av}(A):=
 \frac13\sum_{a=0}^2\langle w_a,Aw_a\rangle.
\]
This value is independent of the containing interval and is
translation invariant. Consequently condition (2) holds, for
example with buffer six and $|I|\le L/2$.

The perturbations already satisfy the symmetry and norm
hypotheses of Proposition~\ref{ml:prop:nsy-stability}.
It gives the assertion with $g_t=\delta_t/2$ for all
sufficiently large rings. To include the finitely many smaller
ones, let $L_*$ bound their lengths and put
$\eta=\sup_L\|\Phi_L\|_*$. Since
$\|s\sum_X\Phi_L(X)\|\le L\eta$, the variational principle
bounds their perturbed gap below by $\delta_t-2L_*\eta$.
Reducing $\epsilon_t$ to at most $\delta_t/(4L_*)$ therefore
gives the same gap for every $L\in6\mathbb N$.
The isolated spectral projection is continuous and has rank three.
\end{proof}

\begin{proof}[Proof of Proposition~\ref{ml:prop:limits}]
Fix $\kappa\ge3$ and set $h(x)=xe^{-\kappa^2x^2}$.
For $1/2\le x\le L/2$, monotonicity of $h$ and
$jL-x\ge(2j-1)x$ give
\begin{equation}\label{ml:eq:periodic-gaussian-bounds}
 |f_{L,\kappa}(x)|
 \le h(x)+2\sum_{j\ge1}h(jL-x)
 \le h(x)\left[1+2\sum_{j\ge1}(2j-1)
 e^{-\kappa^2x^2((2j-1)^2-1)}\right]\le4h(x).
\end{equation}
Indeed the bracket is at most
$3+2\sum_{j\ge1}(2j+1)e^{-18j}<4$.
Thus $V_{L,\kappa}(r)\le32V_\kappa(r)$ for $1\le r\le L/2$,
while $V_{L,\kappa}(r)\to V_\kappa(r)$ at each fixed $r$.
For $r\le L/2$, the periodization in \eqref{pc:periodization}
is at most $2\sum_{j\ge0}V_\kappa((2j+1)r)\le3V_\kappa(r)$.
Dominated convergence therefore shows that the
sum of absolute differences between $V_{L,\kappa}$ and the periodized
$V_\kappa$ over distances $1\le r\le L/2$ tends to zero.
The difference of the corresponding diagonal operators divided by
$L$ is bounded by that sum. Their minimum eigenvalues per site at fixed degree
therefore have the same limit $\beta_{V_\kappa}$, as do their
variational functions $\alpha$. The potential $V_\kappa$ satisfies \eqref{item:U1}--\eqref{item:U3}: Gaussian decay gives the moment bounds, and direct differentiation gives
\[
 V_\kappa''(x)=\frac14e^{-\kappa^2x^2/2}
 (\kappa^4x^4-5\kappa^2x^2+2)>0\qquad(x\ge1).
\]
Proposition~\ref{pro:beta-regularity} now gives the widths. The derivative formula \eqref{eq:beta-onesided} gives \eqref{ml:eq:classical-width}; the first terms of its $q=3$ series give \eqref{ml:eq:classical-one-third}.

It remains to prove convergence for the full operator, for any
fixed $\kappa>0$ and $c\in\mathbb R$. The idea is that cutting
a finite-range interaction into two intervals changes its minimum
by a bounded amount; the Gaussian tail then allows us to remove
the finite-range restriction. Throughout, we use the variational
estimate $|\min\spec A-\min\spec B|\le\|A-B\|$.

For an integer $J\ge1$, let $H_{L,\kappa}^{(J)}$ retain only
$0<x\le J$ in the pair sums of Lemma~\ref{ml:lem:tail}, before
taking their squares. The pair matrices have norm at most one,
and the common phase cancels in each square. Thus
$\|A^*A-B^*B\|\le(\|A\|+\|B\|)\|A-B\|$, summed over the
$2L$ centers, gives
\[
 \frac1L\|H_{L,\kappa}-H_{L,\kappa}^{(J)}\|
 \le4\left(\sum_{x\in\frac12\mathbb N}xe^{-\kappa^2x^2}\right)
       \left(\sum_{\substack{x\in\frac12\mathbb N\\x>J}}
                        xe^{-\kappa^2x^2}\right)
 =:\varepsilon_J\longrightarrow0.
\]
This bound is uniform in $L$.

For fixed $J$, the truncated interaction has range at most $2J$
and coefficients independent of the interval length. Let $e_J(n)$
be its minimum eigenvalue on an open interval of $n$ sites,
including $-c\Num$ and minimizing over all exterior degrees.
Removing the terms crossing a cut leaves two independent
operators, whose minima add. The removed terms have uniformly
bounded total norm, since only sites within distance $2J$ of
the cut contribute. Therefore
\[
 |e_J(n+m)-e_J(n)-e_J(m)|\le C_J
\]
for a constant independent of $n,m$.
For any common multiple $\ell$ of $n,m$, iteration gives
\[
 \left|\frac{e_J(n)}n-\frac{e_J(m)}m\right|
 \le\left|\frac{e_J(n)}n-\frac{e_J(\ell)}\ell\right|
   +\left|\frac{e_J(\ell)}\ell-\frac{e_J(m)}m\right|
 \le C_J\left(\frac1n+\frac1m\right).
\]
Hence $e_J(n)/n$ converges. Closing the interval into a ring
changes the operator by a bounded norm, depending on $J$ but
not on $n$, so the truncated ring minimum per site converges too.
Here it suffices to take $L>4J$, so the truncated local terms
are the same as on the infinite chain.

The full and truncated ring minima per site differ by at most
$\varepsilon_J$. Consequently the limsup and liminf of the full
sequence differ by at most $2\varepsilon_J$ for every $J$.
Letting $J\to\infty$ proves existence of
\[
 \alpha_\kappa(c)=\lim_{\substack{L\to\infty\\6\mid L}}
 \alpha_{H_{L,\kappa},L}(c).
\]
The same proof applies to the diagonal operators: taking the
diagonal does not increase the operator norm and preserves local
supports. For $\kappa\ge3$, the first part of the proof identifies
this diagonal limit with the classical variational function
for $V_\kappa$.

\end{proof}

\begin{proof}[Proof of Lemma~\ref{vb:transport}]
The unitary transport of Lemma~\ref{mg:magnetic-transport} acts
on each variable by the same spatial translation and scalar phase factor.
The projected $V_1$ interaction is covariant under this operation:
in its standard first-relative-derivative contact form
\cite[\S1.1]{WY}, coincidence sets are translated together and
the derivatives of the common scalar phase cancel in
$\partial_{x_i}-\partial_{x_j}$. Thus the quadratic form and its
fixed energy normalization are preserved, proving the first identity
in \eqref{vb:unitary-equivalence}. Exterior powers preserve exterior
degree, proving the second identity. Finally,
\eqref{vb:transport-orbitals} identifies the coefficient matrix in
the rephased frame.
\end{proof}

\section{The Coulomb operator: spectral stability and density plateaux}\label{co:section}

\subsection{The complete operator and three estimates}
\label{co:model}

We use the complete projected Coulomb matrix elements of
\cite[Supplementary Material, \S I, Eq.~(2)]{RotondoEtAl2016},
written in our conventions as \eqref{co:representation}.
The geometry, basis, and scalar normalization are fixed in
Section \ref{eq:coulomb stab}. The change of basis under exchange of the two sides
is described in \cite[Supplementary Material, \S III]{RotondoEtAl2016}.
The model formula is an input; we do not repeat its derivation.
We require the complete coefficients, before the sum-to-integral
and thin-torus approximations in that reference.

Retain the decomposition $H^{\mathrm C}=D^{\mathrm C}+T^{\mathrm C}$
and the kernel $W$ from Section~\ref{eq:coulomb stab}. Put
$C_\kappa=g\kappa/\pi$, and use $J_L$ from
Lemma~\ref{pc:cyclic-kernel}.
The proof uses only the following three estimates on this
decomposition. Explicit coefficients and their estimates are
collected in Appendix~\ref{co:analysis}.

\begin{Lem}[Three estimates for the complete Coulomb operator]
\label{co:bounds}
Let $\kappa\ge4$ and $L\ge6$, and put
\[
 \xi_\kappa=\frac\pi\kappa
 (32e^{-\kappa^2/2}+12e^{-18\kappa^2})<1/100.
\]
\begin{enumerate}
\item The symmetric coefficient $W(L-r)=W(r)$ satisfies
\[
 \Delta^2W(r)\ge\frac{(1-\xi_\kappa)C_\kappa}{r^2},
 \qquad 2\le r\le J_L.
\]
\item The off-diagonal part satisfies
\[
 \|T^{\mathrm C}\|\le\eta_{L,\kappa}
 :=g\kappa L e^{-2\kappa^2}.
\]
\item Let $P_N$ project onto the span of the balanced words of
length $L$ with $N$ occupied indices. Then
\begin{equation}\label{co:compressions}
 P_NT^{\mathrm C}P_N=0,\qquad
 P_N(T^{\mathrm C})^2P_N=\varsigma_N^2P_N,
 \qquad \varsigma_N^2\le4L(g\kappa)^2e^{-4\kappa^2}.
\end{equation}
In particular, $\|T^{\mathrm C}P_N\|^2=\varsigma_N^2$.
\end{enumerate}
\end{Lem}

The first estimate yields the diagonal spectral gap; the second
keeps the complementary spectrum separated; the third controls
the change of the minimum eigenvalue to second order.
The coefficient matrix is independent of the boundary parameters
under Lemma~\ref{mg:magnetic-transport}, so these bounds hold
uniformly over their torus.

\subsection{Minimizing binary words and diagonal spectral gaps}
\label{co:spectral}

We retain the notation $E^{\mathrm C}(N)$ and $E_D^{\mathrm C}(N)$
from Section \ref{eq:coulomb stab} for the lowest energies of $H^{\mathrm C}$ and
$D^{\mathrm C}$ on $\mathcal H_{L,N}$.
Unless otherwise specified, the spectral estimates and comparisons between degrees below assume
$L\ge6$ and $\kappa\ge4$.
To apply Lemma~\ref{pc:cyclic-kernel} to the exact diagonal, put
\begin{align}
 b_\ell&=\Delta^2W(\ell)\quad(2\le\ell<J_L),&
 b_{J_L}&=W(J_L-1)-W(J_L),\nonumber\\
 W^\circ(r)&=\sum_{\ell=2}^{J_L}b_\ell(\ell-r)_+,&
 d_*&=\frac{(1-\xi_\kappa)C_\kappa}{2J_L^2}.
 \label{co:hinge}
\end{align}
Symmetry gives $b_{J_L}=\tfrac12\Delta^2W(J_L)$ for even $L$
and $b_{J_L}=\Delta^2W(J_L)$ for odd $L$.  Thus
\begin{equation}\label{co:hinge-bounds}
 b_\ell\ge\frac{(1-\xi_\kappa)C_\kappa}{\ell^2}
 \quad(\ell<J_L),\qquad
 b_{J_L}\ge d_*,\qquad
 d_*\ge\frac{2(1-\xi_\kappa)C_\kappa}{L^2}.
\end{equation}
Finite summation of the second differences, with the central value
and first difference fixed, yields
$W(r)=W(J_L)+W^\circ(r)+W^\circ(L-r)$ for $1\le r<L$.
This is precisely the hinge representation of
Lemma~\ref{pc:cyclic-kernel}.

\begin{Pro}[Exact minima and two diagonal gap estimates]
\label{co:classical}
For every $0\le N\le L$,
\begin{equation}\label{co:diagonal-energy}
 E_D^{\mathrm C}(N)=W(J_L)\binom N2+
 L\sum_{\ell=2}^{J_L}b_\ell F(\ell N/L).
\end{equation}
If $0<N<L$ and $N/L=p/q$ is reduced, its lowest eigenspace consists
of exactly the $q$ balanced binary-word basis vectors of
$X_{p/q}^{\rm per}$ restricted to the ring.  Use the projection $P_N$ from Lemma~\ref{co:bounds},
and put $Q_N=1-P_N$ on $\mathcal H_{L,N}$.
Then
\begin{equation}\label{co:diagonal-all-gap}
 D^{\mathrm C}|_{\mathcal H_{L,N}}-E_D^{\mathrm C}(N)\ge d_*Q_N.
\end{equation}
The same inequality holds at $N=0,L$ with $Q_N=0$.
If $q<J_L$, the stronger bound holds with
\begin{equation}\label{co:diagonal-q-gap}
 \delta_q=\frac{(1-\xi_\kappa)C_\kappa}{q^2}
\end{equation}
in place of $d_*$.
\end{Pro}

\begin{proof}
The energy formula follows from Lemma~\ref{pc:cyclic-kernel}.
Its block identity \eqref{pc:cyclic-block-identity} gives
\begin{equation}\label{co:integer-defects}
 D^{\mathrm C}(\omega)-E_D^{\mathrm C}(N)
 =\sum_{\ell=2}^{J_L}b_\ell\mathscr D_\ell(\omega).
\end{equation}
The defects are those of \eqref{pc:block-defect}. They are
nonnegative integers: if
$\ell N=Lk+r$, $0\le r<L$, then
$LF(\ell N/L)=(L-r)\binom k2+r\binom{k+1}2$;
Lemma~\ref{lem:integer-Jensen} gives nonnegativity and identifies
zero defect with balanced length-$\ell$ counts.

If all defects for $2\le\ell\le J_L$ vanish, interval counts
at lengths $J_L<\ell<L$ are balanced as well, since they are
$N$ minus the counts in complementary intervals.  Lengths one
and $L$ are automatic.  Periodicity then gives the balanced
condition at every length, and Lemma~\ref{lem:rational-balanced}
identifies the configuration with one of the $q$ minimizing words.
Every other configuration has some integer defect at least one;
\eqref{co:hinge-bounds} proves \eqref{co:diagonal-all-gap}.

For the sharper bound, vanishing defects at lengths $2,\ldots,q$
give $S_q(j)=p$ and hence $\omega_{j+q}=\omega_j$.
The shorter balanced counts then imply all the balanced counts,
as in the proof of \eqref{pc:uniform-gap}.  A nonminimizing word thus has
a positive defect at some $\ell\le q<J_L$, where
$b_\ell\ge\delta_q$.  This proves
\eqref{co:diagonal-q-gap}.
\end{proof}

\subsection{Shift and diagonal symmetries at fixed exterior degree}

Restrict $\mathsf S,\mathsf C$ from \eqref{eq:exterior-clockshift}
to $\mathcal H_{L,N}$. For $N/L=p/q$,
\begin{equation}\label{co:weyl}
 \mathsf C\mathsf S=e^{2\pi ip/q}\mathsf S\mathsf C.
\end{equation}
The Fourier-series representation \eqref{co:representation} shows that
$H^{\mathrm C}$ commutes with both operators: conjugation multiplies
each operator $\mathcal R_{mn}$ by a scalar phase, which cancels in its
squared modulus.  Its diagonal and off-diagonal parts also commute
with them, since the shift permutes the ordered exterior-product basis and the
clock acts diagonally in that basis.

\begin{Lem}[Eigenvalue multiplicities]
\label{co:magnetic}\label{co:variance}
At reduced density $p/q$, every eigenspace of
$D^{\mathrm C}+uT^{\mathrm C}$ has dimension divisible by $q$.
\end{Lem}

\begin{proof}
On any invariant eigenspace $\mathcal E$, taking determinants in
\eqref{co:weyl} gives $e^{2\pi ip\dim\mathcal E/q}=1$.
Since $p$ and $q$ are coprime, $q\mid\dim\mathcal E$.
\end{proof}

The equality $PT^{\mathrm C}P=0$ removes the first-order change
of the minimum eigenvalue. The quantity $\|T^{\mathrm C}P\|^2$
measures coupling to the orthogonal complement. Solving the
eigenvalue equation on that complement and substituting into the
$P$ component gives the following second-order estimate.

\begin{Lem}[Schur-complement estimate for a lowest eigenvalue of multiplicity $q$]
\label{co:feshbach}
Let $D,T$ be self-adjoint operators on a finite-dimensional space,
commuting with unitaries
$\mathsf S,\mathsf C$ satisfying
$\mathsf C\mathsf S=\omega\mathsf S\mathsf C$, where $\omega$
has order $q$.  Suppose $D$ has a lowest-eigenspace projection $P$ of rank
$q$, with minimum eigenvalue $E_0$, and
\[
 D-E_0\ge\delta(1-P),\qquad PTP=0,
 \qquad \|T\|\le\eta<\delta.
\]
Put $Q=1-P$ and $\varsigma=\|TP\|$.
For $0\le u\le1$, $D+uT$ has an exactly $q$-fold minimum eigenvalue
$E_0-e(u)$, and
\begin{align}
 0\le e(u)&\le\frac{u^2\varsigma^2}{\delta-u\eta+e(u)},
 \label{co:schur-energy}\\
 \gap(D+uT)&\ge\delta-u\eta+e(u).\nonumber
\end{align}
For $Q=0$, one has $T=0$, $e=\varsigma=0$, and no eigenvalues above the minimum.
\end{Lem}

\begin{proof}
Subtract $E_0$ and assume $Q\ne0$. The min--max principle and
$PTP=0$ give $\lambda_q(D+uT)\le0$, whereas
$\lambda_{q+1}(D+uT)\ge\delta-u\eta>0$.
Every eigenvalue multiplicity is divisible by $q$, by the
determinant argument above.  The lowest $q$ eigenvalues therefore
coincide at $-e(u)\le0$, proving the gap assertion.

Let $D_Q=Q(D-E_0)Q$ and $A_u=D_Q+uQTQ$ on $Q\mathcal H$.
Then $A_u\ge\delta-u\eta$.  For a lowest-eigenvalue vector
$\psi=p+r$, $p=P\psi$, $r=Q\psi$, the $Q$ equation is
$r=-u(A_u+e(u))^{-1}QTPp$.  The map $\psi\mapsto p$ is
injective and hence an isomorphism onto $P\mathcal H$.
Substitution into the $P$ equation gives
\begin{equation}\label{co:schur-identity}
 e(u)P=u^2PTQ(A_u+e(u))^{-1}QTP.
\end{equation}
Taking norms and using $QTP=TP$ proves the energy bound.
\end{proof}

\subsection{Proof of the main Coulomb theorem}

\begin{proof}[Proof of Theorem~\ref{intro:coulomb}]
We use Proposition~\ref{co:classical} for the diagonal gaps,
Lemma~\ref{co:magnetic} for the multiplicities and
Lemma~\ref{co:bounds} for the bound
$\|T^{\mathrm C}P_N\|^2\le4L(g\kappa)^2e^{-4\kappa^2}$,
and Lemma~\ref{co:feshbach} for the perturbation.

For the first assertion, $q\le Q$ and $L\ge3Q$ imply $q<J_L$.
The size condition gives
\[
 \eta_{L,\kappa}\le\frac{C_\kappa}{4Q^2},\qquad
 \delta_q-\eta_{L,\kappa}
 \ge\frac{C_\kappa}{2q^2}.
\]
Lemma~\ref{co:feshbach} therefore gives an exactly $q$-dimensional
lowest eigenspace throughout $D^{\mathrm C}+uT^{\mathrm C}$,
$0\le u\le1$, with the stated gap. Its energy lowering is at most
\[
 \frac{4L(g\kappa)^2e^{-4\kappa^2}}{C_\kappa/(2q^2)}
 =8\pi g\kappa Lq^2e^{-4\kappa^2}.
\]
The coefficient matrices are independent of the boundary parameters
and commute with the transition maps, by Subsection~\ref{co:model}
and Lemma~\ref{co:magnetic}. The gapped path thus identifies the
lowest-eigenspace bundle with that of the balanced reference words.
Propositions~\ref{sp:continuation} and~\ref{mg:primitive-bundle}
give rank $q$ and Chern number $p$, uniformly in the boundary parameters.

For the second assertion, the stronger size condition yields
$\eta_{L,\kappa}\le C_\kappa/(4L^2)$. By
\eqref{co:hinge-bounds} and $\xi_\kappa<1/100$,
\[
 d_*-\eta_{L,\kappa}>C_\kappa/L^2.
\]
The same Schur-complement bound, now using the gap $d_*$ valid
in every degree, gives
\[
 E^{\mathrm C}(N)=E_D^{\mathrm C}(N)-e_N,
 \qquad 0\le e_N\le
 \frac{4L(g\kappa)^2e^{-4\kappa^2}}{C_\kappa/L^2}=\zeta.
\]
In degrees $N=0,1,L-1,L$, the non-diagonal terms vanish because
they require two occupied and two empty sites, so $e_N=0$.

It remains to compare different degrees. Apply the secant-slope
definitions of $c_{A,\pm}(N)$ in Section~\ref{sec:classical-quantum}
to $A=D^{\mathrm C}$ and $A=H^{\mathrm C}$.
Every secant slope changes by at most $\zeta$, since
$|e_{N'}-e_N|\le\zeta$ and $|N'-N|\ge1$. The same bound holds after
taking maxima or minima. Hence each endpoint moves by at most
$\zeta$, and every classical plateau of width greater than
$2\zeta$ persists, with width at least its original width minus
$2\zeta$. Adding the same $\chi(N)$ to both energy functions leaves
$e_N$ unchanged, proving the remark following the theorem as well.
\end{proof}

\appendix

\section{The invariant-measure formulation and Aubry sets}
\label{app:aubry-mather}

We extend the block identity to invariant measures and use periodic closing to classify Aubry sets.

\subsection{Definitions and the general variational statement}
Use the local energy $A_c$ from \eqref{eq:periodic-energy}.
Retain $\Sigma_2$, $\sigma$, $\mathcal M_\sigma$, $\rho(\mu)$,
and $e_U(\mu)$ from Section~\ref{sec:classical structure}.
In terms of the local energy, $e_U(\mu)=\int_{\Sigma_2}A_0\,d\mu$.
The local energy is continuous by summability; the ergodic theorem and the boundary estimate identify its invariant average with the mean interval energy.

An invariant probability measure minimizing $e_U$ at a prescribed density is called a \emph{Mather measure} for this symbolic variational problem. Its parameter-dependent version minimizes $e_U-c\rho$. Write $\mathscr M_\rho$ and $\mathscr M(c)$ for the unions of supports of the corresponding minimizing measures. The Sturmian subshift $X_\rho^{\rm St}$ is defined in
\eqref{eq:strumian}. It is minimal and uniquely ergodic; conversely, every aperiodic balanced binary sequence of density $\rho$ belongs to it \cite{MorseHedlund1940,lothaire2002algebraic}.

We now prove completeness in the sense specified in the introduction.

\begin{The}[Invariant minimizing measures and the complete staircase]
\label{app:general-variational}\label{thm:Mather-measure-classification}
Assume \eqref{item:U1}--\eqref{item:U3}. For the variational problem
\eqref{eq:beta-def}, the block formula \eqref{eq:beta-formula} holds
at every density. The minimizing measure is unique: it is the uniform orbit measure on $X_{p/q}^{\rm per}$ when $\rho=p/q\in(0,1)$, the unique invariant measure on $X_\rho^{\rm St}$ when $\rho$ is irrational, and the point mass on the constant word at $\rho=0,1$. Moreover, $\alpha_U\in C^1(\mathbb R)$ and its derivative is a continuous nondecreasing function whose rational plateau interiors cover $(0,2\sum_{r\ge1}U(r))$ up to a Lebesgue null set. Thus $\alpha'_U$ is a complete devil's staircase on the transition interval, with derivative zero almost everywhere, and
\[
 \alpha'_U(c)=0\quad(c\le0),\qquad
 \alpha'_U(c)=1\quad\left(c\ge2\sum_{r\ge1}U(r)\right).
\]
\end{The}

\begin{Pro}[Invariant block averages]\label{app:invariant-block-energy}
For $\mu\in\mathcal M_\sigma$,
\begin{equation}\label{eq:block-average}
 \int_{\Sigma_2} B_n\,d\mu=\int_{\Sigma_2} f\bigl(S_n(0;\omega)\bigr)\,d\mu,
\end{equation}
and consequently
\begin{equation}\label{eq:E-block}
 e_U(\mu)=\sum_{n\ge2}\Delta^2U(n)
 \int_{\Sigma_2} f\bigl(S_n(0;\omega)\bigr)\,d\mu.
\end{equation}
\end{Pro}
\begin{proof}
Expand $f(S_n(0;\omega))$ as the sum of pairs in a length-$n$ block. Shift invariance gives
\[
 \int_{\Sigma_2} f\bigl(S_n(0;\omega)\bigr)\,d\mu
 =\sum_{r=1}^{n-1}(n-r)\int_{\Sigma_2}\omega_0\omega_r\,d\mu
 =\int_{\Sigma_2} B_n\,d\mu.
\]
Integrating the uniformly convergent identity \eqref{eq:A0-Bn} proves \eqref{eq:E-block}.
\end{proof}

\begin{Pro}[Distributional derivative]\label{app:interpolation}
On $(0,1)$,
\begin{equation}\label{app:beta-measure}
 D^2\beta_U=\sum_{n\ge2}n\Delta^2U(n)
 \sum_{j=1}^{n-1}\delta_{j/n}.
\end{equation}
Hence $\beta_U$ is differentiable exactly at the irrational interior densities.
\end{Pro}
\begin{proof}
The second derivative of $F(n\rho)$ is
$n\sum_{j=1}^{n-1}\delta_{j/n}$ (in the sense of distribution).
Their weighted total mass is finite by \eqref{eq:sumU-an}, so the series may be differentiated distributionally.
Its atoms are exactly the rational densities; alternatively, use \eqref{eq:beta-onesided}.
\end{proof}

\begin{proof}[Proof of Theorem~\ref{app:general-variational}]
For an invariant measure of density $\rho$, the mean of $S_n(0;\omega)$ is $n\rho$. Proposition~\ref{app:invariant-block-energy} and Lemma~\ref{lem:integer-Jensen} therefore give $e_U(\mu)\ge\beta_U(\rho)$. At rational densities equality is attained by the periodic orbit measure from Theorem~\ref{intro:variational}. For arbitrary $\rho$, choose rational $\rho_k\to\rho$.
A weakly convergent subsequence of the orbit measures has an invariant limit. Continuity of $A_0$, $\omega_0$, and $\beta_U$ gives its density $\rho$ and energy $\beta_U(\rho)$, proving \eqref{eq:beta-def}.

For any minimizing measure, strict positivity of $\Delta^2U(n)$ implies
\begin{equation}\label{eq:all-Jensen-equalities}
 \int_{\Sigma_2} f\bigl(S_n(0;\omega)\bigr)\,d\mu=F(n\rho)\quad(n\ge2).
\end{equation}
The equality cases in Jensen's inequality, shift invariance, and a countable intersection give
\begin{equation}\label{eq:balanced-block-condition}
 S_n(j;\omega)\in\{\lfloor n\rho\rfloor,\lceil n\rho\rceil\}
 \quad(j\in\mathbb Z,\ n\ge1)
\end{equation}
on a set of full measure.
These closed block constraints therefore hold on the support.
At rational density, Lemma~\ref{lem:rational-balanced} puts the support in the single finite orbit $X_{p/q}^{\rm per}$, which carries only its uniform invariant measure. At irrational density the sequences in \eqref{eq:balanced-block-condition} are balanced and aperiodic; the Sturmian characterization and unique ergodicity give the asserted unique measure. The two endpoint cases follow directly from $\omega_j\in\{0,1\}$ and invariance.

The strict convexity of $\beta_U$ by Proposition \ref{pro:beta-regularity} gives a unique maximizer in \eqref{eq:alpha}; compactness makes this maximizer continuous in $c$. The elementary derivative formula for a maximum of affine functions then gives $\alpha'_U(c)=\rho(c)$, the unique selected density. It is nondecreasing by convexity. The subgradient criterion shows that its constant intervals are exactly those computed in Section~\ref{sec:minimizing-sets}. Grouping the nonnegative slope jumps by their reduced rational location gives
\[
 \sum_{\substack{p/q\in(0,1)\cap\mathbb Q\\\gcd(p,q)=1}}W_{U,p/q}
 =\sum_{q=2}^\infty\varphi(q)\sum_{l=1}^\infty lq\Delta^2U(lq)=\sum_{n=2}^\infty\sum_{q\mid n,q\geq 2}n\Delta^2U(n)\varphi(q)=\sum_{n\ge2}n(n-1)\Delta^2U(n)
 =2\sum_{r\ge1}U(r).
\]
Here we use the notation of Euler's totient function $\varphi(q)$ and the fact $\sum_{q\mid n}\varphi(q)=n$ and $\varphi(1)=1$. 
Their interiors are disjoint, so they fill the transition interval up to a null set. Hence the derivative of $\alpha'_U$ is zero almost everywhere. The endpoint derivatives of $\beta_U$ in Proposition \ref{pro:beta-regularity} give the remaining claims. 
\end{proof}

The preceding theorem identifies the original invariant formulation as
\[
 \alpha_U(c)=-\min_{\mu\in\mathcal M_\sigma}\int_{\Sigma_2} A_c\,d\mu,
 \qquad
 \mathscr M(c)=\mathscr M_{\alpha'_U(c)},\qquad
 \mathscr M_\rho=\bigcup_{c\in D^-\beta_U(\rho)}\mathscr M(c).
\]
Thus density is the rotation parameter and $\beta_U$ the minimum average action. The symbolic Peierls barrier \cite{GaribaldiLopes2008} also detects nonrecurrent minimizing configurations.

\subsection{Aubry sets defined by return actions}
An isolated transition between periodic tails need not support a new invariant measure.
We detect such configurations by approximate returns of vanishing normalized action.
Closing these returns requires control of interactions crossing the cut, so we impose one additional moment:

We assume here, in addition to \eqref{item:U1}--\eqref{item:U3}, that
\begin{enumerate}\setcounter{enumi}{3}\renewcommand{\theenumi}{U\arabic{enumi}}
\item\label{item:U4} $\sum_{r\ge1}rU(r)<\infty$.
\end{enumerate}
The discrete convex decomposition gives
\begin{equation}\label{app:first-moment}
 \sum_{r\ge1}rU(r)
 =\frac16\sum_{n\ge2}n(n-1)(n+1)\Delta^2U(n).
\end{equation}
Thus \eqref{item:U4} is equivalent to $\sum_{n\ge2}n^3\Delta^2U(n)<\infty$.

Let $\Sigma_2^+=\{0,1\}^{\mathbb N_0}$ with $d(\mathbf x,\mathbf y)=\sum_{j\ge0}2^{-j-1}|x_j-y_j|$. The action $A_c$ of \eqref{eq:periodic-energy} depends only on these nonnegative coordinates. Set
\begin{equation}\label{eq:Abar}
 \bar A_c:=A_c+\alpha_U(c),\qquad
 \mathcal S_N\bar A_c(\mathbf z):=
 \sum_{j=0}^{N-1}\bar A_c(\sigma^j\mathbf z).
\end{equation}
For $\mathbf x,\mathbf y\in\Sigma_2^+$, the Peierls barrier is
\begin{equation}\label{eq:Peierls-barrier}
 h_c(\mathbf x,\mathbf y):=
 \lim_{\varepsilon\to 0^+}\liminf_{N\to\infty}
 \inf_{\substack{d(\mathbf z,\mathbf x)<\varepsilon\\
 d(\sigma^N\mathbf z,\mathbf y)<\varepsilon}}
 \mathcal S_N\bar A_c(\mathbf z).
\end{equation}
Define
\begin{equation}\label{eq:projected-Aubry}
 \mathscr A^+(c)=\{\mathbf x:h_c(\mathbf x,\mathbf x)=0\}.
\end{equation}
Writing $\pi^+\omega=(\omega_0,\omega_1,\ldots)$, its two-sided version is
\begin{equation}\label{eq:natural-Aubry}
 \mathscr A(c)=\{\omega\in\Sigma_2:
 \pi^+\sigma^j\omega\in\mathscr A^+(c)\ \hbox{for every }j\in\mathbb Z\}.
\end{equation}
Use the endpoint classes $\mathcal C_{p/q}^\pm$ defined in
\eqref{eq:Cplus}--\eqref{eq:Cminus}.
Write $X_{p/q}^{\pm}=\mathcal C_{p/q}^{\pm}\setminus X_{p/q}^{\rm per}$. These are the two classes of rational skew-balanced sequences.

At a rational endpoint these block constraints allow, respectively, one extra or one missing entry equal to one in a resonant block.

\begin{The}[Aubry sets at interior densities]\label{thm:classification-aubry}
Assume \eqref{item:U1}--\eqref{item:U4} and let $\rho=\alpha'_U(c)\in(0,1)$. Then
\[
 \mathscr A(c)=
 \begin{cases}
 X_\rho^{\rm St},&\rho\notin\mathbb Q,\\
 X_{p/q}^{\rm per},&\rho=p/q,\quad c_{U,-}(p/q)<c<c_{U,+}(p/q),\\
 \mathcal C_{p/q}^{+},&\rho=p/q,\quad c=c_{U,+}(p/q),\\
 \mathcal C_{p/q}^{-},&\rho=p/q,\quad c=c_{U,-}(p/q).
 \end{cases}
\]
At either rational endpoint, $\mathscr M(c)=X_{p/q}^{\rm per}\subsetneq\mathscr A(c)$. Every additional configuration is eventually $q$-periodic in both directions, with possibly different translates on its two tails. There are no additional invariant minimizing measures.
\end{The}

\subsubsection{Block defects and closing estimates}
Fix $\rho=\alpha'_U(c)\in(0,1)$ and choose $t\in[0,1]$ so that
\begin{equation}\label{eq:c-t}
 c=(1-t)\beta'_{U,-}(\rho)+t\beta'_{U,+}(\rho).
\end{equation}
For $n\ge2$, set
\begin{equation}\label{eq:tau-rational}
 \tau_n:=(1-t)(\lceil n\rho\rceil-1)+t\lfloor n\rho\rfloor,
 \qquad
 d_n(j;\omega):=f(S_n(j;\omega))-\tau_n S_n(j;\omega)
 +n\rho\tau_n-F(n\rho).
\end{equation}
The quantity $d_n$ is the difference between the pair count and a supporting affine function of $F$ at $n\rho$.

\begin{Pro}[Nonnegative sums for periodic words]\label{app:periodic-defects}
For every $n\ge2$, $0\le d_n(j;\omega)\le3n^2$. For an $N$-periodic word $\omega$,
\begin{equation}\label{app:periodic-action}
 \mathcal S_N\bar A_c(\pi^+\omega)
 =N\bigl[e_U(\omega)-c\rho(\omega)+\alpha_U(c)\bigr]
 =\sum_{n\ge2}\Delta^2U(n)\sum_{j=0}^{N-1}d_n(j;\omega).
\end{equation}
For $n\rho\notin\mathbb Z$, its block defect $d_n(j;\omega)$ vanishes exactly when $S_n(j;\omega)\in\{\lfloor n\rho\rfloor,\lceil n\rho\rceil\}$. For $n\rho=m\in\mathbb Z$, it vanishes exactly when
\begin{equation}\label{eq:dn-resonant-zero}
 \begin{cases}
 S_n(j;\omega)=m,&0<t<1,\\
 S_n(j;\omega)\in\{m-1,m\},&t=0,\\
 S_n(j;\omega)\in\{m,m+1\},&t=1.
 \end{cases}
\end{equation}
\end{Pro}
\begin{proof}
The number $\tau_n$ belongs to the subdifferential of $F$ at $n\rho$, so convexity and $F=f$ on integers give $d_n\ge0$. Since $|\tau_n|\le n$, $0\le S_n\le n$, and $0\le F(n\rho)\le n^2/2$, the displayed definition also gives $d_n\le3n^2$. The derivative formula \eqref{eq:beta-onesided} and convex duality yield
\[
 c=\sum_{n\ge2}n\Delta^2U(n)\tau_n,\qquad
 \alpha_U(c)=\sum_{n\ge2}\Delta^2U(n)
 [n\rho\tau_n-F(n\rho)].
\]
Insert these equalities and the cyclic block-energy identity \eqref{eq:periodic-block-energy}; use $\sum_{j=0}^{N-1}S_n(j;\omega)=nN\rho(\omega)$. All series converge absolutely by Lemma~\ref{lem:U}, proving \eqref{app:periodic-action}. If $n\rho\notin\mathbb Z$ and $k=\lfloor n\rho\rfloor$, direct expansion of $f$ and $F$ gives
\[
 d_n(j;\omega)=\tfrac12(S_n(j;\omega)-k)(S_n(j;\omega)-k-1).
\]
If $n\rho=m\in\mathbb Z$, put $r=S_n(j;\omega)-m$ to obtain $d_n(j;\omega)=\tfrac12r(r+1-2t)$. Their integer zeros are exactly those stated.
\end{proof}

\begin{Lem}[Closing an approximate return]\label{app:closing-return}
Suppose $N\ge m$ and $\mathbf z\in\Sigma_2^+$ satisfies $z_{N+i}=z_i$ for $0\le i<m$. Let $\omega$ be the periodic extension of $z_0\cdots z_{N-1}$. Then
\begin{equation}\label{app:closing-error}
 \left|\mathcal S_N\bar A_c(\mathbf z)
 -\mathcal S_N\bar A_c(\pi^+\omega)\right|
 \le R_m:=\sum_{r>m}rU(r)\longrightarrow0.
\end{equation}
In particular, $h_c(\mathbf x,\mathbf x)\ge0$ for every $\mathbf x$. If this self-barrier is zero, then $d_n(0;\mathbf x)=0$ for every $n\ge2$.
\end{Lem}
\begin{proof}
The sequences $\mathbf z$ and $\pi^+\omega$ agree through index $N+m-1$. Hence their linear terms agree in the first $N$ sites, and so do all pair terms with separation $r\le m$. For $r>m$, at most $\min(N,r-m)$ starting sites can contribute to the difference. Each contribution has absolute value at most $U(r)$. Summing gives \eqref{app:closing-error}.

If both ends of a return are sufficiently close to $\mathbf x$, they have the same first $m$ symbols as $\mathbf x$. For $N\ge m$, the estimate just proved and \eqref{app:periodic-action} imply $\mathcal S_N\bar A_c(\mathbf z)\ge-R_m$. Taking the limits in \eqref{eq:Peierls-barrier} and then $m\to\infty$ shows that every self-barrier is nonnegative. For fixed $n\le m$, the periodic closure also satisfies $d_n(0;\omega)=d_n(0;\mathbf x)$. Keeping this single nonnegative term in \eqref{app:periodic-action} yields
\[
 \mathcal S_N\bar A_c(\mathbf z)
 \ge\Delta^2U(n)d_n(0;\mathbf x)-R_m.
\]
The same limiting argument gives $h_c(\mathbf x,\mathbf x)\ge\Delta^2U(n)d_n(0;\mathbf x)$. Strict positivity of $\Delta^2U(n)$ proves the final assertion.
\end{proof}

\begin{Lem}[Recurrent sequences with zero block defects]\label{app:recurrent-return}
Suppose $d_n(j;\omega)=0$ for all $n\ge2$ and $j\in\mathbb Z$. If $\pi^+\omega$ is recurrent, then $h_c(\pi^+\omega,\pi^+\omega)=0$.
\end{Lem}
\begin{proof}
For each $m$, recurrence gives $N_m\ge m$ tending to infinity with
$\omega_{N_m+i}=\omega_i$ for $0\le i<m$.
Let $\omega^{(m)}$ be the $N_m$-periodic closure of its initial word.
Every block of length $n\le m$ agrees with a block of $\omega$, including those crossing the cut, and has zero defect.
For $m<n\le N_m$, at most $n-1$ blocks cross the cut; for $n>N_m$, there are only $N_m$ starting sites.
Thus Proposition~\ref{app:periodic-defects} gives
\[
 0\le\mathcal S_{N_m}\bar A_c(\pi^+\omega^{(m)})
 \le3\sum_{n>m}n^2\min(n,N_m)\Delta^2U(n)
 \le3\sum_{n>m}n^3\Delta^2U(n)\to 0\qquad (m\to\infty)
\]
by \eqref{app:first-moment}.
Both ends of these periodic returns converge to $\pi^+\omega$, so the barrier is at most zero.
Lemma~\ref{app:closing-return} supplies the reverse inequality.
\end{proof}

\begin{Pro}[Classification of the block constraints]\label{app:block-classification}
The two-sided sequences satisfying $d_n(j;\omega)=0$ for all $n\ge2$ and $j\in\mathbb Z$ are $X_\rho^{\rm St}$ at irrational density, $X_{p/q}^{\rm per}$ at rational density with $0<t<1$, $\mathcal C_{p/q}^+$ at $t=1$, and $\mathcal C_{p/q}^-$ at $t=0$. Every sequence in $X_{p/q}^{+}\cup X_{p/q}^{-}$ has $q$-periodic tails.
\end{Pro}
\begin{proof}
At irrational density, all indices are nonresonant; the zero-defect  condition ($d_n=0$)  is exactly \eqref{eq:balanced-block-condition}, which characterizes $X_\rho^{\rm St}$. At rational density with $0<t<1$, the same condition follows from \eqref{eq:dn-resonant-zero}, and Lemma~\ref{lem:rational-balanced} identifies $X_{p/q}^{\rm per}$. At $t=1$ or $t=0$, the zero-defect conditions are exactly \eqref{eq:Cplus} or \eqref{eq:Cminus}  by Proposition \ref{app:periodic-defects}.

For $\omega\in\mathcal C_{p/q}^+$, put $e_j=S_q(j;\omega)-p\in\{0,1\}$. For every $\ell\ge1$,
\[
 \sum_{r=0}^{\ell-1}e_{j+rq}
 =S_{\ell q}(j;\omega)-\ell p\in\{0,1\}.
\]
Consequently each residue class modulo $q$ contains at most one index with $e_j=1$. The sequence $\{e_j\}_{j\in\mathbb Z}$ has finite support, and $e_{j+1}-e_j=\omega_{j+q}-\omega_j$ shows that $\omega$ is eventually $q$-periodic in either direction. Its periodic tails satisfy the rational block condition and thus belong to $X_{p/q}^{\rm per}$. For $\mathcal C_{p/q}^-$ the same argument applies to $e_j=p-S_q(j;\omega)$.
\end{proof}

\begin{Lem}[Periodic approximation of rational skew sequences]\label{lem:skew-approximation}
Let $p/q\in(0,1)$ be reduced. For every $\omega\in\mathcal C_{p/q}^{\pm}$, there are $q_k=M_kq\to\infty$ and periodic mechanical sequences $\omega^{(k)}$ of density
\[
 \rho_k=\frac{p}{q}\pm\frac1{q_k}
\]
that converge locally to $\omega$, after suitable shifts.
\end{Lem}
\begin{proof}
Write $\rho=p/q$ and fix a length-$L$ block $\omega_a\cdots\omega_{a+L-1}$ of $\omega\in\mathcal C_{p/q}^+$.
Set $P_0=0$, $P_j=\sum_{i=0}^{j-1}\omega_{a+i}$ and $T_j(\varrho)=P_j-j\varrho$.
This block is mechanical at density $\varrho$ precisely when some $\vartheta\in[0,1)$ satisfies
\begin{equation}\label{eq:theta-intersection}
 P_j\le j\varrho+\vartheta<P_j+1\qquad(0\le j\le L).
\end{equation}
The corresponding half-open intervals intersect precisely when
\begin{equation}\label{eq:Tj}
 \max_{0\le j\le L}T_j(\varrho)
 <\min_{0\le j\le L}T_j(\varrho)+1;
\end{equation}
the $j=0$ interval enforces $0\le\vartheta<1$.
For $r<s$, the upper block constraints give
\begin{equation}\label{eq:T-difference-at-h}
 T_s(\rho)-T_r(\rho)
 \in\{-\{(s-r)\rho\},1-\{(s-r)\rho\}\}.
\end{equation}
Nonresonant differences lie in $(-1,1)$; resonant differences are $0$ or $1$.
Replacing $\rho$ by $\rho+\varepsilon$ decreases each difference by $(s-r)\varepsilon$.
There are finitely many pairs, so all differences lie in $(-1,1)$ for every sufficiently small $\varepsilon>0$.
This proves \eqref{eq:Tj}.
Choose $\varepsilon=1/(Mq)$ with $M$ large; the resulting mechanical sequence is $Mq$-periodic and contains the prescribed block.
Applying this to $\omega|_{[-k,k]}$, with $M_k\to\infty$, and shifting the block to its original position proves local convergence.

For $\mathcal C_{p/q}^-$, resonant differences are $-1$ or $0$.
Replacing $\rho$ by $\rho-\varepsilon$ increases them by $(s-r)\varepsilon$; the same finite-intersection argument and diagonal choice prove the lower approximation.
\end{proof}

\subsubsection{Classification of the Aubry sets}
\begin{proof}[Proof of Theorem~\ref{thm:classification-aubry}]
If $\omega\in\mathscr A(c)$, apply Lemma~\ref{app:closing-return} to every forward tail. It gives $d_n(j;\omega)=0$ for all $n,j$. Proposition~\ref{app:block-classification} therefore proves the required inclusion into the sets displayed in the theorem.

Every point of a Sturmian subshift is recurrent, as is every point of a periodic orbit. Their block defects vanish by \eqref{eq:block-telescope}, the block characterization, and Proposition~\ref{app:periodic-defects}. Lemma~\ref{app:recurrent-return}, applied after every shift, gives the reverse inclusion for irrational densities and for the periodic orbit at every point of a rational plateau.

At the upper rational endpoint, let $\omega\in\mathcal C_{p/q}^+$ and choose the periodic approximants of Lemma~\ref{lem:skew-approximation}. With $\rho=p/q$ and $c=c_{U,+}(\rho)=\beta'_{U,+}(\rho)$, their periods and minimum energies give
\begin{align}\label{eq:endpoint-periodic-action}
\mathcal S_{q_k}\bar A_c(\pi^+\omega^{(k)})
 =q_k[\beta_U(\rho_k)-\beta_U(\rho)-c(\rho_k-\rho)]=\frac{\beta_U(\rho+1/q_k)-\beta_U(\rho)}{1/q_k}-c
\to0
\end{align}
as $k\to\infty$. Since the approximants are periodic and converge locally, both ends approach $\pi^+\omega$ along periods tending to infinity.
They are admissible in \eqref{eq:Peierls-barrier}, so $h_c(\pi^+\omega,\pi^+\omega)\le0$.
Lemma~\ref{app:closing-return} gives the reverse inequality. Applying the same construction to every shift proves $\omega\in\mathscr A(c)$. At the lower endpoint the corresponding action is
\[
 \frac{\beta_U(\rho-1/q_k)-\beta_U(\rho)}{1/q_k}+c
 \to-\beta'_{U,-}(\rho)+c=0,\qquad (k\to\infty)
\]
which proves the remaining inclusion.

To see that the endpoint inclusions are strict, take coordinatewise limits
\[
 \omega^{\pm}=\lim_{\varepsilon\to 0^+}
 \omega^{\rho\pm\varepsilon,0},\quad  S_n(j;\omega^\pm)
 =\lim_{\varepsilon\to 0^+}
 S_n(j;\omega^{\rho\pm\varepsilon,0}).
\]
The floor formulas \eqref{eq:lower-mechanical}--\eqref{eq:per} ensure existence and gives $\omega^{\pm}\in\mathcal C_{p/q}^{\pm}$, while when $q\varepsilon\in(0,1)$
\begin{align*}
     S_q(-q;\omega^{\rho+\varepsilon,0})=\lfloor0\rfloor-\lfloor-p-q\varepsilon\rfloor=p+1,\qquad
 S_q(0;\omega^{\rho-\varepsilon,0})=\lfloor p-q\varepsilon\rfloor-\lfloor0\rfloor=p-1.
\end{align*}
by passing to the coordinatewise limits. Thus these sequences are not in $X_{p/q}^{\rm per}$. The periodic-tail assertion follows from Proposition~\ref{app:block-classification}, and the uniqueness of invariant minimizing measures follows from Theorem~\ref{app:general-variational}.
\end{proof}

\section{Analytic estimates for the complete Coulomb coefficients}
\label{co:analysis}

Starting from the exact matrix elements in
\cite[Supplementary Material, \S I, Eq.~(2)]{RotondoEtAl2016},
this appendix proves only the three estimates in
Lemma~\ref{co:bounds}. The thin-torus approximations in
\cite[Supplementary Material, \S II]{RotondoEtAl2016} do not
supply these finite-size error bounds. We therefore retain all
periodic sums and estimate their remainders explicitly.

\subsection{Exact coefficients}

For $m\ne0$ set
\begin{equation}\label{co:Am}
 A_m(r)=\sum_{n\in\Z}w_{mn}e^{2\pi i nr/L},\qquad a_m=A_m(0),
 \qquad r\in\Z/L\Z.
\end{equation}
These coefficients are real and nonnegative, as established below.
Define
\[
 \mathcal I_L=\{(m,d):m\in\Z\setminus L\Z,\quad
       d\in\Z/L\Z,\quad d\not\equiv0,m,-m\pmod L\}.
\]

\begin{Pro}[Complete Coulomb coefficients]\label{co:coefficients}
For every $L\ge2$ and $\kappa>0$, the coefficient of the full
diagonal $D^{\mathrm C}$ satisfies $W(L-r)=W(r)$ and
\begin{align}
 K^{\mathrm C}_{L,\kappa}(r)
 &=C_\kappa\sum_{n\ge1}
       \frac{e^{-2\pi^2n^2/(\kappa^2L^2)}}n
                      \cos\frac{2\pi nr}{L},\label{co:K}\\
 W(r)&=K^{\mathrm C}_{L,\kappa}(r)
       -2\sum_{m\equiv r\ (L)}a_m
       +2\sum_{m\in L\Z\setminus\{0\}}A_m(r),
       \qquad 1\le r<L.\label{co:W}
\end{align}
The off-diagonal part is exactly
\begin{equation}\label{co:hopping}
 T^{\mathrm C}=\sum_{i=0}^{L-1}\sum_{(m,d)\in\mathcal I_L}
 A_m(d)c_i^*c_{[i+d+m]}^*c_{[i+d]}c_{[i+m]}.
\end{equation}
In particular, every displayed nonzero monomial has four distinct
sites; all integer lifts of $m$ remain in the sum.

\end{Pro}

\begin{proof}
Moving exterior multiplications to the left of contractions, using
$c_jc_i^*=\delta_{ij}I-c_i^*c_j$, gives the exact identity
\begin{equation}\label{co:quartic}
 \mathcal R_{mn}^*\mathcal R_{mn}-\Num
 =\sum_{i,j=0}^{L-1}e^{2\pi i n(j-i)/L}
                  c_i^*c_{[j+m]}^*c_jc_{[i+m]}.
\end{equation}
The terms with $m=0$ give \eqref{co:K}. For
$m\not\equiv0\pmod L$, the terms $j=i$ are
$-a_m n_i n_{[i+m]}$; collecting unordered pairs gives the exchange
term in \eqref{co:W}. For $m\in L\Z\setminus\{0\}$ every term
is diagonal and gives the last, winding contribution in
\eqref{co:W}. The remaining terms vanish when $j-i=\pm m$ modulo
$L$, since then two creation or two annihilation indices coincide.
All other terms have four distinct sites and give
\eqref{co:hopping}.

For later use, put $h=\sum_{(m,d)\in\mathcal I_L}A_m(d)$.
Positivity of these coefficients is verified below.
On a binary word of degree $N$, a monomial in
\eqref{co:hopping} acts only if $i$ is empty and $i+m$ is occupied.
For fixed $(m,d)$ there are at most $\min\{N,L-N\}$ such anchors.
The absolute column sum, and by self-adjointness the absolute
row sum, is therefore at most $h\min\{N,L-N\}$. Hence
\[
 \bigl\|T^{\mathrm C}|_{\cH_{L,N}}\bigr\|
 \le h\min\{N,L-N\},\qquad
 \|T^{\mathrm C}\|\le Lh/2.
\]
\end{proof}

We abbreviate $K^{\mathrm C}_{L,\kappa}$ to $K$.

\subsection{The Gaussian-smoothed logarithmic kernel}

\begin{Lem}\label{co:logarithmic-curvature}
For $L\ge2$ and $\kappa\ge3$,
\[
 K''(x)\ge\frac{C_\kappa}{x^2}\quad(1\le x\le L-1).
\]
\end{Lem}
\begin{proof}
Let $Z$ be a standard Gaussian, and put
$f_\sigma(x)=\mathbb E[-\log|x+\sigma Z|]$.
Gaussian convolution and integration by parts give
$\sigma^2f_\sigma''+xf_\sigma'=-1$ and $f_\sigma'(0)=0$.
With $z=x/(\sqrt2\sigma)$ and
$\mathcal D(z)=e^{-z^2}\int_0^ze^{t^2}\,dt$, this gives
\[
 f_\sigma''(x)=\sigma^{-2}(2z\mathcal D(z)-1).
\]
For $z\ge2$, define
$b(z)=2z\mathcal D(z)-1-(2z^2)^{-1}$.
The change of variable $t=2-u$ yields
\[
 \mathcal D(2)=\int_0^2e^{-4u+u^2}\,du
 \ge\int_0^2e^{-4u}(1+u^2+u^4/2)\,du
 =\frac{75-1283e^{-8}}{256}>\frac9{32}.
\]
Thus $b(2)>0$. At any zero of $b$, the identity
$\mathcal D'=1-2z\mathcal D$ gives $b'=3/(2z^3)>0$.
Consequently $b$ cannot cross from positive to nonpositive, and
\begin{equation}\label{co:gaussian-curvature}
 f_\sigma''(x)\ge x^{-2}\qquad
                (|x|\ge2\sqrt2\sigma).
\end{equation}
The same differential equation, or Gaussian Fourier inversion, gives
\[
 f_\sigma''(x)=-\int_0^\infty
       t e^{-\sigma^2t^2/2}\cos(tx)\,dt,
 \qquad
 \widehat{f_\sigma''}(\xi)=-\pi|\xi|
                                e^{-\sigma^2\xi^2/2},
\]
where $\widehat f(\xi)=\int_\R f(x)e^{-i\xi x}\,dx$.
Two integrations by parts in the integral show
$f_\sigma''(x)=O_\sigma(x^{-2})$. Periodization and comparison of
Fourier coefficients therefore give the absolutely convergent identity
\begin{equation}\label{co:periodized-curvature}
 K''(x)=C_\kappa\sum_{s\in\Z}f_{1/\kappa}''(x+sL).
\end{equation}
For $1\le x\le L-1$, every $|x+sL|\ge1>2\sqrt2/\kappa$.
All terms are positive; retaining $s=0$ proves the bound.
\end{proof}

\subsection{The remaining diagonal coefficient terms}

For $m\ne0$, the Gaussian integral for the inverse square root and
Poisson summation give
\begin{align}
 A_m(r)&=\sum_{s\in\Z}\alpha_m(r+sL),\label{co:poisson}\\
 \alpha_m(t)&=\frac{g\kappa}{4\pi}
   \int_{1/2}^\infty
   \frac{\exp\{-\kappa^2[m^2u+t^2/(4u)]\}}
             {\sqrt{u(u-1/2)}}\,du\ge0.
   \label{co:positive-integral}
\end{align}
For example, this follows by inserting
\[
 \frac1{\sqrt{b^2+x^2}}
 =\frac1{\sqrt\pi}\int_0^\infty
       t^{-1/2}e^{-t(b^2+x^2)}\,dt
\]
with $b=\kappa|m|$ and $x=2\pi n/(\kappa L)$, combining the
Gaussian lowest eigenspace factor, and setting $u=t+1/2$ before summing in $n$.
All exchanges are justified by the positive integrals and the
convergent Gaussian sums. In particular $A_m(r)\ge0$, while its
Fourier series gives $A_m(r)\le a_m$.

Assume $\kappa\ge4,L\ge6$, and set $a_*=\kappa^2/2$.
A Gaussian sum bounded by one plus its integral gives
\begin{align*}
 a_m&\le\frac{g}{2L\kappa|m|}e^{-a_*m^2}
          \sum_{n\in\Z}e^{-2\pi^2n^2/(\kappa^2L^2)}\\
 &\le g\left(\frac1{2L\kappa}+
                       \frac1{2\sqrt{2\pi}}\right)
                  \frac{e^{-a_*m^2}}{|m|}
 \le\frac g4\frac{e^{-a_*m^2}}{|m|}.
\end{align*}
Consequently, for every integer $M\ge1$,
\begin{equation}\label{co:exchange-tail}
 S_M:=\sum_{|m|\ge M}a_m
 \le\frac{g}{2M}
       \frac{e^{-a_*M^2}}{1-e^{-a_*(2M+1)}}
 \le\frac gM e^{-a_*M^2}.
\end{equation}
Write $W-K=X+Y$, with $X$ and $Y$ the exchange and winding terms
in \eqref{co:W}. If $2\le r\le J_L$, every integer in one of the
residue classes $r-1,r,r+1$ has absolute value at least $r-1$.
It follows that
\[
 |\Delta^2X(r)|\le8S_{r-1},\qquad
 r^2|\Delta^2X(r)|\le32g e^{-a_*}.
\]
For the second bound, use that
$(s+1)^2s^{-1}e^{-a_*s^2}$ decreases for $s\ge1,a_*\ge8$.
Similarly,
\[
 |\Delta^2Y(r)|\le8S_L,\qquad
 r^2|\Delta^2Y(r)|\le2gL e^{-a_*L^2}
                         \le12g e^{-36a_*}.
\]
Lemma~\ref{co:logarithmic-curvature} and triangular averaging give
\[
 \Delta^2K(r)=\int_{-1}^1(1-|t|)K''(r+t)\,dt
 \ge C_\kappa\int_{-1}^1\frac{1-|t|}{(r+t)^2}\,dt
 \ge\frac{C_\kappa}{r^2}.
\]
Combining these inequalities proves the first estimate of
Lemma~\ref{co:bounds}, with the stated $\xi_\kappa$. Finally $\xi_4<0.00844$, and
$\xi_\kappa$ decreases for $\kappa\ge4$.

\subsection{Summing the complete off-diagonal part}

We show $h\le2g\kappa e^{-2\kappa^2}$ for
$\kappa\ge3$ and $L\ge2$.
Every lift $(m,r)$ occurring in
$A_m(d)=\sum_s\alpha_m(d+sL)$ for $(m,d)\in\mathcal I_L$
has $m,r\ne0$ and $|mr|\ge2$.
Indeed $|mr|=1$ would imply $m=\pm1$, $r=\pm1$, and hence
$d\equiv m$ or $-m$, which is excluded.

Put $z=\kappa^2|mr|$. The inequality
$m^2u+r^2/(4u)\ge|mr|$ bounds the integral in
\eqref{co:positive-integral} over $[1/2,1]$ by
$2\operatorname{arcosh}(\sqrt2)e^{-z}$.
On $[1,\infty)$ use
$[u(u-1/2)]^{-1/2}\le\sqrt2/u$ and extend the resulting integral
to $(0,\infty)$. Its value is $2\sqrt2 K_0(z)$, where
\[
 K_0(z)=\int_0^\infty e^{-z\cosh t}\,dt
       \le\sqrt{\frac\pi{2z}}e^{-z}.
\]
Here the equality follows from $u=(|r|/(2|m|))e^t$ and symmetry
in $t$, and the inequality follows from $\cosh t\ge1+t^2/2$.
Thus
\begin{equation}\label{co:alpha-hopping-bound}
 \alpha_m(r)\le
 g\left(\frac{\kappa\operatorname{arcosh}(\sqrt2)}{2\pi}
        +\frac1{2\sqrt{\pi|mr|}}\right)e^{-\kappa^2|mr|}.
\end{equation}
Set $x=e^{-\kappa^2}$ and
\[
 B_\kappa=\frac{\kappa\operatorname{arcosh}(\sqrt2)}{2\pi}
                  +\frac1{2\sqrt{2\pi}}.
\]
Dropping the residue restrictions only enlarges the positive sum.
If $\tau(k)$ denotes the number of positive divisors of $k$, then
\begin{align*}
 h&\le\sum_{\substack{m,r\in\Z\setminus\{0\}\\|mr|\ge2}}
                       \alpha_m(r)
 \le4gB_\kappa\sum_{k\ge2}\tau(k)x^k\le4gB_\kappa\frac{x^2(2-x)}{(1-x)^2}
 \le2g\kappa e^{-2\kappa^2}.
\end{align*}
The second line uses $\tau(k)\le k$. For the last inequality,
$B_\kappa/\kappa<0.22$ and $(2-x)/(1-x)^2<2.01$ for
$\kappa\ge3$, and $4(0.22)(2.01)<2$.
The bound $\|T^{\mathrm C}\|\le Lh/2$ proved above now
gives the second estimate of Lemma~\ref{co:bounds}.

\subsection{Compression to the balanced subspace}

\begin{proof}[Proof of the third estimate in Lemma~\ref{co:bounds}]
For $N=0,L$, the off-diagonal operator vanishes and the assertion
is immediate. Otherwise write $N/L=p/q$ in lowest terms.
On the translated balanced basis vectors $\Psi_a$, the eigenvalues
of $\mathsf C$ are $z_0e^{2\pi iap/q}$, $0\le a<q$,
and are distinct. An operator on this subspace commuting with
$\mathsf S$ and $\mathsf C$ is therefore scalar. Both
$P_NT^{\mathrm C}P_N$ and $P_N(T^{\mathrm C})^2P_N$ have this
property. Since $T^{\mathrm C}$ has zero diagonal, the first
compression is zero, while the second equals
$\|T^{\mathrm C}\Psi_a\|^2P_N$.

The four-site expansion \eqref{co:hopping} and the column estimate
in the proof of Proposition~\ref{co:coefficients} give
\[
 \sum_\tau|\langle\tau,T^{\mathrm C}\Psi_a\rangle|
 \le\min\{N,L-N\}h.
\]
For a fixed resulting configuration $\tau$, the index $i$ of the leftmost exterior multiplication must be one
of the two indices added to the original word. For each such index, the
sum of all contributing coefficients, including all periodic
lifts, is at most $h$.  Hence every matrix element in this column
has absolute value at most $2h$.  Multiplying these two bounds
gives \eqref{co:compressions}, using
$h\le2g\kappa e^{-2\kappa^2}$.
\end{proof}

\end{document}